\documentclass[12pt]{amsart}
\usepackage{currvita, amssymb, amsmath, amsthm, 
graphicx, subfigure, pifont, wrapfig, pinlabel, amscd,
latexsym ,exscale, enumerate, amsfonts, times, color, hyphenat, pinlabel, 
mathtools, a4wide, fullpage, array, bbm, url}
\usepackage{float}
\usepackage[colorlinks=false, pdfborder={0.9 0.9 1}]{hyperref}
\usepackage{stmaryrd}
\usepackage{multirow}
\usepackage[all]{xy}

\newcommand{\bn}[1]{\llbracket #1 \rrbracket}

\DeclareMathOperator{\Mor}{Mor}

\newcommand{\bN}{\mathbb{N}}
\newcommand{\bZ}{\mathbb{Z}}
\newcommand{\bQ}{\mathbb{Q}}
\newcommand{\bR}{\mathbb{R}}

\DeclareMathOperator{\Tor}{Tor}
\DeclareMathOperator{\Kh}{Kh}

\newcommand{\refequal}[1]{\xy {\ar@{=}^{#1}
(-1,0)*{};(1,0)*{}};
\endxy}

\DeclareMathOperator{\Ob}{Ob}
\DeclareMathOperator{\mat}{\textbf{Mat}}

\DeclareMathOperator{\Cb}{\textbf{Cb}}
\DeclareMathOperator{\cob}{\textbf{Cob}^2}

\DeclareMathOperator{\ucob}{u\textbf{Cob}^2}
\DeclareMathOperator{\ukob}{\textbf{Kob}_b(\emptyset)}
\DeclareMathOperator{\ukobk}{\textbf{Kob}_b(k)}

\DeclareMathOperator{\RMOD}{R-\textbf{Mod}}

\DeclareMathOperator{\KAR}{\textbf{Kar}}

\DeclareMathOperator{\Kom}{\textbf{Kom}}

\newtheorem{prop}{Proposition}[section]
\newtheorem{thm}[prop]{Theorem}
\newtheorem{lem}[prop]{Lemma}
\newtheorem{cor}[prop]{Corollary}

\theoremstyle{remark}
\newtheorem{rem}[prop]{Remark}

\newtheorem*{nota}{Notation}

\theoremstyle{remark}
\newtheorem{ex}[prop]{\textbf{Example}}

\theoremstyle{remark}

\theoremstyle{definition}
\newtheorem{defn}[prop]{Definition}

\numberwithin{equation}{section}

\newcommand{\jpg}[2]{{\hspace{-3pt}\begin{array}{c}%
  \raisebox{-2.5pt}{\includegraphics[width=#1]{res/figs/topcom/#2.eps}}%
\end{array}\hspace{-3pt}}}

\newcommand{\backoverslash}{\xy
 (0,0)*{\includegraphics[height=.02\textheight]{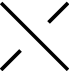}};
 \endxy}
\newcommand{\slashoverback}{\xy
 (0,0)*{\includegraphics[height=.02\textheight]{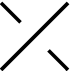}};
 \endxy}
\newcommand{\smoothing}{\xy
 (0,0)*{\includegraphics[height=.02\textheight]{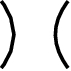}};
 \endxy}
\newcommand{\hsmoothing}{\xy
 (0,0)*{\includegraphics[height=.02\textheight]{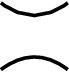}};
 \endxy}
\newcommand{\virtual}{\xy
 (0,0)*{\includegraphics[height=.02\textheight]{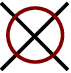}};
 \endxy}
 \newcommand{\overcrossing}{\xy
 (0,0)*{\includegraphics[height=.02\textheight]{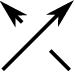}};
 \endxy}
\newcommand{\undercrossing}{\xy
 (0,0)*{\includegraphics[height=.02\textheight]{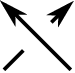}};
 \endxy}
 
\newcommand{\virtualoriented}{\xy
 (0,0)*{\includegraphics[height=.02\textheight]{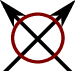}};
 \endxy}

\newcommand{\uu}{\xy
 (0,0)*{\includegraphics[height=.02\textheight]{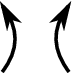}};
 \endxy}
\newcommand{\dd}{\xy
 (0,0)*{\includegraphics[height=.02\textheight]{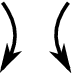}};
 \endxy}
\newcommand{\du}{\xy
 (0,0)*{\includegraphics[height=.02\textheight]{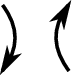}};
 \endxy}
\newcommand{\ud}{\xy
 (0,0)*{\includegraphics[height=.02\textheight]{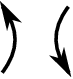}};
 \endxy}
 \newcommand{\lel}{\xy
 (0,0)*{\includegraphics[height=.02\textheight, angle=90]{res/figs/orsmoothings/uu.eps}};
 \endxy}
\newcommand{\rir}{\xy
 (0,0)*{\includegraphics[height=.02\textheight, angle=90]{res/figs/orsmoothings/dd.eps}};
 \endxy}
\newcommand{\ril}{\xy
 (0,0)*{\includegraphics[height=.02\textheight, angle=90]{res/figs/orsmoothings/du.eps}};
 \endxy}
\newcommand{\ler}{\xy
 (0,0)*{\includegraphics[height=.02\textheight, angle=90]{res/figs/orsmoothings/ud.eps}};
 \endxy}
 
\newcommand{\uup}{\xy
 (0,0)*{\includegraphics[height=.025\textheight]{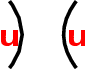}};
 \endxy}
\newcommand{\uupp}{\xy
 (0,0)*{\includegraphics[width=.025\textheight]{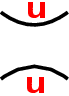}};
 \endxy}
\newcommand{\udp}{\xy
 (0,0)*{\includegraphics[height=.025\textheight]{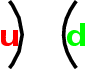}};
 \endxy}
\newcommand{\udpp}{\xy
 (0,0)*{\includegraphics[width=.025\textheight]{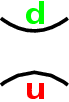}};
 \endxy}
\newcommand{\dup}{\xy
 (0,0)*{\includegraphics[height=.025\textheight]{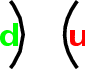}};
 \endxy}
\newcommand{\dupp}{\xy
 (0,0)*{\includegraphics[width=.025\textheight]{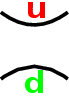}};
 \endxy}
\newcommand{\ddp}{\xy
 (0,0)*{\includegraphics[height=.025\textheight]{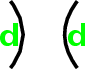}};
 \endxy}
\newcommand{\ddpp}{\xy
 (0,0)*{\includegraphics[width=.025\textheight]{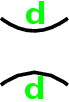}};
 \endxy}
\newcommand{\none}{\xy
 (0,0)*{\includegraphics[height=.03\textheight]{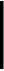}};
 \endxy}
\newcommand{\down}{\xy
 (0,0)*{\includegraphics[height=.03\textheight]{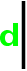}};
 \endxy}
\newcommand{\up}{\xy
 (0,0)*{\includegraphics[height=.03\textheight]{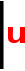}};
 \endxy}
\newcommand{\smoothinga}{\xy
 (0,0)*{\includegraphics[height=.025\textheight]{res/figs/intro/smoothing.eps}};
 \endxy}
\newcommand{\hsmoothinga}{\xy
 (0,0)*{\includegraphics[height=.025\textheight]{res/figs/intro/hsmoothing.eps}};
 \endxy} 
 
 \newcounter{In}

\def\In{\stepcounter{In}\vskip 6pt \par\noindent\scriptsize
  In[\theIn]$:=$\tt
  \catcode`\#=12 \catcode`\&=12 \catcode`\~=12 \catcode`\%=12
  \catcode`\$=12 \catcode`\^=12 \catcode`\_=12
  \obeylines\obeyspaces}
\def\Out{\vskip 6pt \par\noindent\scriptsize\parindent=12pt
  Out[\theIn]$=$\tt
  \catcode`\#=12 \catcode`\&=12 \catcode`\~=12 \catcode`\%=12
  \catcode`\$=12 \catcode`\^=12 \catcode`\_=12 \catcode`\{=12
  \catcode`\}=12 \obeyspaces
}

\title{Virtual Khovanov homology using cobordisms}
\author{D.~Tubbenhauer}
\thanks{The author was supported by the German Research 
Foundation (Deutsche Forschungsgemeinschaft (DFG)) 
through the Institutional Strategy of the University of G\"{o}ttingen and by the Graduiertenkolleg 1493 ``Mathematische Strukturen in der modernen Quantenphysik''.}
\date{Last compiled \today; last edited \today}

\begin{document}
\begin{abstract}
We extend Bar-Natan's cobordism based categorification of the Jones polynomial to virtual links. Our topological complex allows a direct extension of the classical Khovanov complex ($h=t=0$), the variant of Lee ($h=0,t=1$) and other classical link homologies. We show that our construction allows, over rings of characteristic two, extensions with no classical analogon, e.g. Bar-Natan's $\bZ/2$-link homology can be extended in two non-equivalent ways.

Our construction is computable in the sense that one can write a computer program to perform calculations, e.g. we have written a MATHEMATICA based program.

Moreover, we give a classification of all unoriented TQFTs which can be used to define virtual link homologies from our topological construction. Furthermore, we prove that our extension is combinatorial and has semi-local properties. We use the semi-local properties to prove an application, i.e. we give a discussion of Lee's degeneration of virtual homology.
\end{abstract}

\maketitle

\paragraph*{Acknowledgements}
The author thanks Vassily Manturov for corrections and many helpful comments. Moreover, I wish to thank Aaron Kaestner and Louis Kauffman for suggestions and helpful comments, hopefully helping me to write a faster computer program for calculations in future work. I also thank Thomas Schick, Marco Mackaay and A referee whose further suggestions greatly improved the presentation. Moreover, special thanks to A referee for spotting a crucial typo in one of the relations.

It seems that the remaining typos and mistakes are the main contribution of the author.
\tableofcontents
\section{Introduction}\label{sec-introa}
We consider \textit{virtual link diagrams} $L_D$ in this paper, i.e. planar graphs of valency four where every vertex is either an overcrossing $\slashoverback$, an undercrossing $\backoverslash$ or a virtual crossing $\virtual$, which is marked with a circle. We also allow circles, i.e. closed edges without any vertices.

We call the crossings $\slashoverback$ and $\backoverslash$ \textit{classical crossings} or just \textit{crossings}. For a virtual link diagram $L_D$ we define the \textit{mirror image} $\overline L_D$ of $L_D$ by switching all classical crossings from an overcrossing to an undercrossing and vice versa.

\begin{figure}[ht]
      \centerline{\includegraphics[scale=0.55]{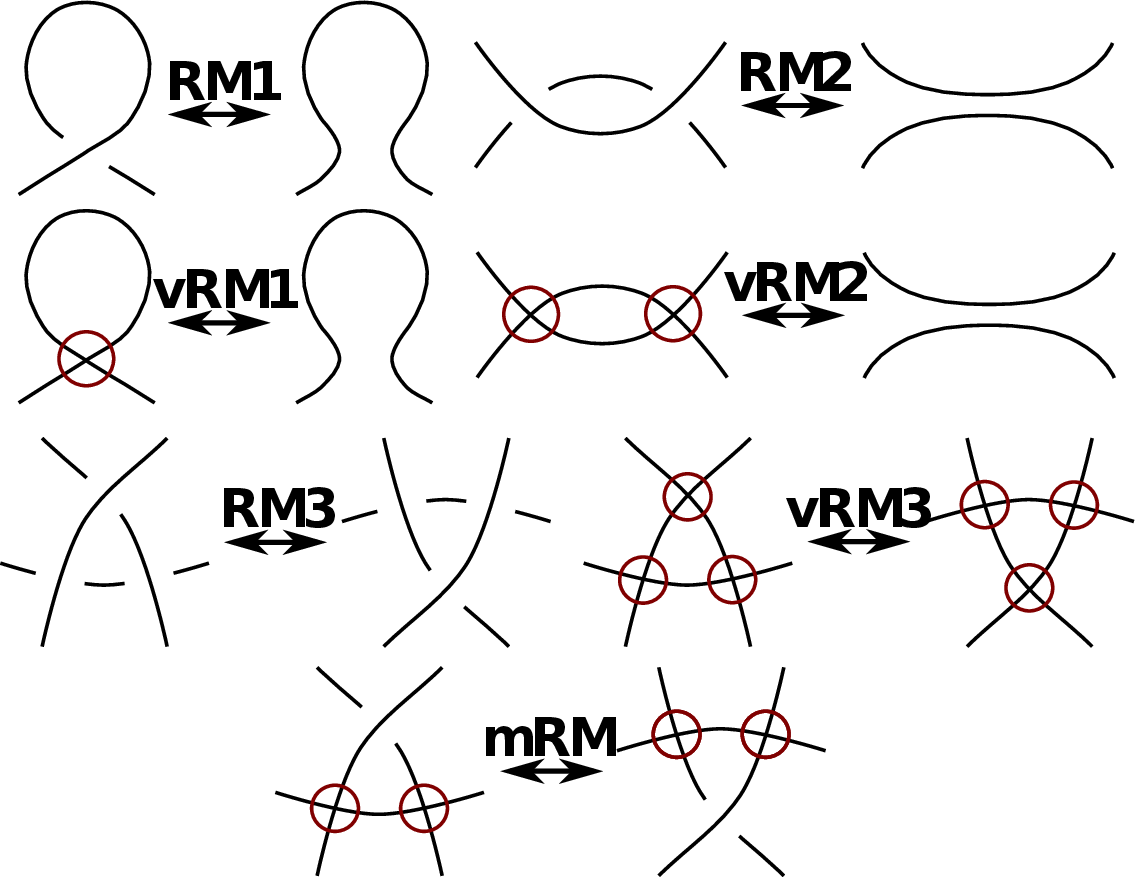}}
  \caption{The generalised Reidemeister moves are the moves pictured plus mirror images.}
  \label{figureintroa-1}
\end{figure}

A \textit{virtual link} $L$ is an equivalence class of virtual link diagrams modulo planar isotopies and \textit{generalised Reidemeister moves}, see Fig.~\ref{figureintroa-1}.

We call the moves RM1, RM2 and RM3 the \textit{classical Reidemeister moves}, the moves vRM1, vRM2 and vRM3 the \textit{virtual Reidemeister moves} and the move mRM the \textit{mixed Reidemeister move}. We call a virtual link diagram $L_D$ \textit{classical} if all crossings of $L_D$ are classical crossings. Furthermore, we say a that virtual link $L$ is \textit{classical}, if the set $L$ contains a classical link diagram.

The notions of an \textit{oriented} virtual link diagram and of an \textit{oriented} virtual link are defined analogously. The latter modulo isotopies and \textit{oriented} generalised Reidemeister moves. Note that an \textit{oriented virtual link diagram} is a diagram together with a choice of an orientation of the diagram such that every crossing is of the form $\overcrossing$, $\undercrossing$ or $\virtualoriented$. Furthermore, we use the shorthand notations c- and v- for everything that starts with classical or virtual, e.g. c-knot means classical knot and v-crossing means virtual crossing. 
\vspace*{0.5cm}

Virtual links are an essential part of modern knot theory and were proposed by Kauffman in~\cite{ka3}. They arise from the study of links which are embedded in a thickened, orientable surface $\Sigma_g$ of genus $g\geq 0$. These links were studied by Jaeger, Kauffman and Saleur in~\cite{jks}. Note that for c-links the surface is $\Sigma_g=S^2$, i.e. v-links are a generalisation of c-links and they should for example have analogous ``applications'' in quantum physics.

\begin{figure}[ht]
     \centerline{\includegraphics[scale=0.5]{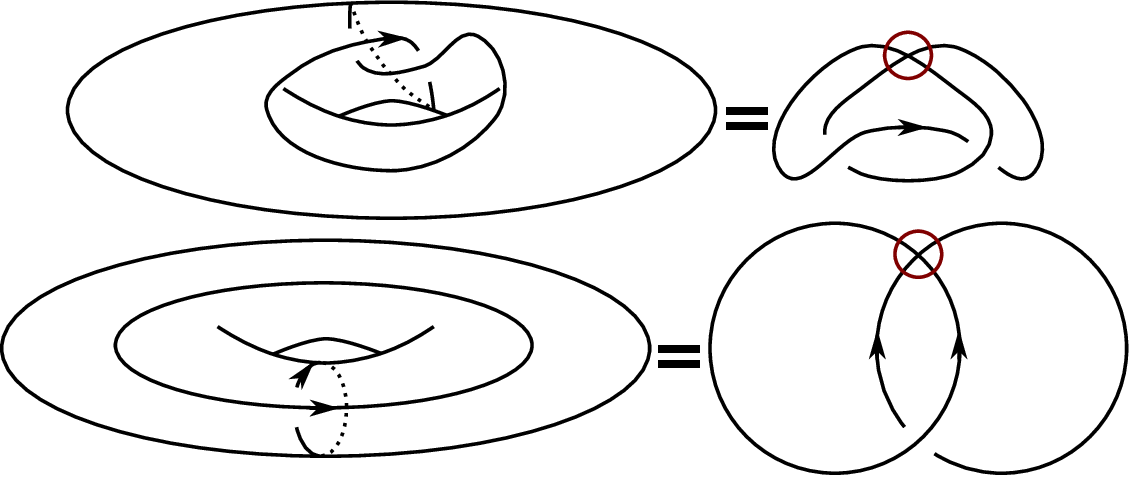}}
  \caption{Two knot diagrams on a torus. The first virtual knot is called the \textit{virtual trefoil}.}
  \label{figureintroa-2}
\end{figure}

From this perception v-links are a combinatorial interpretation of projections on $\Sigma_g$. It is well-known that two v-link diagrams are equivalent iff their corresponding surface embeddings are \textit{stably equivalent}, i.e. equal modulo:
\begin{itemize}
\item The Reidemeister moves RM1, RM2 and RM3 and isotopies.
\item Adding/removing handles which do not affect the link diagram.
\item Homeomorphisms of surfaces.
\end{itemize}
For a proof see for example Proposition 6 in~\cite{cks}. For an example see Fig.~\ref{figureintroa-2}.

We are also interested in \textit{virtual tangle diagrams} and \textit{virtual tangles}. The first ones are graphs embedded in a disk $D^2$ such that each vertex is either one valent or four valent. The four valent vertices are, as before, labelled with an \textit{overcrossing} $\slashoverback$, an \textit{undercrossing} $\backoverslash$ or a \textit{virtual crossing} $\virtual$. The one valent vertices are part of the boundary of $D^2$ and we call them \textit{boundary points} and a virtual tangle diagram with $k$ one valent vertices a virtual tangle diagram with \textit{$k$-boundary points}.

A virtual tangle with $k$-boundary points is an equivalence class of virtual tangle diagrams with $k$-boundary points modulo the generalised Reidemeister moves and boundary preserving isotopies. We note that all of the moves in Fig.~\ref{figureintroa-1} can be seen as moves among virtual tangle diagrams. Examples are given later, e.g. in Sec.~\ref{sec-vkhcat}. As before, the notions of \textit{oriented} virtual tangle diagrams and \textit{oriented} virtual tangles can be defined analogously, but modulo \textit{oriented} generalised Reidemeister moves and boundary preserving isotopies.

If the reader is unfamiliar with the notion v-link or v-tangle, we refer to some introductory papers of Kauffman and Manturov, e.g.~\cite{ka2} and~\cite{kama}, and the references therein.
\vspace*{0.5cm}

\begin{figure}[ht]
     \centerline{\includegraphics[scale=0.8]{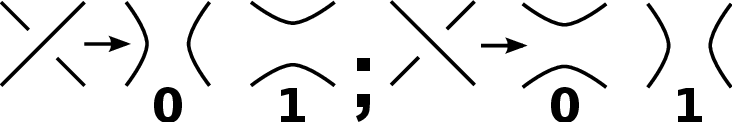}}
  \caption{The two possible resolutions of a crossing called \textit{0-resolution} and \textit{1-resolution}.}
  \label{figureintroa-3}
\end{figure}

Suppose one has a (classical) crossing $c$ in a diagram of a v-link (or an oriented v-link). We call a substitution of a crossing as shown in Fig.~\ref{figureintroa-3} a \textit{resolution} of the crossing $c$.

Furthermore, if we have a v-link diagram $L_D$, a \textit{resolution} of the v-link diagram $L_D$ is a diagram where all (classical) crossings of $L_D$ are replaced by one of the two resolutions from Fig.~\ref{figureintroa-3}. We use the same notions for v-tangle diagrams.

One of the greatest developments in modern knot theory was the discovery of \textit{Khovanov homology} by Khovanov in his famous paper~\cite{kh1} (Bar-Natan gave an exposition of Khovanov's construction in~\cite{bn1}). Khovanov homology is a \textit{categorification} of the Jones polynomial in the sense that the graded Euler characteristic of the \textit{Khovanov complex}, which we call the \textit{classical Khovanov complex}, is the Jones polynomial (up to normalisation). 

Recall that the Jones polynomial is known to be related to various parts of modern mathematics and physics, e.g. it origin lies in the study of von Neumann algebras. We note that the Jones polynomial can be extended to v-links in a rather straightforward way, see e.g. Sec. 5 in~\cite{ka1}. We call this extension the \textit{virtual Jones polynomial} or \textit{virtual $\mathfrak{sl}_2$ polynomial}. 

As a categorification, Khovanov homology reflects these connections on a ``higher level''. Moreover, the Khovanov homology of c-links is strictly stronger than its decategorification, e.g. see Sec. 4 in~\cite{bn1}. Another great development was the \textit{topological interpretation} of the Khovanov complex by Bar-Natan in~\cite{bn2}. This topological interpretation is a generalisation of the classical Khovanov complex for c-links and one of its modifications has functorial properties, see Theorem 1.1 in~\cite{cmw}. He constructed a \textit{topological complex} whose chain groups are formal direct sums of c-link resolutions and whose differentials are formal matrices of cobordisms between these resolutions.

\begin{figure}[ht]
     \centerline{\includegraphics[scale=0.55]{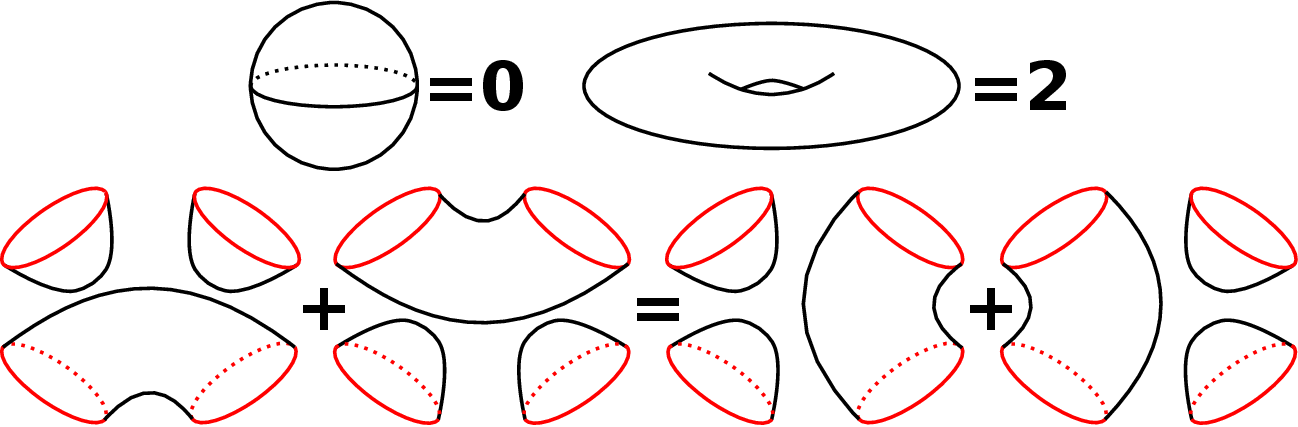}}
  \caption{A cobordism that contains a sphere $S$ is zero, a cobordism that contains a torus $T$ is two times the cobordism without the torus and the four tubes relation.}
  \label{figureintroa-4}
\end{figure}

Bar-Natan's construction modulo chain homotopy and the \textit{local relations} $S,T,4Tu$, also called \textit{Bar-Natan relations}, see Fig.~\ref{figureintroa-4}, is an invariant of c-links.

It is possible with this construction to classify all TQFTs which can be used to define c-link homologies from this approach, see Proposition 5 in~\cite{kh2}. Moreover, it is algorithmic, i.e. computable in less than exponential time (depending on the number of crossings of a given diagram), see~\cite{bn3}. So it is only natural to search for such a topological categorification of the virtual Jones polynomial.

An algebraic categorification of the virtual Jones polynomial over the ring $\bZ/2$ is rather straightforward and was done by Manturov in~\cite{ma1}. Moreover, he also published a version over the integers $\bZ$ later in~\cite{ma2}. A topological categorification was done by Turaev and Turner in~\cite{tutu}, but their version does not generalise Khovanov homology, since their complex is not bi-graded (see Sec. 4.2 in~\cite{tutu}). Another problem with their version is that it is not clear how to \textit{compute} the homology.
\vspace*{0.5cm}

We give a topological categorification which generalises the version of Turaev and Turner in the sense that a restriction of the version given here gives the topological complex of Turaev and Turner, another restriction gives a bi-graded complex that agrees with the Khovanov complex for c-links, and another restriction gives the so-called \textit{Lee complex}, i.e. a variant of the Khovanov complex that can be used to define the \textit{Rasmussen invariant} of a c-knot, see~\cite{ra}, which is also not included in the version of Turaev and Turner. Moreover, the version given here is computable and also strictly stronger than the virtual Jones polynomial.

Another restriction of the construction in his paper gives a different version than the one given by Manturov~\cite{ma2} in the sense that we conjecture it to be strictly stronger than his version. Moreover, in Secs.~\ref{sec-vkhtan} and~\ref{sec-vkhca}, we extend the construction to v-tangles in a ``good way'', something that is not known for Manturov's construction.

To be more precise, the categorification extends from c-tangles to v-tangles in a trivial way (by setting open saddles to be zero). This has an obvious disadvantage, i.e. it is neither a ``good'' invariant of v-tangles nor can it be used to calculate bigger complexes by ``tensoring'' smaller pieces. We give a local notion that is a strong invariant of v-tangles and allows ``tensoring''. We note that the construction for v-links is more difficult (combinatorially) than the classical case.
\section{A brief summary}\label{sec-vkhsum}
Let us give a brief, informal summary of the constructions in this paper. We will assume that the reader is not completely unfamiliar with the notion of the classical Khovanov complex as mentioned before, e.g. the construction of the Khovanov cube (more about cubes in Sec.~\ref{sec-techcube}) based on so-called \textit{resolutions of crossings} as shown in Fig.~\ref{figureintroa-3}. There are many good introductions to classical Khovanov homology, e.g. a nice exposition of the classical Khovanov homology can be found in Bar-Natan's paper~\cite{bn1}. We hope to demonstrate that the main ideas of the construction are easy, e.g. the construction is given by an algorithm, general, e.g. the construction extends all the ``classical'' homologies, but if one works over a ring $R$ of characteristic two, then, by setting $\theta\neq 0$, one obtains ``non-classical'' homologies. Moreover, the construction has other nice properties, e.g. it should have, up to a sign (?), functorial properties.
\vspace*{0.5cm}

Let $a$ be a word in the alphabet $\{0,1\}$. We denote by $\gamma_a$ the resolution of a v-link diagram $L_D$ with $|a|$ crossings, where the $i$-th crossing of $L_D$ is resolved $a_i\in\{0,1\}$ as indicated in Fig.~\ref{figureintroa-3}. Beware that we \textit{only} resolve classical crossings. We denote the number of v-circles, that is closed circles with only v-crossings, in the resolution $\gamma_a$ by $|\gamma_a|$.

Moreover, suppose we have two words $a,b$ with $a_k=b_k$ for $k=1,\dots,|a|=|b|,k\neq i$ and $a_i=0,b_i=1$ for a fixed $i\in\bN$. Then we call the expression $S\colon\gamma_a\to\gamma_b$ a \textit{(formal) saddle} between the resolutions.

Furthermore, suppose we have a v-link diagram $L_D$ with at least two crossings $c_1,c_2$. We call a quadruple $F=(\gamma_{00},\gamma_{01},\gamma_{10},\gamma_{11})$ of four resolutions of the v-diagram $L_D$ a \textit{face} of the diagram $L_D$, if in all four resolutions $\gamma_{00},\gamma_{01},\gamma_{10},\gamma_{11}$ all crossings of $L_D$ are resolved in the same way except that $c_1$ is resolved $i$ and $c_2$ is resolved $j$ in $\gamma_{ij}$ (with $i,j\in\{0,1\}$). Furthermore, there should be an oriented arrow from $\gamma_{ij}$ to $\gamma_{kl}$ if $i=j=0$ and $k=0,l=1$ or $k=1,l=0$ or if $i=0,j=1$ and $k=l=1$ or if $i=1,j=0$ and $k=l=1$. That is, faces look like
\[
\xymatrix{
 & \gamma_{01}\ar[rd]^{S_{*1}} &\\
 \gamma_{00}\ar[rd]_{S_{*0}}\ar[ru]^{S_{0*}} &  & \gamma_{11},\\
 & \gamma_{10}\ar[ru]_{S_{1*}} &}
\]
where the * for the saddles should indicate the change $0\to 1$. 
\vspace*{0.5cm}

We also consider \textit{algebraic faces} of a resolution. That is the same as above, but we replace $\gamma_a$ with $\bigotimes_n A$, if $\gamma_a$ has $n$ components. Here $A$ is an $R$-module and $R$ is a commutative, unital ring (which is for us usually of arbitrary characteristic).

Moreover, recall that the differential in the classical Khovanov complex consists of a multiplication $m\colon A\otimes A\to A$ and a comultiplication $\Delta\colon A\to A\otimes A$ for the $R$-algebra $A=R[X]/(X^2)$ with gradings $\deg 1=1,\deg X=-1$. The comultiplication $\Delta$ is given by
\[
\Delta\colon A\to A\otimes A;\begin{cases}

  1\mapsto 1\otimes X+X\otimes 1,\\
  X\mapsto X\otimes X.
\end{cases}
\]
The problem in the case of v-links is the emergence of a new map. This happens, because it is possible for v-links that a saddle $S\colon\gamma_a\to\gamma_b$ between two resolutions does not change the number of v-circles, i.e. $|\gamma_a|=|\gamma_b|$. This is a difference between c-links and v-links, i.e. in the first case one always has $|\gamma_a|=|\gamma_b|\pm 1$.
\vspace*{0.5cm}

Thus, in the algebraic complex we need a new map $\cdot\theta\colon A\rightarrow A$ together with the classical multiplication and comultiplication $m\colon A\otimes A\rightarrow A$ and $\Delta\colon A\rightarrow A\otimes A$. As we will see later the only possible way to extend the classical Khovanov complex to v-links is to set $\theta=0$ (for $R=\bZ$). But then a face could look like (maybe with extra signs).
\begin{equation}\label{probcube}
\begin{gathered}
\begin{xy}
  \xymatrix{
 & A\otimes A\ar[rd]^{m} &\\
 A\ar[ru]^{\Delta}\ar[rd]_{\cdot\theta} &  & A.\\
 & A\ar[ru]_{\cdot\theta} &}
\end{xy}
\end{gathered}
\end{equation}
We call such a face a \textit{problematic face}. With $\theta=0$ and the classical $\Delta,m$, this face does not commute (for $R=\bZ$). Therefore, there is no straightforward extension of the Khovanov complex to v-links. Moreover, in the cobordism based construction of the classical Khovanov complex, there is no corresponding cobordism for $\theta$.
\vspace*{0.5cm}

To solve these problems we consider a certain category called $\ucob_R(\emptyset)$, i.e. a category of (possibly non-orientable) cobordisms with boundary decorations $\{+,-\}$. Roughly, a punctured M\"obius strip plays the role of $\theta$ and the decorations keep track of how (orientation preserving or reversing) the surfaces are glued together. Hence, in our category we have different (co)multiplications, depending on the different decorations.
\vspace*{0.5cm}

Furthermore, in order to get the right signs, one has to use constructions related to $\wedge$-products. Note that this is rather surprising, since such constructions are not needed for Khovanov homology in the c-case where a ``usual'' sign placement suffices, see for example Sec. 3 in~\cite{bn1}. And furthermore, such constructions are in the c-case related to so-called \textit{odd} Khovanov homology introduced by Ozsv\'{a}th, Rasmussen and Szab\'{o} in~\cite{ors}. But we show that in fact our construction agrees for c-links with the (even) Khovanov homology (see Theorem~\ref{thm-classic}).
\vspace*{0.5cm}

The following table summarises the connection between the classical and the virtual case.
\vspace*{0.5cm}

\begin{center}
\begin{tabular}{|c|c|c|}
\hline $\phantom{.}$ & \textbf{Classical} & \textbf{Virtual} \\ 
\hline \textbf{Objects} & c-link resolutions & v-link resolutions \\ 
\hline \textbf{Morphisms} & Orientable cobordisms & Possibly non-orientable cobordisms \\ 
\hline \textbf{Cobordisms} & Embedded & Immersed \\ 
\hline \textbf{Decorations} & None & $+,-$ at the boundary \\ 
\hline \textbf{Signs} & Usual & Related to $\wedge$-products \\ 
\hline 
\end{tabular} 
\end{center}
\vspace*{0.5cm}

Hence, a main point in the construction of the virtual Khovanov complex is to say which saddles, i.e. morphisms, are orientable and which are non-orientable, how to place the decorations and how to place the signs.
\vspace*{0.5cm}

This is roughly done in the following way.
\begin{itemize}
\item Every saddle either splits one circle (orientable, called \textit{comultiplication}, denoted by $\Delta$. See Fig.~\ref{figure1-1} - fourth column), glues two circles (orientable, called \textit{multiplication}, denoted by $m$. See Fig.~\ref{figure1-1} - fifth column) or does not change the number of circles at all (non-orientable, called \textit{M\"obius cobordism}, denoted by $\theta$. See Fig.~\ref{figure1-1} - rightmost morphism).
\item Every saddle $S$ can be locally denoted (up to a rotation) by a formal symbol $S\colon\smoothing\rightarrow\hsmoothing$ (both smoothings are neighbourhoods of the crossing).
\item The \textit{glueing numbers}, i.e. the decorations, are now spread by \textit{choosing a formal orientation} for the resolution. We note that the construction will not depend on this choice (or on any other choices involved).
\item After all resolutions have an orientation, a saddle $S$ could for example be of the form $S\colon\du\rightarrow\ler$. This is (our choice for) the \textit{standard form}, i.e. in this case all glueing number will be $+$.
\item Now spread the decorations as follows. Every boundary component gets a $+$ iff the orientation is as in the standard form and a $-$ otherwise.
\item The \textit{degenerated} cases (everything non-alternating), e.g. $S\colon\dd\rightarrow\rir$, are the non-orientable surfaces and do not get any decorations. Compare to Table 1 in Definition~\ref{defn-deco}.
\item The signs are spread based on a numbering of the v-circles in the resolutions and on a special x-marker for the crossings. Note that without the x-marker one main lemma, i.e. Lemma~\ref{lem-virtualisation}, would not work.
\end{itemize}
\vspace*{0.5cm}

\begin{figure}[ht]
     \centerline{\includegraphics[scale=0.525]{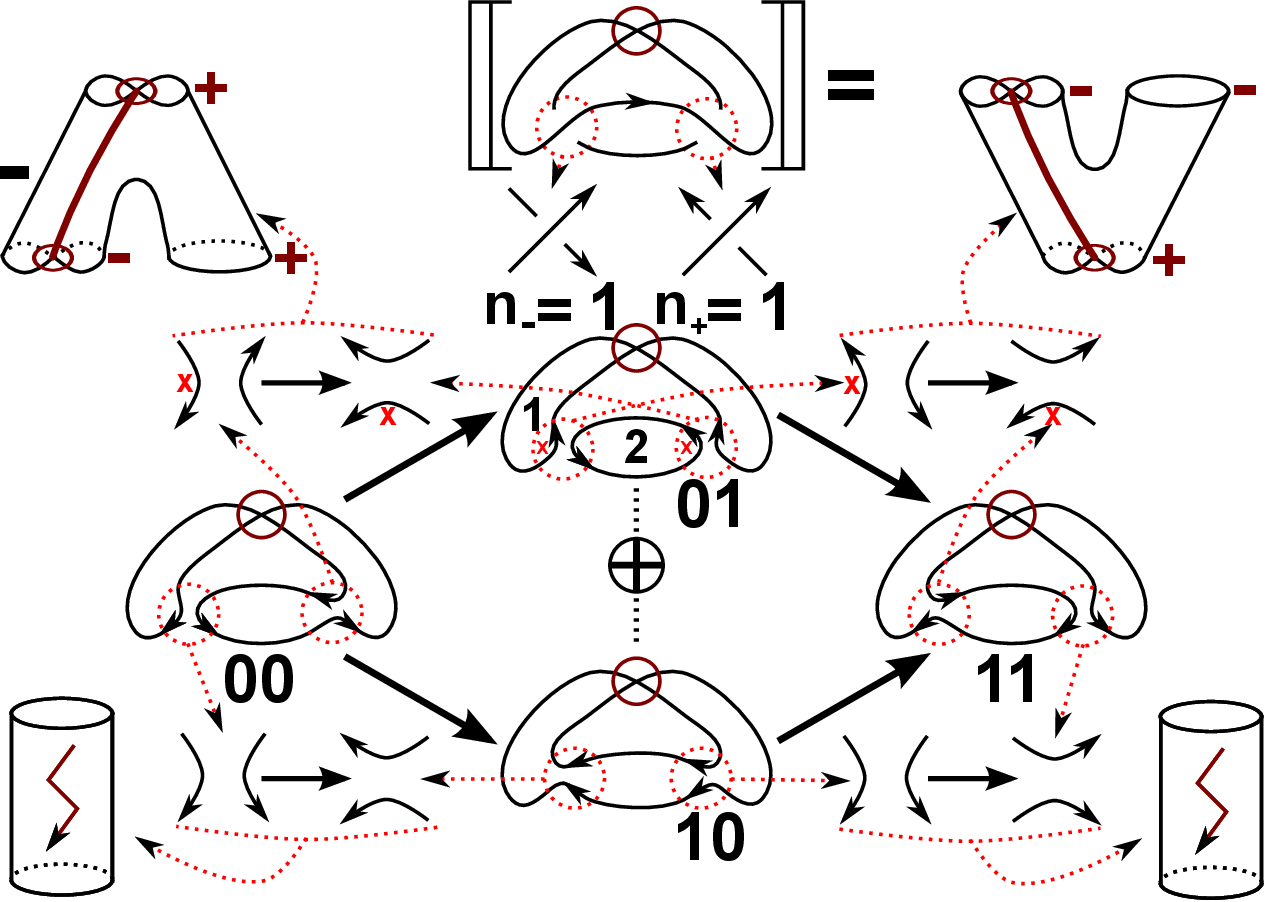}}
  \caption{The virtual Khovanov complex of the unknot.}\label{figure0-big}
\end{figure}

Or summarised in Fig.~\ref{figure0-big}. The complex below is the complex of a trivial v-link diagram.
\vspace*{0.5cm}

By our later construction, the homology degree zero part will be the direct sum in the middle. The bolt symbol indicates that the cobordism is non-orientable.
\vspace*{0.5cm}

We should note that this complex is exactly the problematic face from~\ref{probcube}.

\begin{figure}[ht]
     \centerline{\includegraphics[scale=0.75]{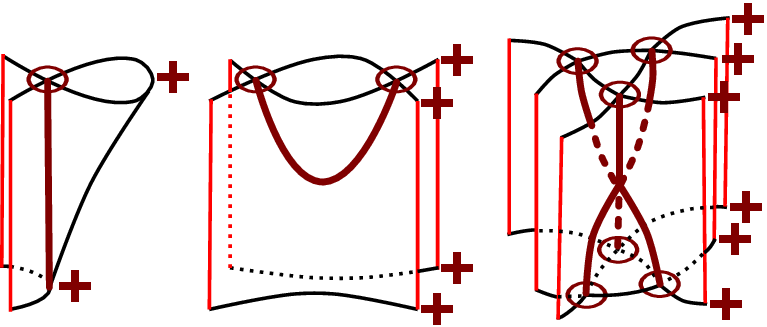}}
  \caption{The virtual Reidemeister cobordisms.}
  \label{figure0-reide}
\end{figure}

To construct the virtual Khovanov complex for v-tangles we need to extend these notions in such a way that they still work for ``open'' cobordisms. A first generalisation is easy, i.e. we will still use immersed, possibly non-orientable surfaces with decorations, but we allow vertical boundary components, e.g. the three \textit{v-Reidemeister cobordisms vRM1, vRM2 and vRM3} in Fig.~\ref{figure0-reide}. We note that we always read cobordisms from top to bottom, i.e. the first two cobordisms simplify the v-tangle diagram.

One main point is the question what to do with the \textit{``open'' saddles}, i.e. saddles with no closed boundary. A possible solution is to define them to be zero.
\vspace*{0.5cm}

\begin{figure}[ht]
     \centerline{\includegraphics[scale=0.5]{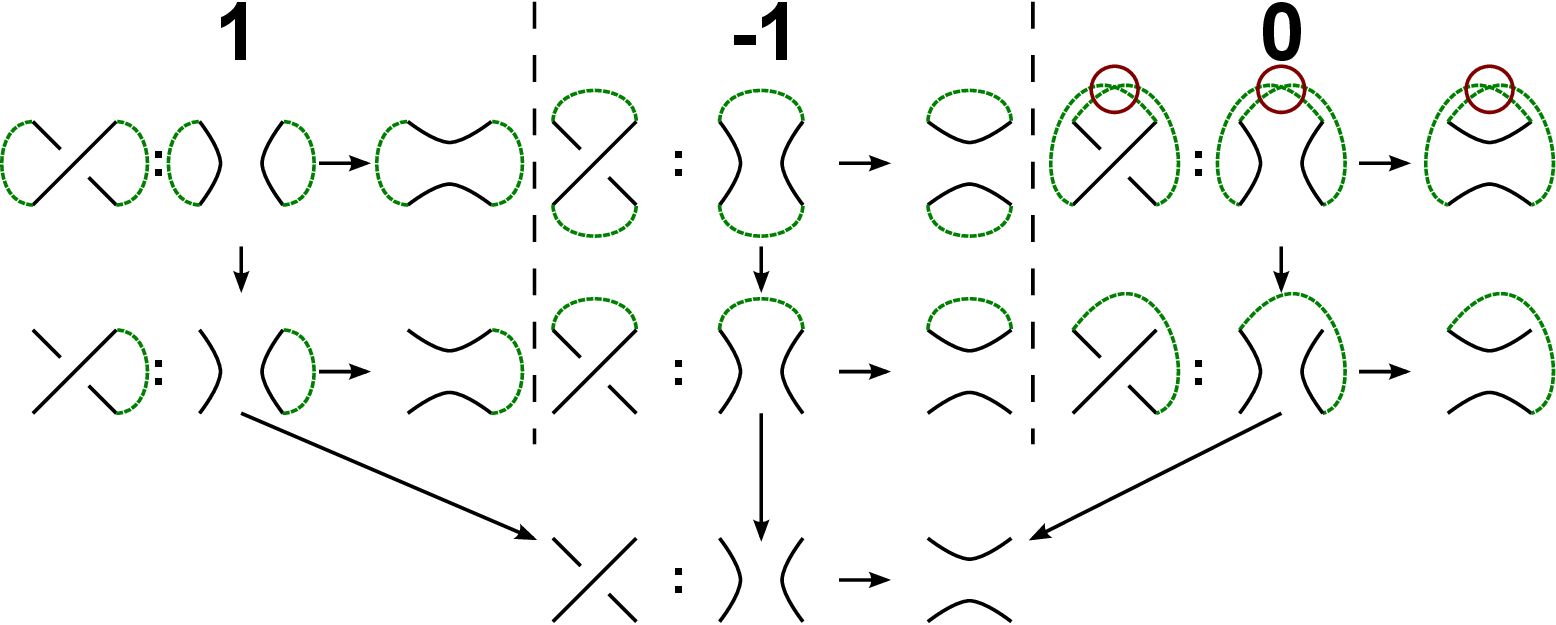}}
  \caption{All of the closed cases give rise to the unclosed.}
  \label{figure0-order}
\end{figure}

But this has two major problems. First the loss of information is big and second we would not have local properties as in the classical case (``tensoring'' of smaller parts), since an open saddle can, after closing some of his boundary circles, become either $m$, $\Delta$ or $\theta$, see Fig.~\ref{figure0-order}. This figure illustrates that we can never be sure what ``type'' of saddle a local saddle will be after glueing it inside a bigger piece.
Hence, an information mod 3 is missing. We therefore consider morphisms with an \textit{indicator}, i.e. an element of the set $\{0,+1,-1\}$. Then, after taking care of some technical difficulties, the concept extends from c-tangles to v-tangles in a suitable way. It should be noted that we do not know how to spread signs locally. But we get a so-called \textit{projective complex}, i.e. faces commute up to a unit. We have collected some of the technical points in Sec.~\ref{sec-techcube}.

\begin{figure}[ht]
     \centerline{\includegraphics[scale=0.5]{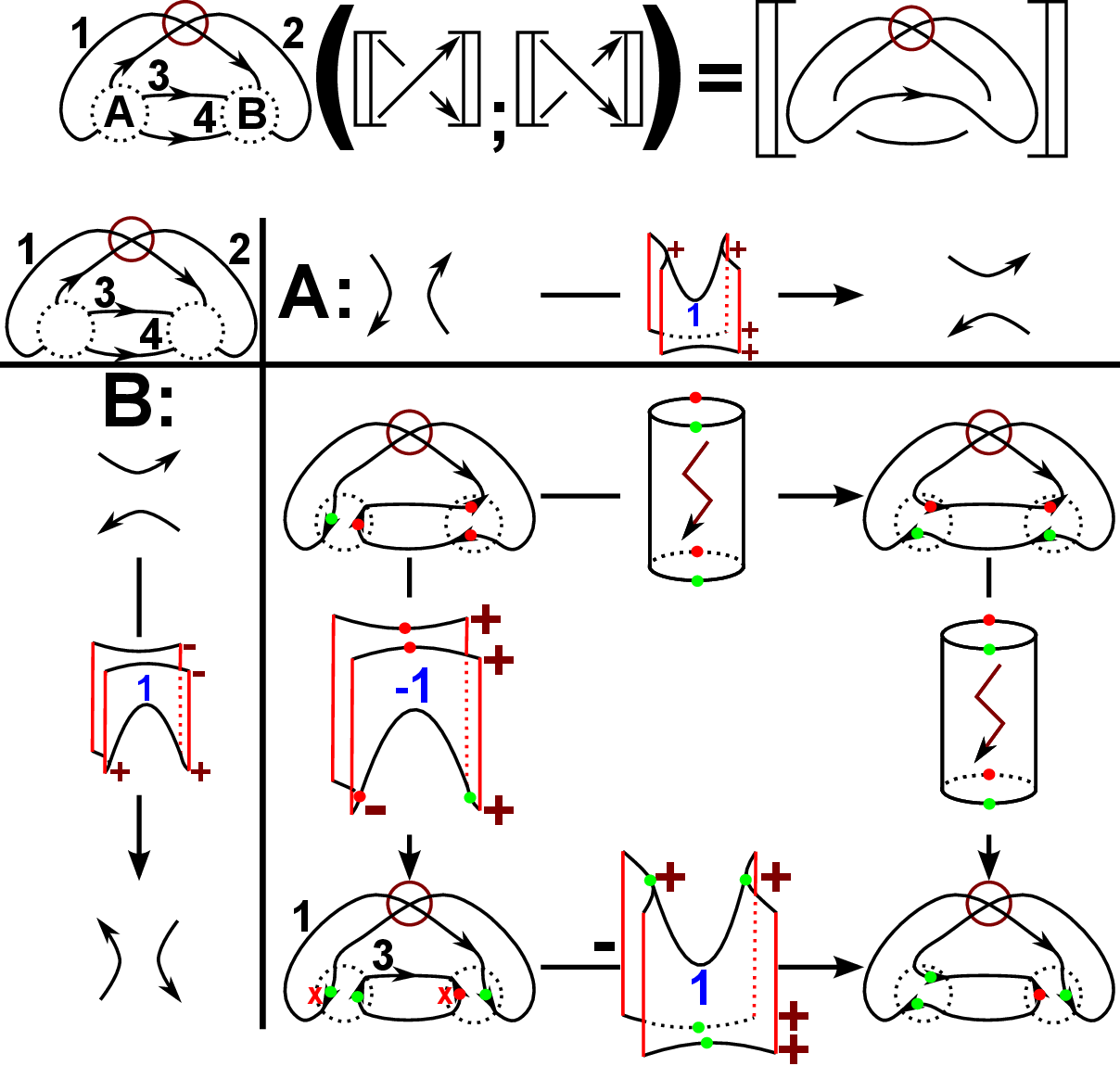}}
  \caption{After we have fixed an orientation/numbering of the circuit diagram, we only have to compare whether the local orientations match (green) or mismatch (red) and compose if necessary with $\Phi^-_+$ (red). Iff we have a double mismatch at the top and bottom, then we add a bolt symbol.}
  \label{figure0-main}
\end{figure}

Then, after taking care of some difficulties again, we can ``tensor'' smaller pieces together as indicated in the Fig.~\ref{figure0-main}.
It should be noted that there are some technical points that make our construction only semi-local (a disadvantage that arises from the fact that ``non-orientability'' is not a local property). Note that indicators, if necessary, are pictured on the surfaces.
\vspace*{0.5cm}

The outline of the paper is as follows.
\begin{itemize}
\item In Sec.~\ref{sec-vkhcat} we define the category of (possibly non-orientable) cobordisms with boundary decorations. First in the ``closed'' case in Definition~\ref{defn-category} and then more generally in the ``open'' case in Definition~\ref{defn-category3}. We also prove/recall some basic facts in Sec.~\ref{sec-vkhcat}.
\item In Sec.~\ref{sec-vkhcom} we define the virtual Khovanov complex for v-links in Definition~\ref{defn-topcomplex}. We show that it is a v-link invariant (see Theorem~\ref{thm-geoinvarianz}) and agrees with the construction in the c-case (see Theorem~\ref{thm-classic}). There are two important things about the construction. 
\begin{itemize}
\item The first is that there are many choices in the definition of the virtual complex, but we show in Lemma~\ref{lem-commutativeindependence} that different choices give isomorphic complexes.
\item Second, it is not clear that the complex is a well-defined chain complex, but we show this fact in Theorem~\ref{theo-facescommute} and Corollary~\ref{cor-chaincomplex}. In order to show that the construction gives a well-defined chain complex we have to use a ``trick'', i.e. we use a move called \textit{virtualisation}, as shown in Fig.~\ref{figure0-virt}, to reduce the question whether the faces of the virtual Khovanov cube are anti-commutative to a finite and small number of so-called \textit{basic faces} (see Fig.~\ref{figure-basic}).
\begin{figure}[ht]
     \centerline{\includegraphics[scale=0.6]{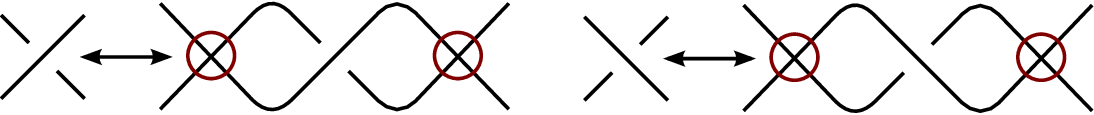}}
  \caption{The virtualisation of a crossing.}
  \label{figure0-virt}
\end{figure}
\end{itemize}
\item In Sec.~\ref{sec-vkhfa} we show that our constructions can be compared to so-called \textit{skew-extended Frobenius algebras} (Theorem~\ref{thm-tututheo}). With this we are able to classify all possible v-link homologies from our approach, see Theorem~\ref{thm-class}. We note that all the classical homologies are included. And we can therefore show in Corollary~\ref{cor-khovhomocat} that our construction is a categorification of the virtual Jones polynomial.
\item Secs.~\ref{sec-vkhtan} and~\ref{sec-vkhca} are analoga of the earlier sections, but for v-tangles.
\item Sec.~\ref{sec-vkhapp} uses that our construction is semi-local (see Theorem~\ref{thm-semiloc}). As a result, we can still show that Lee's variant of Khovanov homology is in some sense degenerated, see Theorem~\ref{thm-leedeg}. This fact is one of the main ingredients to define Rasmussen's invariant in the classical case.
\item Sec.~\ref{sec-vkhpro} gives some calculation results with a MATHEMATICA program written by the author. It is worth noting that we give examples of v-links with seven crossings which can not be distinguished by the virtual Jones polynomial, but by virtual Khovanov homology.
\item We have collected some known facts about strong deformation retracts and cube complexes in the two Secs.~\ref{sec-techhomalg} and~\ref{sec-techcube} before the final section.
\item Moreover, we have collected some open questions in Sec.~\ref{sec-vkhend}. 
\end{itemize}
\subsection{Notation and Remarks}\label{notation}
We call $\smoothing$ the \textit{$0$-} and $\hsmoothing$ the \textit{$1$-resolution} of the crossing $\slashoverback$ for a given v-link diagram $L_D$ or v-tangle diagram $T^k_D$. For an oriented v-link diagram $L_D$ or v-tangle diagram $T^k_D$ we call $\overcrossing$ a \textit{positive} and $\undercrossing$ a \textit{negative} crossing. The \textit{number of positive crossings} is denoted by $n_+$ and the \textit{number of negative crossings} is denoted by $n_-$.

For a given v-link diagram $L_D$ or v-tangle diagram $T^k_D$ with $n$ numbered crossings we define a collection of closed curves and open strings $\gamma_a$ in the following way. Let $a$ be a word of length $n$ in the alphabet $\left\{0,1\right\}$. Then $\gamma_a$ is the collection of closed curves and open strings which arise, when one performs a $a_i$-resolution at the $i$-th crossing for all $i=1,\dots,n$. We call such a collection $\gamma_a$ the \textit{$a$-th resolution} of $L_D$ or $T^k_D$. All appearing v-circles should be numbered with consecutive numbers from $1,\dots,k_a$ in these resolutions, where $k_a$ is the total number of v-circles of the resolution $\gamma_a$.

We can choose an orientation for the different components of $\gamma_a$. We call such a $\gamma_a$ an \textit{orientated resolution}, i.e. every v-crossing of the resolution $\gamma_a$ should look like $\virtualoriented$. Then a local neighbourhood of a $0,1$-resolved crossing could for example look like $\uu$. We call these neighbourhoods \textit{orientated crossing resolutions}.

If we ignore orientations, then there are $2^n$ different resolutions $\gamma_a$ of $L_D$ or $T^k_D$. We say a resolution has length $m$ if it contains exactly $m$ 1-letters. That is $m=\sum_{i=1}^na_i$.

For two resolutions $\gamma_a$ and $\gamma_{a^{\prime}}$ with $a_r=0$ and $a^{\prime}_r=1$ for one fixed $r$ and $a_i=a^{\prime}_i$ for $i\neq r$ we define a \textit{saddle between the resolutions $S$}. This means: Choose a small (no other crossing, classical or virtual, should be involved) neighbourhood $N$ of the $r$-th crossing and define a cobordism between $\gamma_a$ and $\gamma_{a^{\prime}}$ to be the identity outside of $N$ and a saddle inside of $N$. Note that we, by a slight abuse of terminology, call these cobordisms saddles although they contain in general some cylinder components.

From now on we consider \textit{faces} $F=(\gamma_{00},\gamma_{01},\gamma_{10},\gamma_{11})$ of four resolutions, as mentioned above, always \textit{together with the saddles} between the resolutions. We denote the saddles for example by $S_{0*}\colon\gamma_{00}\to\gamma_{01}$, where the position of the $*$ indicates the change $0\to 1$.

It should be noted that any v-link or v-tangle diagram should be oriented in the usual sense. But with a slight abuse of notation, we will suppress this orientation throughout the whole paper, since the afore mentioned oriented resolutions are main ingredients of our construction and easy to confuse with the usual orientations. Recall that these usual orientations are needed for the shifts in homology gradings, see for example Sec. 3 in~\cite{bn1}.
\vspace*{0.5cm}

\begin{figure}[ht]
     \centerline{\includegraphics[width=0.75\linewidth]{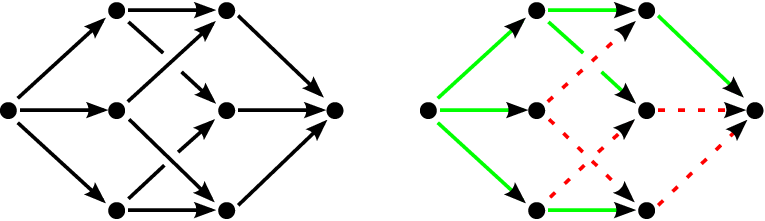}}
  \caption{A Khovanov cube and a spanning tree of the cube (green edges).}
  \label{figure0-spantree}
\end{figure}

Sometimes we need a so-called \textit{spanning tree argument}, i.e. choose a spanning tree of a cube (as in Fig.~\ref{figure0-spantree}) and change e.g. orientations of resolutions such that the edges of the tree change in a suitable way, starting at the rightmost leaves, then remove the rightmost leaves and repeat. Notice that two cubes together with a chain map between them form again a bigger cube. It is worth noting that most of the spanning tree arguments work out in the end because of certain preconditions, e.g. the anti-commutativity of faces.

Moreover, we have collected some facts from homological algebra that we need in this paper in Sec.~\ref{sec-techhomalg} and in Sec.~\ref{sec-techcube}.

\begin{rem}\label{rem-consistency}
We note that it is a priori not clear why Definition~\ref{defn-category} gives something non-trivial. We show later in Sec.~\ref{sec-vkhfa} that there is a model of the topological category $\ucob_R(\emptyset)$. One particular blueprint model takes values in $\RMOD$ (the category of $R$-modules) for $R=\mathbb{Z}$, see Table 2 with $a=1$, $\alpha=\beta=\gamma=t=0$. The reader may think of $\Phi_+^-$ as the base change matrix $1\mapsto 1$ and $X\mapsto -X$ for the algebra $A=\mathbb{Z}[X]/(X^2)$. But this is just one possible model for $\ucob_R(\emptyset)$.
\end{rem}

\begin{rem}(\textbf{Remark about colours})
We use colours in this paper. If the reader does not use colours or has only access to the uncoloured version:

The main difference is between things coloured red and green. These two colours can be distinguished in uncoloured versions because the red one will be shaded darker. For example in
\[
\smoothinga\simeq\ddp\oplus\dup\oplus\udp\oplus\uup\;\;\text{ and }\;\;\hsmoothinga\simeq\ddpp\oplus\dupp\oplus\udpp\oplus\uupp.
\]
one can distinguish between red=u and green=d by their darkness. That is, every appearance of sentences similar to ``the red'' should be read as ``the darker shaded''.  
\end{rem}
\section{The topological category}\label{sec-vkhcat}
\subsection{The topological category for v-links}
In this section we describe our topological category which we call $\ucob_R(\emptyset)$. This is a category of cobordisms between v-link resolutions in the spirit of Bar-Natan~\cite{bn2}, but we admit that the cobordisms are non-orientable as in~\cite{tutu}.
\vspace*{0.5cm}

The basic idea of the construction is that the usual pantsup- and pantsdown-cobordisms do not satisfy the relation $m\circ\Delta=\theta^2$. But we need this relation for the face from~\ref{probcube}. This is the case, because we need an extra information for v-links, namely \textit{how} two cobordisms are glued together.

To deal with this problem, we decorate the boundary components of a cobordism with a formal sign $+,-$. With this construction $m_i\circ\Delta_j$ is sometimes $=\theta^2$ and sometimes $\neq\theta^2$, depending on the boundary decorations, which are here represented by indices $i,j=1,\dots,8$. The first case will occur iff $m_i\circ\Delta_j$ is a non-orientable surface.

\begin{figure}[ht]
  	 \centerline{\xy
     (0,0)*{\includegraphics[scale=0.5]{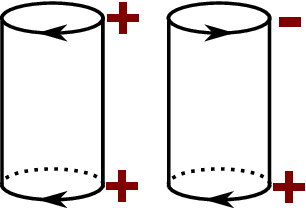}};
     (-8,0)*{\mathrm{id}^+_+};
     (6,0)*{\Phi^-_+};
     \endxy}
  \caption{Glueing the boundary together as indicated can not be done without immersion in the case on the right.}
  \label{figure-idphi}
\end{figure}

One main idea of this construction is the usage of a cobordism $\Phi^-_+$ between two circles \textit{different} from the identity $\mathrm{id}^+_+$, see Fig.~\ref{figure-idphi}.

Furthermore, we need relations between the decorated cobordisms. One of these relations identifies all boundary-preservingly homeomorphic cobordisms if their boundary decorations are all equal or are all different (up to a sign). Moreover, some of the standard relations of the category $\cob_R(\emptyset)$ (see for example in the book of Kock~\cite{ko}, Chapter 1.4) should hold. We denote the category with the extra signs by $\ucob_R(\emptyset)$ and the category without the extra signs by $\ucob_R(\emptyset)^*$. Therefore, there will be two different cylinders in these categories.

Note that most of the constructions are easier for $\ucob_R(\emptyset)^*$ than for $\ucob_R(\emptyset)$. That is why we will only focus on the latter category and hope the reader does not have too many difficulties to do similar constructions for $\ucob_R(\emptyset)^*$ while reading this subsection. 

At the end of this subsection we will prove some basic relations (see Lemma~\ref{lem-basiscalculations}) between the generators of our category. We also characterise the cobordisms of the face from~\ref{probcube} (see Proposition~\ref{prop-nonorientablefaces}).
\vspace*{0.5cm}

It should be noted that, in order to extend the construction to v-tangle diagrams, we need some more extra notions. We will define them after Definition~\ref{defn-category} in an extra subsection in Definition~\ref{defn-category3} to avoid too many notions at once.

We start with the following definition. Beware that we consider v-circles as objects and cobordisms together with decorations. We denote the decorations by $+,-$ and illustrate them next to boundary components.

We note again that $R$ denotes a commutative, unital ring of \textit{arbitrary} characteristic.
\begin{defn}\label{defn-category}\textbf{(The category of cobordisms with boundary decorations)} We describe the category \textit{$\ucob_R(\emptyset)$} in six steps. Note that our category is $R-$pre-additive\footnote{Sometimes also called $R$-category, i.e. the set of morphisms form a $R$-module and composition is $R$-linear. Or said otherwise, the category should be enriched over $\RMOD$}. The symbol $\amalg$ denotes the disjoint union.

\textbf{The objects:}

The \textit{objects} \textit{$\Ob(\ucob_R(\emptyset)$}) are disjoint unions of numbered \textit{v-circles} (recall that v-circles are circles without c-crossings). We denote the objects by $\mathcal O=\amalg_{i\in I}\mathcal O_i$. Here $\mathcal O_i$ are the v-circles and $I$ is a finite, ordered index set. Note that, by a slight abuse of notation, we denote the objects by $\mathcal O$ to point out that the category can be seen as a $2$-category with v-circles as $1$-morphisms between the empty set (but this is inconvenient for our purpose). The objects of the category are equivalence (modulo \textit{planar isotopies}) classes of four-valent graphs.

\textbf{The generators:}

The \textit{generators} of $\Mor(\ucob_R(\emptyset))$ are the eight cobordisms from Fig.~\ref{figure1-1} plus topologically equivalent cobordisms, but with all other possible boundary decorations (we do not picture them because one can obtain them using the ones shown after taking the relations below into account). Every orientable generator has a decoration from the set $\{+,-\}$ at the boundary components. We call these decorations the \textit{glueing number} (of the corresponding boundary component).
\begin{figure}[ht]
  	 \centerline{\xy
     (0,0)*{\includegraphics[scale=0.5]{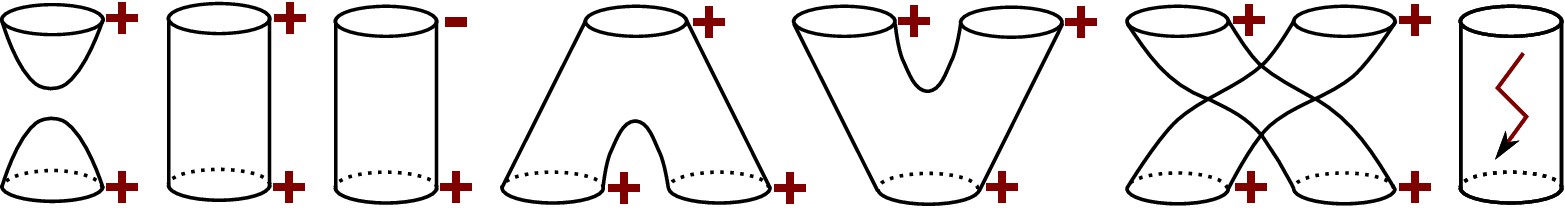}};
     (-67.75,-3.8)*{\varepsilon^+};
     (-67.75,3.8)*{\iota_+};
     (-47.15,-11)*{\mathrm{id}^+_+};
     (-31.75,-11)*{\Phi^-_+};
     (14.5,-11)*{m^{++}_+};
     (-10,-11)*{\Delta_{++}^+};
     (42.5,-11)*{\tau^{++}_{++}};
     (61.5,-11)*{\theta};
     \endxy}
  \caption{The generators of the set of morphisms. The cobordism on the right is the M\"obius cobordism, i.e. a two times punctured projective plane.}
  \label{figure1-1}
\end{figure}

We consider these cobordisms up to boundary preserving homeomorphisms (as abstract surfaces). Hence, between circles with v-crossings the (not pictured) generators are the same up to boundary preserving homeomorphisms, but immersed into $\bR^2\times[-1,1]$.

The eight cobordisms are (from left to right): a \textit{cap-cobordism} and a \textit{cup-cobordism} between the empty set and one circle and vice versa. Both are homeomorphic to a disc $D^2$ and both have a positive glueing number. We denote them by $\iota_+$ and $\varepsilon^+$ respectively.

Two \textit{cylinders} from one circle to one circle. The first has two positive glueing numbers and we denote this cobordism by $\mathrm{id}^+_+$. The second has a negative upper glueing number and a positive lower glueing number and we denote it by $\Phi^-_+$.

A \textit{multiplication-} and a \textit{comultiplication-cobordism} with only positive glueing numbers. Both are homeomorphic to a three times punctured $S^2$. We denote them by $m^+_{++}$ and $\Delta^+_{++}$.

A \textit{permutation-cobordism} between two upper and two lower boundary circles with only positive glueing numbers. We denote it by $\tau^{++}_{++}$.

A \textit{two times punctured projective plane}, also called \textit{M\"obius cobordism}. This cobordism is not orientable, hence it has no glueing numbers. We denote it by $\theta$.

The composition of the generators is given by glueing them together along their common boundary. In all pictures the upper cobordism is the $C$ in the composition $C^{\prime}\circ C$. The decorations are not changing at all (except that we remove the decorations if any connected component is non-orientable) \textit{before} taking the relations as in the equations~\ref{eq-combrel12},~\ref{eq-combrel3},~\ref{eq-commasss},~\ref{eq-unitrel},~\ref{eq-permrel1},~\ref{eq-permrel2} and~\ref{eq-frobandco} into account. Formally, \textit{before} taking quotients, the composition of the generators also needs internal decorations to remember if the generators were glued together alternatingly, i.e. minus to plus or plus to minus, or non-alternatingly. But after taking the quotients as indicated, these internal decorations are not needed any more. Hence, we suppress these internal decorations to avoid a too messy notation.

The reader should keep the informal slogan ``Composition with $\Phi^-_+$ changes the decoration'' in mind. 

\textbf{The morphisms:}

The \textit{morphisms} \textit{$\Mor(\ucob_R(\emptyset))$} are cobordisms between the objects in the following way. Note that we call a morphism non-orientable if any of its connected components is non-orientable.

We identify the collection of numbered v-circles with circles immersed into $\bR^2$. Given two objects $\mathcal O_1,\mathcal O_2$ with $k_1,k_2$ numbered v-circles, a morphism $C\colon\mathcal O_1\to\mathcal O_2$ is a surface immersed in $\bR^2\times[-1,1]$ whose boundary lies only in $\bR^2\times\{-1,1\}$ and is the disjoint union of the $k_1$ numbered v-circles from $\mathcal O_1$ in $\bR^2\times\{1\}$ and the disjoint union of the $k_2$ numbered v-circles from $\mathcal O_2$ in $\bR^2\times\{-1\}$. The morphisms are generated (as abstract surfaces) by the generators from above. It is worth noting that all possible boundary decorations can occur.

\textbf{The decorations:}

Given a $C\colon\mathcal O_1\to\mathcal O_2$ in $\Mor(\ucob_R(\emptyset))$, let us say that the v-circles of $\mathcal O_1$ are numbered from $1,\dots,k$ and the v-circles of $\mathcal O_2$ are numbered from $k+1,\dots,l$.

Every orientable cobordism has a decoration on the $i$-th boundary circle. This decoration is an element of the set $\{+,-\}$. We call this decoration of the $i$-th boundary component the \textit{$i$-th glueing number} of the cobordism.

Hence, the morphisms of the category are pairs $(C,w)$. Here $C\colon\mathcal O_1\to\mathcal O_2$ is a cobordism from $\mathcal O_1$ to $\mathcal O_2$ immersed in $\bR^2\times[-1,1]$ and $w$ is a string of length $l$ in such a way that the $i$-th letter of $w$ is the $i$-th glueing number of the cobordism or $w=0$ if the cobordism is non-orientable.

\textbf{Shorthand notation:}

We denote an orientable, connected morphism $C$ by $C^{u}_{l}$. Here $u,l$ are words in the alphabet $\{+,-\}$ in such a way that the $i$-th character of $u$ (of $l$) is the glueing number of the $i$-th circle of the upper (of the lower) boundary. The construction above ensures that this notation is always possible. Therefore, we denote an arbitrary orientable morphism $(C,w)$ by
\[
C=C^{u_1}_{l_1}\amalg\cdots\amalg C^{u_k}_{l_k}.
\]
Here $C^{u_i}_{l_i}$ are its connected components and $u_i,l_i$ are words in $\{+,-\}$. For a non-orientable morphism we do not need any boundary decorations.

\textbf{The relations:}

There are two different types of relations, namely \textit{topological relations} and \textit{combinatorial relations}. The latter relations have to do with the glueing numbers and the glueing of the cobordisms. The relations between the morphisms are the relations pictured below, i.e. the three \textit{combinatorial}~\ref{eq-combrel12} for the orientable and~\ref{eq-combrel3} for non-orientable cobordisms, \textit{commutativity} and \textit{cocommutativity} relations~\ref{eq-commasss}, \textit{associativity} and \textit{coassociativity} relations~\ref{eq-commasss}, \textit{unit and counit} relations~\ref{eq-unitrel}, \textit{permutation} relations~\ref{eq-permrel1} and~\ref{eq-permrel2}, a \textit{Frobenius relation} and the \textit{torus and M\"obius} relations~\ref{eq-frobandco} and different \textit{commutation} relations. Latter ones are not pictured, but all of them should hold with a plus sign. If the reader is unfamiliar with these relations, then we refer to the book of Kock~\cite{ko} (Chapter 1.4) and hope (we really do) that it should be clear how to translate his pictures to our context (by adding some decorations). 

Beware that we have pictured several relations in some figures at once. We have separated them by a thick line.

Moreover, some of the relations contain several cases at once, e.g. in the right part of Equation~\ref{eq-frobandco}. In those cases it should be read: If the conditions around the equality sign are satisfied, then the equality holds.

The first combinatorial relations are (read the right part as ``$\Phi_+^-$ changes the glueing numbers'')
\begin{align}\label{eq-combrel12}
\xy(0,0)*{\includegraphics[scale=0.25]{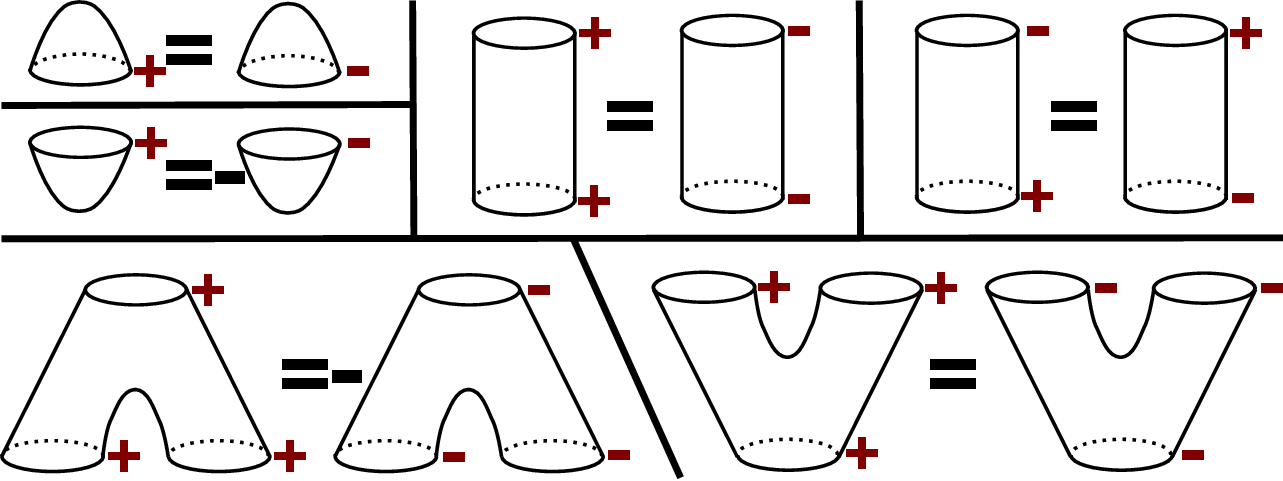}};\endxy\;\;\;\;\;\;\;\;\xy(0,0)*{\includegraphics[scale=0.25]{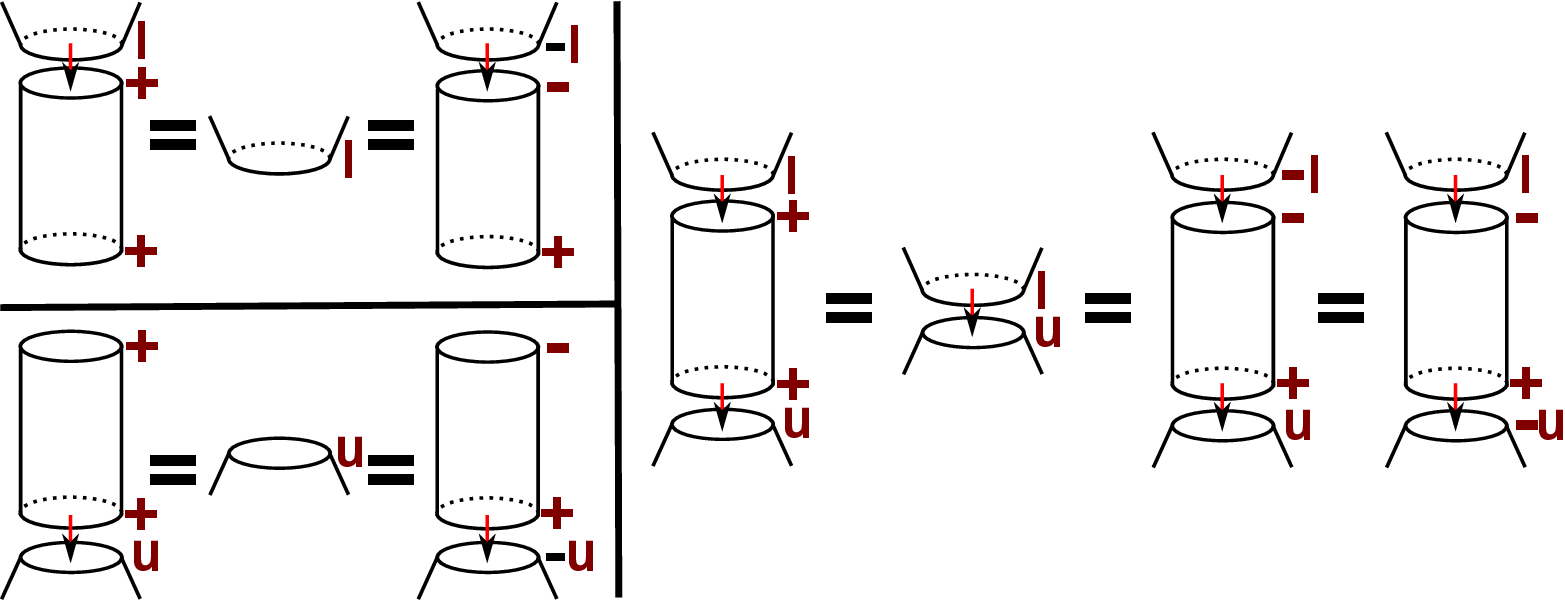}};\endxy
\end{align}
and the third for the non-orientable cobordisms is
\begin{align}\label{eq-combrel3}
\xy(0,0)*{\includegraphics[scale=0.25]{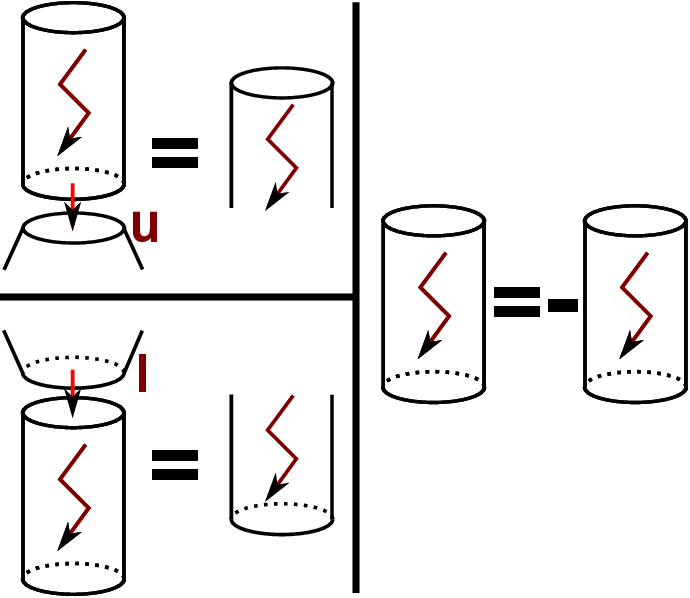}};\endxy
\end{align}
Note that the relation~\ref{eq-combrel3} above is not the same as $\theta=0$, since we work over rings of arbitrary characteristic. The (co)commutativity and (co)associativity relations are
\begin{align}\label{eq-commasss}
\xy(0,0)*{\includegraphics[scale=0.25]{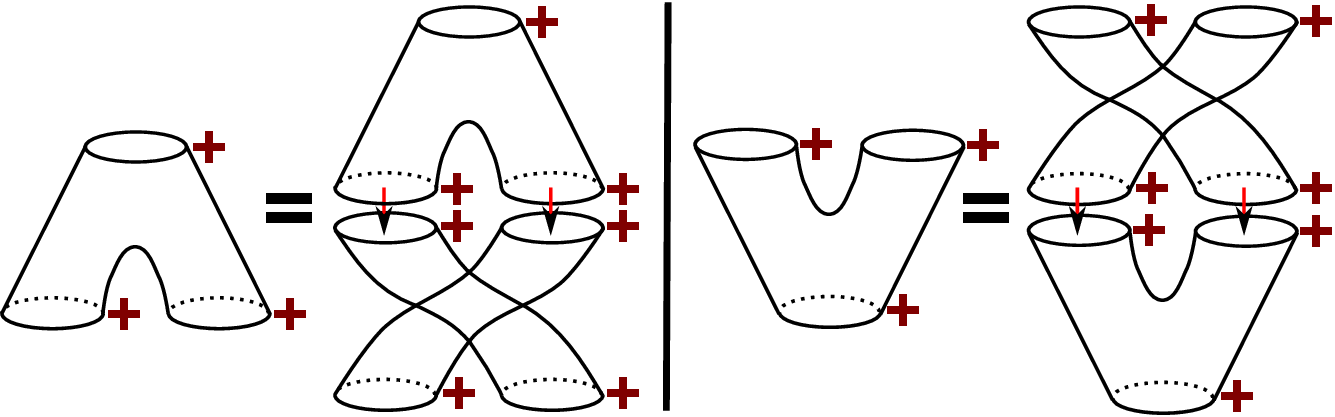}};\endxy\;\;\;\;\;\;\;\;\xy(0,0)*{\includegraphics[scale=0.25]{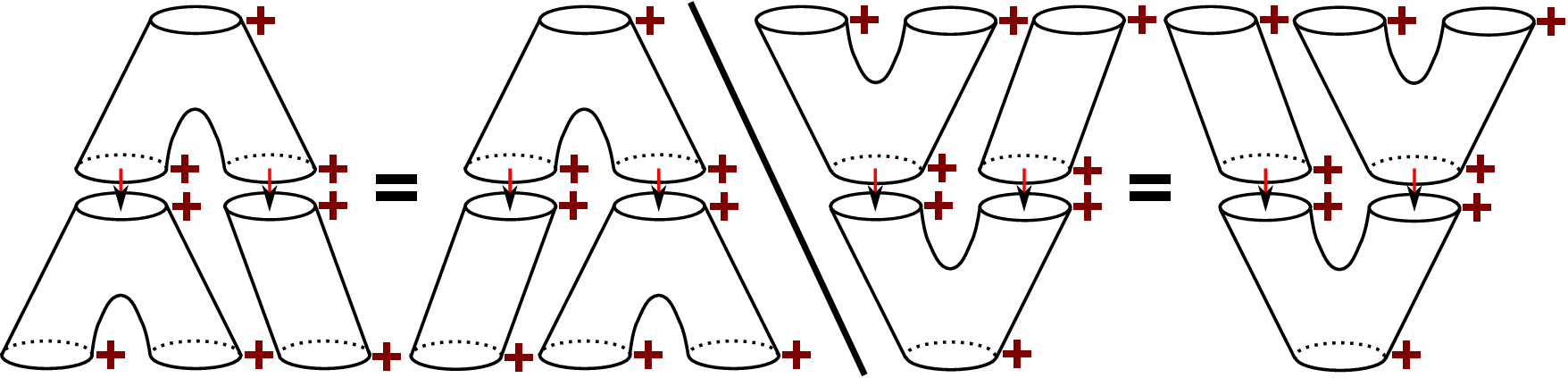}};\endxy
\end{align}
and the (co)unit relations are
\begin{align}\label{eq-unitrel}
\xy(0,0)*{\includegraphics[scale=0.25]{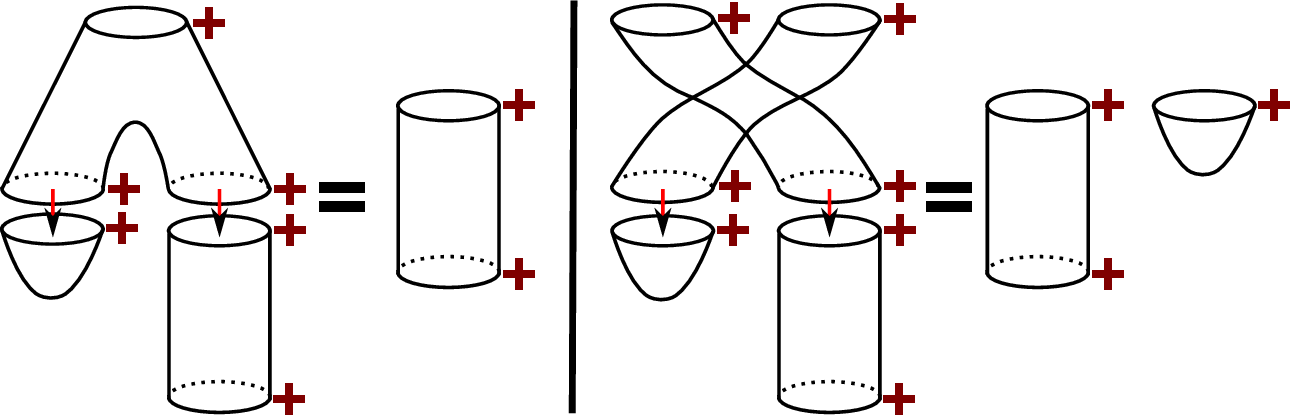}};\endxy\;\;\;\;\;\;\;\;\xy(0,0)*{\includegraphics[scale=0.25]{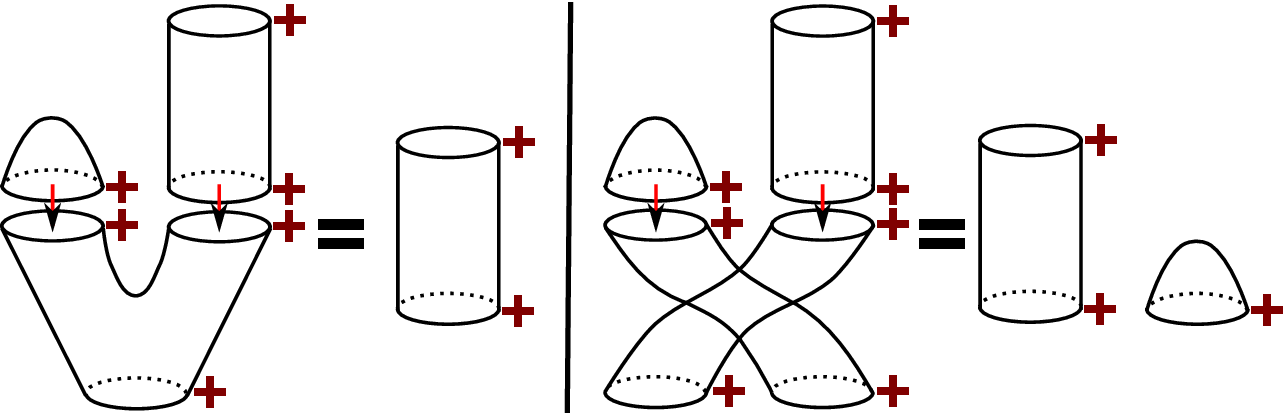}};\endxy
\end{align}
The first and second permutation relations are
\begin{align}\label{eq-permrel1}
\xy(0,0)*{\includegraphics[scale=0.25]{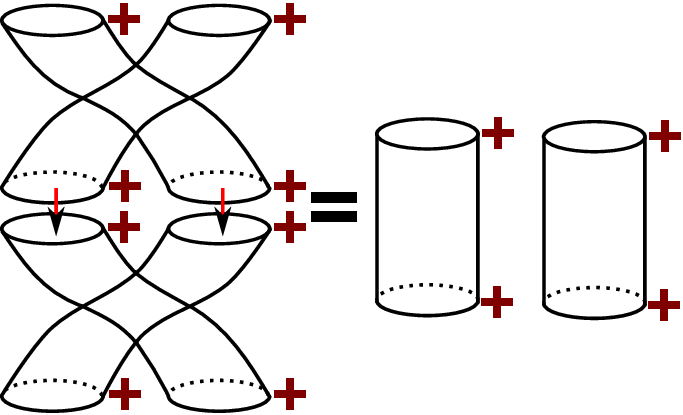}};\endxy\;\;\;\;\;\;\;\;\xy(0,0)*{\includegraphics[scale=0.25]{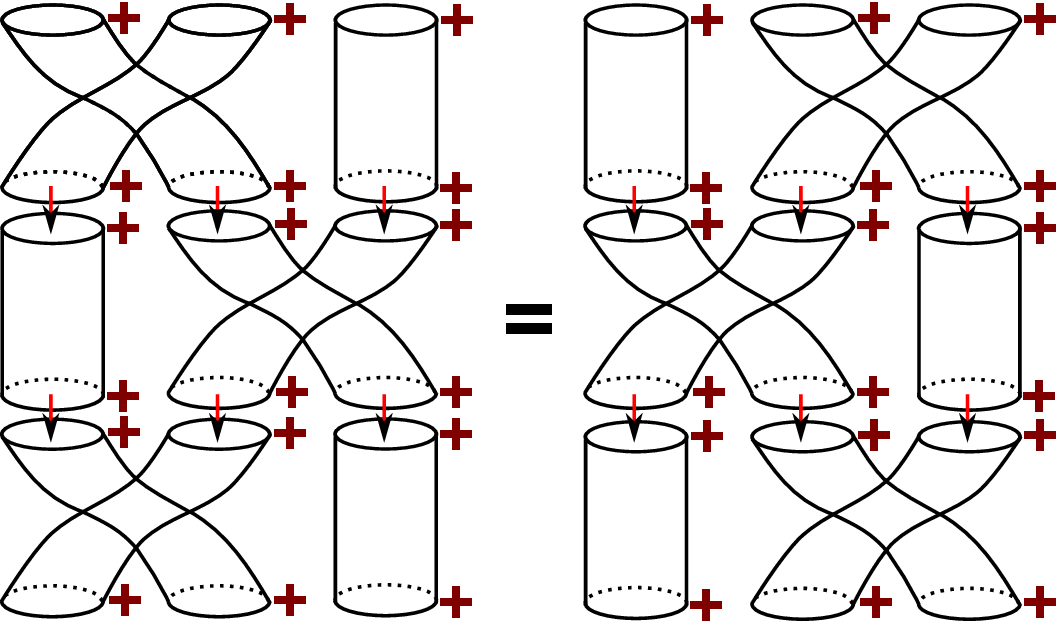}};\endxy
\end{align}
while the third permutation relation is
\begin{align}\label{eq-permrel2}
\xy(0,0)*{\includegraphics[scale=0.3]{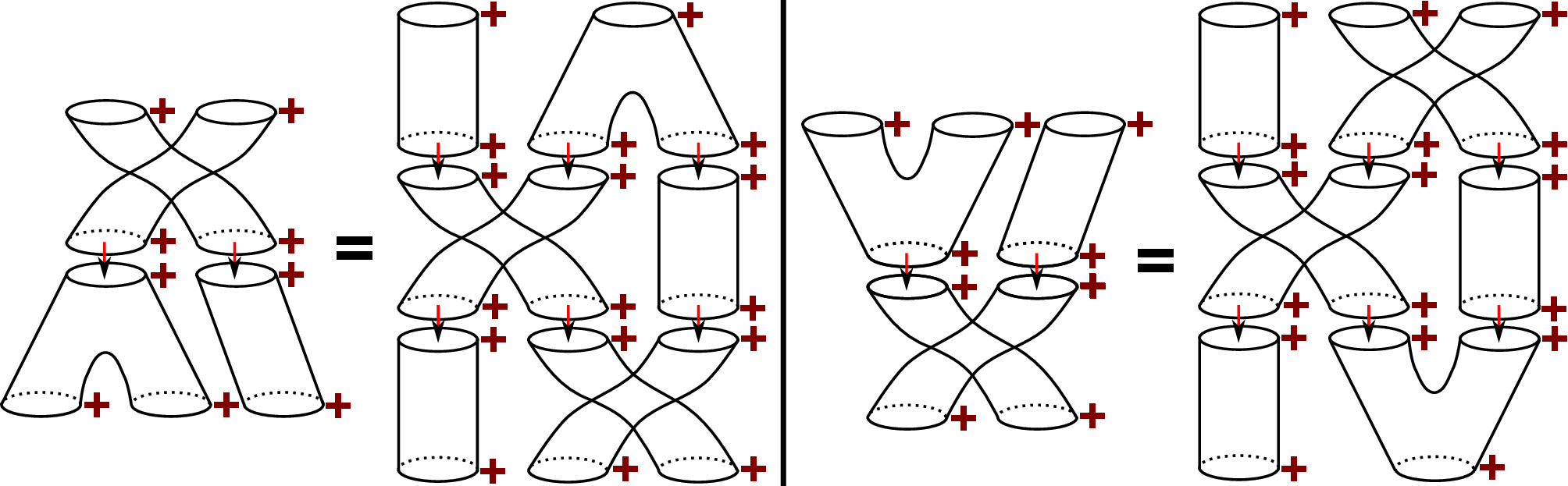}};\endxy
\end{align}
The important \textit{Frobenius, torus and M\"obius} relations are
\begin{align}\label{eq-frobandco}
\xy(0,0)*{\includegraphics[scale=0.25]{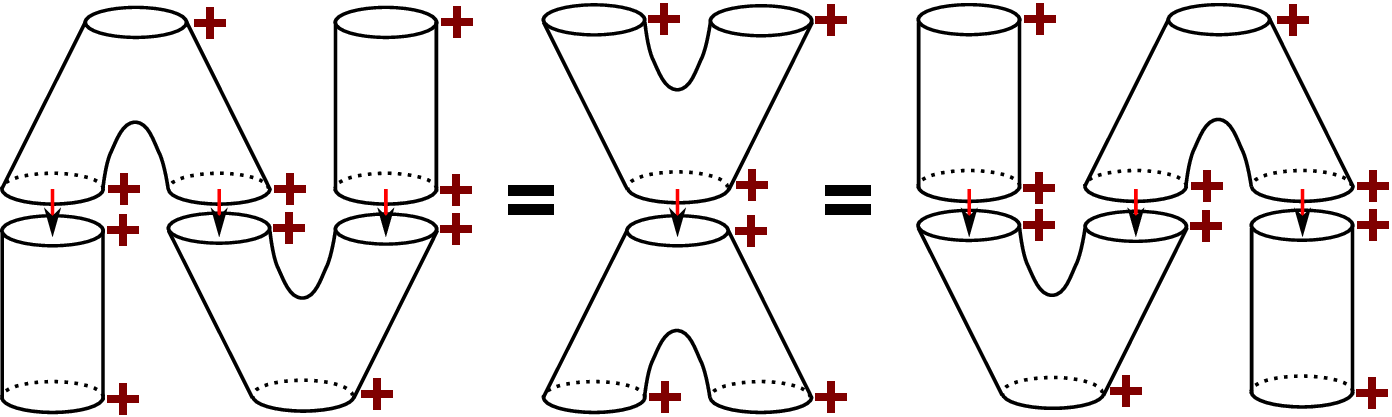}};\endxy\;\;\;\;\;\;\;\;\xy(0,0)*{\includegraphics[scale=0.3]{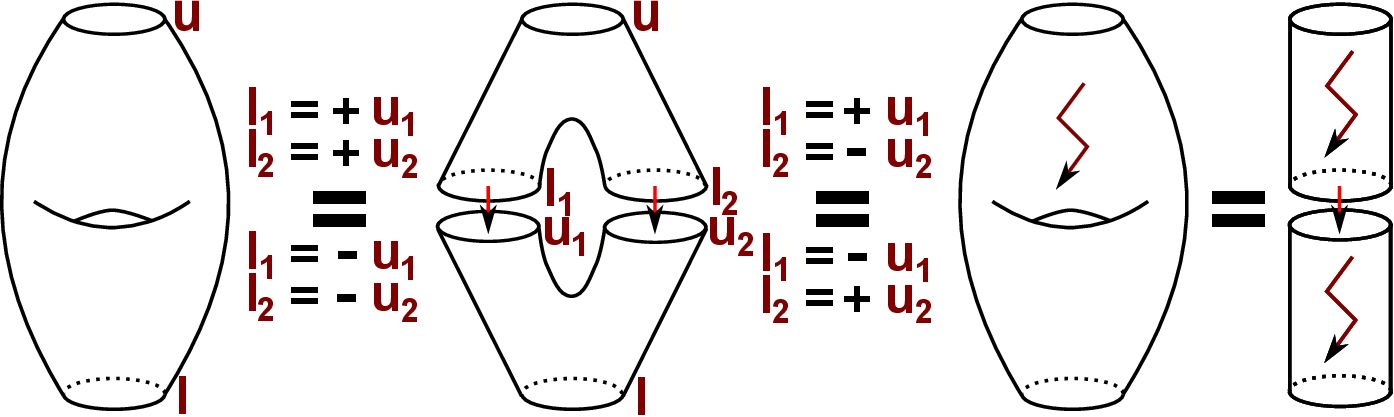}};\endxy
\end{align}
An $u$ or an $l$ mean an arbitrary glueing number and $-u,-l$ are the glueing numbers $u$ or $l$ multiplied by $-1$. Furthermore, the bolt represent a non-orientable surfaces and not illustrated parts are arbitrary.

It follows from these relations, that the cobordism $\amalg_{i\in I}\mathrm{id}^+_+$ is the identity morphism between $|I|$ v-circles. The cobordism $\Phi^-_+$ changes the boundary decoration of a morphism. Hence, the category above contains all possibilities for the decorations of the boundary components.

The category $\ucob_R(\emptyset)^*$ is the same as above, but without all minus signs in the relations (we mean ``honest'' minus signs, i.e. the minus-decorations are still in use).

Both categories are strict monoidal categories, since we are working with isomorphism classes of cobordisms. The monoidal structure is induced by the disjoint union $\amalg$. Moreover, both categories are symmetric. Note that they can be seen as $2$-categories, but it is more convenient to see them as monoidal $1$-category.
\end{defn}
\vspace*{0.5cm}

It is worth noting that the rest of this subsection can also be done for the category $\ucob_R(\emptyset)^*$ by dropping all the corresponding minus signs.

As for example in Definition 3.2 in~\cite{bn2}, we define the category $\mat(\mathcal C)$ to be the \textit{category of formal matrices} over a pre-additive category $\mathcal C$, i.e. the objects $\Ob(\mat(\mathcal C))$ are ordered, formal direct sums of the objects $\Ob(\mathcal C)$ and the morphisms $\Mor(\mat(\mathcal C))$ are matrices of morphisms $\Mor(\mathcal C)$. The composition is defined by the standard matrix multiplication. This category is sometimes called the \textit{additive closure} of the pre-additive category $\mathcal C$.

Furthermore, we define the category $\Kom_b(\mathcal C)$ to be the \textit{category of formal, bounded chain complexes} over a pre-additive category $\mathcal C$. Denote the category modulo formal chain homotopy by $\Kom_b(\mathcal C)^h$. More about such categories is collected in Sec.~\ref{sec-techhomalg}.
\vspace*{0.5cm}

Furthermore, we define $\ucob_R(\emptyset)^l$, which has the same objects as the category $\ucob_R(\emptyset)$, but morphisms modulo the local relations from Fig.~\ref{figureintroa-4}. We make the following definition.
\begin{defn}\label{defn-category2}
We denote by $\ukob_R$ the category $\Kom_b(\mat(\ucob_R(\emptyset)))$. Here our objects are formal, bounded chain complexes in the additive closure of the category of (possibly non-orientable) cobordisms with boundary decorations. We define $\ukob_R^h$ to be the category $\ukob_R$ modulo formal chain homotopy. Furthermore, we define $\ukob_R^l$ and $\ukob_R^{hl}$ in the obvious way. The notations $\ucob_R(\emptyset)^{(l)}$ or $\ukob_R^{(h)(l)}$ mean that we consider all possible cases, namely with or without an $h$ and with or without an $l$.
\end{defn}
\vspace*{0.5cm}

One effective way of calculation in $\ucob_R(\emptyset)$ is the usage of the \textit{Euler characteristic}\footnote{Here we consider our morphisms as surfaces.}. It is well-known that the Euler characteristic is invariant under homotopies and that it satisfies 
\[
\chi(C_2\circ C_1)=\chi(C_1)+\chi(C_2)-\chi(\mathcal O_2)\;\;\text{ and }\;\;\chi(C_1\amalg C_2)=\chi(C_1)+\chi(C_2)
\]
for any two cobordisms $C_1\colon\mathcal O_1\to\mathcal O_2$ and $C_2\colon\mathcal O_2\to\mathcal O_3$. Because the objects of $\ucob_R(\emptyset)$ are disjoint unions of v-circles, we have the following lemmata.
\begin{lem}\label{cor-euler1}
The Euler characteristic satisfies $\chi(C_1\circ C_2)=\chi(C_1)+\chi(C_2)$ for all morphisms $C_1,C_2$ of the category $\ucob_R(\emptyset)$.
\end{lem}
\begin{lem}\label{cor-euler2}
The generators of the category $\ucob(\emptyset)_R$ satisfy $\chi(\mathrm{id}^+_+)=\chi(\mathrm{id}^-_-)=0$ and $\chi(\Phi^-_+)=\chi(\Phi^+_-)=0$ and $\chi(\Delta^+_{++})=\chi(m^{++}_+)=\chi(\theta)=-1$. The composition of a cobordism $C$ with $\mathrm{id}^+_+$ or $\Phi^-_+$ does not change $\chi(C)$.
\end{lem}
\begin{proof}
Both Lemmata are well-known (we hope) statements. They can be found in various textbooks and we skip to make a recommendation to any.
\end{proof}
It is worth noting that the Lemmata~\ref{cor-euler1} and~\ref{cor-euler2} ensure that the category $\ucob_R(\emptyset)$ can be seen as a \textit{graded category}, that is the grading of morphisms is the Euler characteristic. Recall that a saddle between v-circles is a saddle inside a certain neighbourhood and the identity outside of it.
\begin{lem}\label{cor-saddleclass}
Every saddle is homeomorphic to one of the following three cobordisms (and some extra cylinders for not affected components). Hence, after decorating the boundary components, we get nine different possibilities, if we fix the decorations of the cylinders to be $+$.
\begin{itemize}
\item[(a)] A two times punctured projective plane $\theta=\mathbb{RP}^2_2$ iff the saddle has two boundary circles.
\item[(b)] A pantsup-morphism $m$ iff the saddle is a cobordism from two circles to one circle.
\item[(c)] A pantsdown-morphism $\Delta$ iff the saddle is a cobordism from one circle to two circles.
\end{itemize}
\end{lem}
\begin{proof}
We note that an open saddle $S$ has $\chi(S)=-1$. Hence, after closing its boundary components, we get the statement.
\end{proof}
\vspace*{0.5cm}

Now we deduce some basic relations between the basic cobordisms. Afterwards, we prove a proposition which is a key point for the understanding of the problematic face from~\ref{probcube}. Note the difference between the relations (b),(c) and (d),(e). Moreover, (k) and (l) are also very important.
\begin{lem}\label{lem-basiscalculations}
The following rules hold.
\begin{itemize}
\item[(a)] $\Phi^-_+\circ\Phi^-_+=\mathrm{id}^+_+\circ\mathrm{id}^+_+=\mathrm{id}^+_+$, $\tau^{++}_{++}\circ\tau^{++}_{++}=\mathrm{id}^{++}_{++}$.
\item[(b)] $(\Phi^-_+\amalg\Phi^-_+)\circ\Delta^+_{++}=\Delta^+_{--}=-\Delta^-_{++}=-\Delta^+_{++}\circ\Phi^-_+$.
\item[(c)] $(\Phi^-_+\amalg\mathrm{id}^+_+)\circ\Delta^+_{++}=\Delta^+_{-+}=-\Delta^-_{+-}=-(\mathrm{id}^+_+\amalg\,\Phi^-_+)\circ\Delta^+_{++}\circ\Phi^-_+$.
\item[(d)] $m^{++}_+\circ(\Phi^-_+\amalg\Phi^-_+)=m^{--}_+=m^{++}_-=\Phi^-_+\circ m^{++}_+$.
\item[(e)] $m^{++}_+\circ(\Phi^-_+\amalg\mathrm{id}^+_+)=m^{-+}_+=m^{+-}_-=\Phi^+_+\circ m^{++}_+\circ(\mathrm{id}^+_+\amalg\,\Phi^-_+)$.
\item[(f)] $m^{++}_+\circ\Delta^+_{++}=(\mathrm{id}^+_+\amalg\Delta^+_{++})\circ(m^{++}_+\amalg\mathrm{id}^+_+)$ (Frobenius relation).
\item[(g)] $m^{++}_{+}\circ(m^{++}_+\amalg\,\mathrm{id}^+_+) =m^{++}_{+}\circ(\mathrm{id}^+_+\amalg\,m^{++}_+)$ (associativity relation).
\item[(h)] $(\Delta^+_{++}\amalg\,\mathrm{id}^+_+)\circ\Delta^+_{++}=(\mathrm{id}^+_+\amalg\,\Delta^+_{++})\circ\Delta^+_{++}$ (associativity relation).
\item[(i)] $m^{++}_+\circ\tau^{++}_{++}\circ(\Phi^-_+\amalg\mathrm{id}^+_+)=m^{+-}_+$ (first permutation $\Phi$ relation).
\item[(j)] $(\Phi^-_+\amalg\mathrm{id}^+_+)\circ\tau^{++}_{++}\circ\Delta^+_{++}=\Delta^+_{+-}$ (second permutation $\Phi$ relation).
\item[(k)] $\theta\circ\Phi^-_+=\Phi^-_+\circ\theta=\theta$, $\theta=-\theta$ ($\theta$ relations).
\item[(l)] $\mathcal K=\theta^2$. Here $\mathcal K$ is a two times punctured Klein bottle.
\end{itemize}
\end{lem}
\begin{proof}
Most of the equations follow directly from the relations in Definition~\ref{defn-category} above. The rest are easy to check and their proofs are therefore omitted.
\end{proof}
The following example illustrates that some cobordisms are in fact isomorphisms.
\begin{ex}\label{ex-vrmremoval}
The two cylinders $\mathrm{id}^+_+,\Phi^-_+$ are the only isomorphisms between two equal objects. Let us denote $\mathcal O_1$ and $\mathcal O_2$ two objects which differ only through a finite sequence of the virtual Reidemeister moves. The vRM-cobordisms from Fig.~\ref{figure0-reide} induce isomorphisms $C\colon\mathcal O_1\to\mathcal O_2$. To see this we mention that the three cobordisms are isomorphisms, i.e. their inverses are the cobordisms which we obtain by turning the pictures upside down. Then use statement (a) of Lemma~\ref{lem-basiscalculations}.
\end{ex}
\vspace*{0.5cm}
\begin{prop}\label{prop-nonorientablefaces}\textbf{(Non-orientable faces)}
Let $\Delta^{u}_{l_1l_2}$ and $m^{u^{\prime}_1u^{\prime}_2}_{l^{\prime}}$ be the surfaces from Fig.~\ref{figure1-1}. Then the following are equivalent.
\begin{itemize}
 \item[(a)] $m^{u^{\prime}_1u^{\prime}_2}_{l^{\prime}}\circ\Delta^{u}_{l_1l_2}=\mathcal K$. Here $\mathcal K$ is a two times punctured Klein bottle.
 \item[(b)] $l_1=u^{\prime}_1$ and $l_2=-u^{\prime}_2$ or $l_1=-u^{\prime}_1$ and $l_2=u^{\prime}_2$.
\end{itemize}
Otherwise $m^{u^{\prime}_1u^{\prime}_2}_{l^{\prime}}\circ\Delta^{u}_{l_1l_2}$ is a two times punctured torus $\mathcal T$. We call this the \textit{M\"obius relation}.
\end{prop}
\begin{proof}
Let us call $C$ the composition $C=m^{u^{\prime}_1u^{\prime}_2}_{l^{\prime}_1}\circ\Delta^{u_1}_{l_1l_2}$. A quick computation shows $\chi(C)=-2$. Because $C$ has two boundary components, $C$ is either a 2-times punctured torus or a 2-times punctured Klein bottle and the statement follows from the torus and M\"obius relations in~\ref{eq-frobandco}.
\end{proof}
\subsection{The topological category for v-tangles}
In this part of Sec.~\ref{sec-vkhcat} we extend the notions above in such a way that they can be used for v-tangles as well. As explained in Sec.~\ref{sec-vkhsum}, the most important difference is the usage of an extra decoration which we call the \textit{indicator}. The rest is (almost) the same as above. Again all definitions and statements can be done for an analogue of the category $\ucob_R(\emptyset)^*$. First we define/recall the notion of a \textit{virtual tangle (diagram)}, called \textit{v-tangle (diagram)}.
\begin{defn}\label{defn-vtangle}(\textbf{Virtual tangles}) A \textit{virtual tangle diagram with $k\in\bN$ boundary points} $T^k_D$ is a planar graph embedded in a disk $D^2$. This planar graph is a collection of \textit{usual vertices} and $k$-\textit{boundary vertices}. We also allow circles, i.e. closed edges without any vertices.
\vspace*{0.5cm}

The usual vertices are all of valency four. Any of these vertices is either an overcrossing $\slashoverback$ or an undercrossing $\backoverslash$ or a virtual crossing $\virtual$. Latter is marked with a circle. The boundary vertices are of valency one and are part of the boundary of $D^2$.

As before, we call the crossings $\slashoverback$ and $\backoverslash$ \textit{classical crossings} or just \textit{crossings} and a virtual tangle diagram without virtual crossings a \textit{classical tangle diagram}.

A \textit{virtual tangle with $k\in\bN$ boundary points} $T^k$ is an equivalence class of virtual tangle diagrams $T^k_D$ module boundary preserving isotopies and \textit{generalised Reidemeister moves}.

We call a virtual tangle $T^k$ \textit{classical} if the set $T^k$ contains a classical tangle diagram. A v-string is a string starting and ending at the boundary without classical crossings. Moreover, we call a v-circle/v-string without virtual crossings a \textit{c-circle/c-string}.
\vspace*{0.5cm}

The \textit{closure of a v-tangle diagram with *-marker} $\mathrm{Cl}(T^k_D)$ is a v-link diagram which is constructed by capping of neighbouring boundary points (starting from a fixed point marked with the *-marker and going counterclockwise) without creating new virtual crossings. For an example see Fig.~\ref{figure1-2}.

There are exactly two, sometimes not equivalent, closures of any v-tangle diagram. In the figure below the two closures are pictured using green, dashed edges.
\begin{figure}[ht]
     \centerline{\includegraphics[scale=0.45]{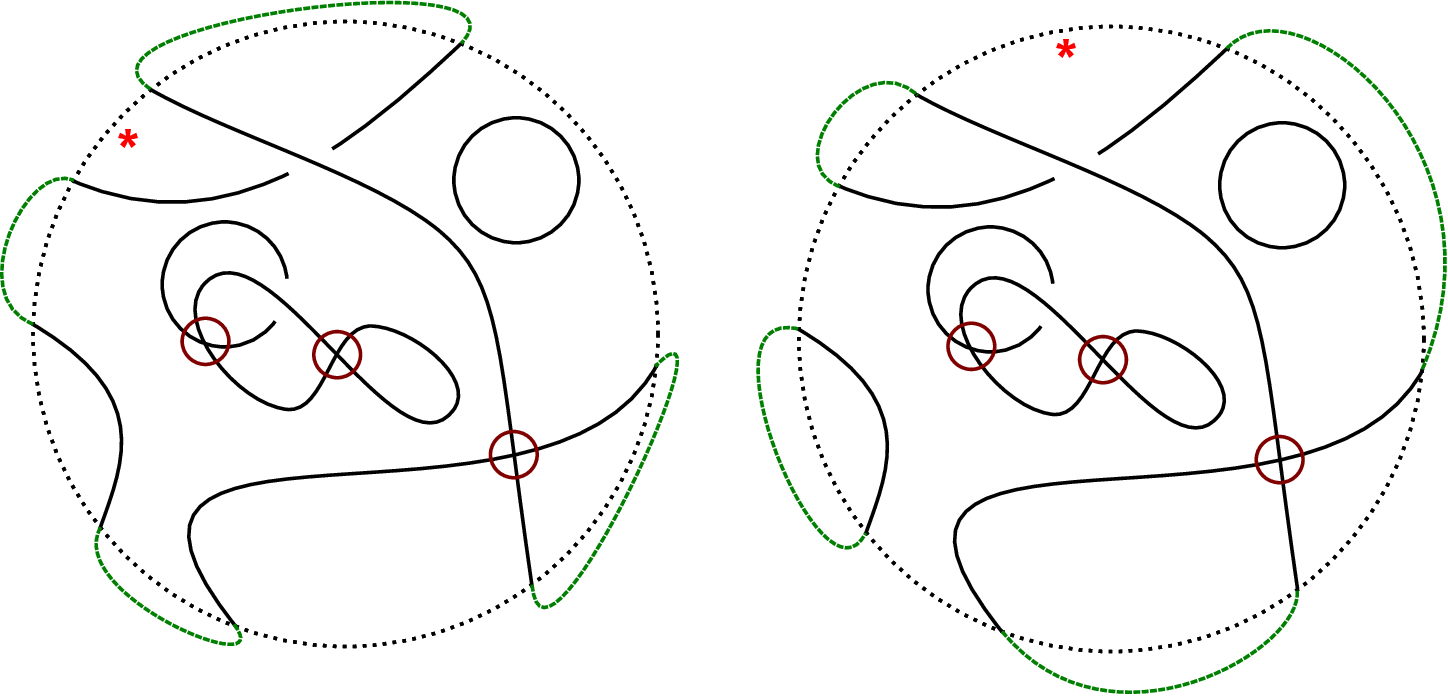}}
  \caption{A *-marked v-tangle and two different closures.}
  \label{figure1-2}
\end{figure}

The notions of an \textit{oriented} virtual tangle diagram and of an \textit{oriented} virtual tangle are defined analogously (see also Sec.~\ref{sec-introa}). The latter modulo \textit{oriented} generalised Reidemeister moves and boundary preserving isotopies. From now on every v-tangle (diagram) is oriented. But we suppress this notion to avoid confusion with other (more important) notations.
\end{defn}
We define the category of \textit{open cobordisms with boundary decorations}. It is almost the same as in Definition~\ref{defn-category}, but the corresponding cobordisms are allowed to be open, i.e. they could have vertical boundary components, and are decorated with an extra information: A number in the set $\{0,+1,-1\}$ (exactly one, even for non-connected cobordisms). We picture the number $0$ as a bolt.
\begin{defn}\label{defn-category3}(\textbf{The category of open cobordisms with boundary decorations}) Let $k\in\bN$ and let $R$ be a commutative and unital ring. The category $\ucob_R(k)$ should be $R-$pre-additive. The symbol $\amalg$ denotes the disjoint union.

\textbf{The objects:}

The \textit{objects} \textit{$\Ob(\ucob_R(k))$} are numbered (all components are labelled with numbers) v-tangle diagrams with $k$ boundary points without classical crossings. Objects are denoted by $\mathcal O=\coprod_{i\in I}\mathcal O_i$. Here $\mathcal O_i$ are the v-circles/v-strings and $I$ is a finite, ordered index set. The objects of the category are equivalence (modulo \textit{boundary preserving, planar isotopies}) classes of four-valent graphs.

\textbf{The generators:}

The \textit{generators} of $\Mor(\ucob_R(k))$ are the cobordisms in Fig.~\ref{figure1-3}. The cobordisms pictured are all between c-circles or c-strings. As before, we do not picture all the other possibilities, but we include them in the list of generators. 
\begin{figure}[ht]
	 \centerline{\xy
     (0,0)*{\includegraphics[scale=0.5]{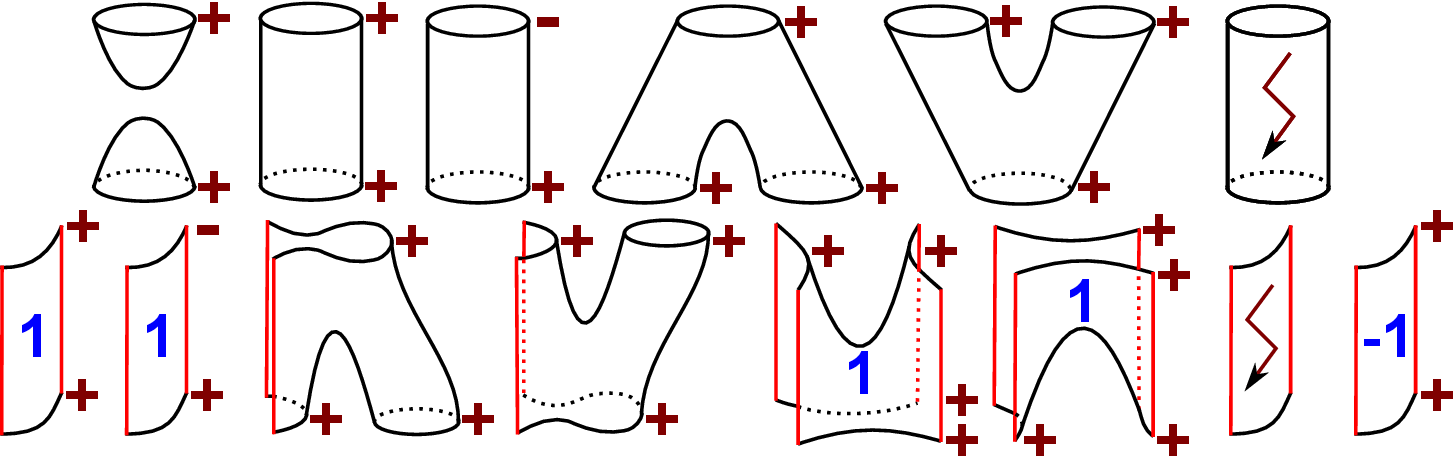}};
     (-55,7.8)*{\varepsilon^+};
     (-55,13.8)*{\iota_+};
     (-34,21.5)*{\mathrm{id}^+_+};
     (-20.5,21.5)*{\Phi^-_+};
     (27,21.5)*{m^{++}_+};
     (2.5,21.5)*{\Delta_{++}^+};
     (46.8,21.5)*{\theta};
     (-58.5,-20)*{\mathrm{id}(1)^+_+};
     (-47.5,-20)*{\Phi(1)^-_+};
     (-31,-20)*{S^{+}_{++}};
     (-10,-20)*{S^{++}_{+}};
     (31,-20)*{S^{++}_{++}};
     (12.5,-20)*{S^{++}_{++}};
     (45.5,-20)*{\theta};
     (58.5,-20)*{\mathrm{id}(-1)^+_+};
     \endxy}
  \caption{The generators for the set of morphisms.}
  \label{figure1-3}
\end{figure}

Every generator has a decoration from the set $\{0,+1,-1\}$. We call this decoration the \textit{indicator} of the cobordism. If no indicator is pictured, then it is $+1$. Indicators behave multiplicatively.

Every generator with a decoration $\{+1,-1\}$ has extra decorations from the set $\{+,-\}$ at every horizontal boundary component. We call these decorations the \textit{glueing numbers} of the cobordism. The vertical boundary components are pictured in red.

We consider these cobordisms up to boundary preserving homeomorphisms (as abstract surfaces). Hence, between circles or strings with v-crossings the generators are the same up to boundary preserving homeomorphisms, but immersed into $D^2\times[-1,1]$.

We denote the different generators (from left to right; top row first) by $\iota_+$ and $\varepsilon^+$, $\mathrm{id}^+_+$ and $\Phi^-_+$, $\Delta^+_{++}$, $m^{++}_+$ and $\theta$, $\mathrm{id}(1)^+_+$ and $\Phi(1)^-_+$, $S^{+}_{++}$ and $S^{++}_{+}$, $S(1)^{++}_{++}$, $\theta$ and $\mathrm{id}(-1)^+_+$.

The composition of the generators formally needs again internal decorations to remember how they were glued together. But again we suppress them and hope the reader does not get confused (at least not more than the author). Moreover, as before, cobordisms with a $0$-indicator do not have any boundary decorations, i.e. they are deleted after glueing.

\textbf{The morphisms:}

The \textit{morphisms} \textit{$\Mor(\ucob_R(k))$} are cobordisms between the objects in the following way. We identify the collection of numbered v-circles/v-strings with circles/strings immersed into $D^2$.

Given two objects $\mathcal O_1,\mathcal O_2$ with $k_1,k_2\in\bN$ numbered v-circles or v-strings, then a morphism $\mathcal C\colon\mathcal O_1\to\mathcal O_2$ is a surface immersed in $D^2\times[-1,1]$ whose non-vertical boundary lies only in $D^2\times\{-1,1\}$ and is the disjoint union of the $k_1$ numbered v-circles or v-strings from $\mathcal O_1$ in $D^2\times\{1\}$ and the disjoint union of the $k_2$ numbered v-circles or v-strings from $\mathcal O_2$ in $D^2\times\{-1\}$. The morphisms are generated (as abstract surfaces) by the generators from above (see Fig.~\ref{figure1-3}).

\textbf{The decorations:}

Every morphism has an \textit{indicator} from the set $\{0,+1,-1\}$.

Moreover, every morphism $C\colon\mathcal O_1\to\mathcal O_2$ in $\Mor(\ucob_R(k))$ is a cobordism between the numbered v-circles or v-strings of $\mathcal O_1$ and $\mathcal O_2$. Let us say that the v-circles or v-strings of $\mathcal O_1$ are numbered $i\in\{1,\dots,l_1\}$ and the v-circles or v-strings of $\mathcal O_2$ are numbered for $i\in\{l_1+1,\dots,l_2\}$.

Every cobordism with $+1,-1$ as an indicator has a decoration on the $i$-th boundary circle. This decoration is an element of the set $\{+,-\}$. We call the decoration of the $i$-th boundary component the \textit{$i$-th glueing number} of the cobordism.

Hence, the morphisms of the category are pairs $(C,w)$. Here $C\colon\mathcal O_1\to\mathcal O_2$ is a cobordism from $\mathcal O_1$ to $\mathcal O_2$ immersed into $D^2\times[-1,1]$ and $w$ is a string of length $l_2$ in such a way that the $i$-th letter of $w$ is the $i$-th glueing number of the cobordism and the last letter is the indicator or $w=0$ if the cobordism has $0$ as an indicator.

\textbf{Shorthand notation:}

We denote a morphism $C$ with an indicator from $\{+1,-1\}$ which is a connected surface by $C^{u}_{l}(\mathrm{in})$. Here $u,l$ are words in the alphabet $\{+,-\}$ in such a way that the $i$-th character of $u$ (of $l$) is the glueing number of the $i$-th circle of the upper (of the lower) boundary. The number $\mathrm{in}$ is the indicator. The construction above ensures that this notation is always possible. Therefore, we denote an arbitrary morphism as before by ($C^{u_i}_{l_i}$ are its connected components and $u_i,l_i$ are words in $\{+,-\}$)
\[
C(\pm 1)=(C^{u_1}_{l_1}\amalg\cdots\amalg C^{u_{k}}_{l_{k}})(\pm 1).
\]
For a morphism with $0$ as indicator we do not need any boundary decorations. With a slight abuse of notation, we denote all these cobordisms as the non-orientable cobordism $\theta$.

\textbf{The relations:}

There are different relations between the cobordisms, namely \textit{topological relations} and \textit{combinatorial relations}. The latter relations are described by the glueing numbers and indicators of the cobordisms and the glueing of the cobordisms. The topological relations are not pictured but it should be clear how they should work. Moreover, we have only pictured the most important new relations below, but there should hold analogous relations as in Definition~\ref{defn-category}. The reader should read these relations in the same vein as before.

The most interesting new relations are the three combinatorial
\begin{align}\label{eq-combrel123}
\xy(0,0)*{\includegraphics[scale=0.25]{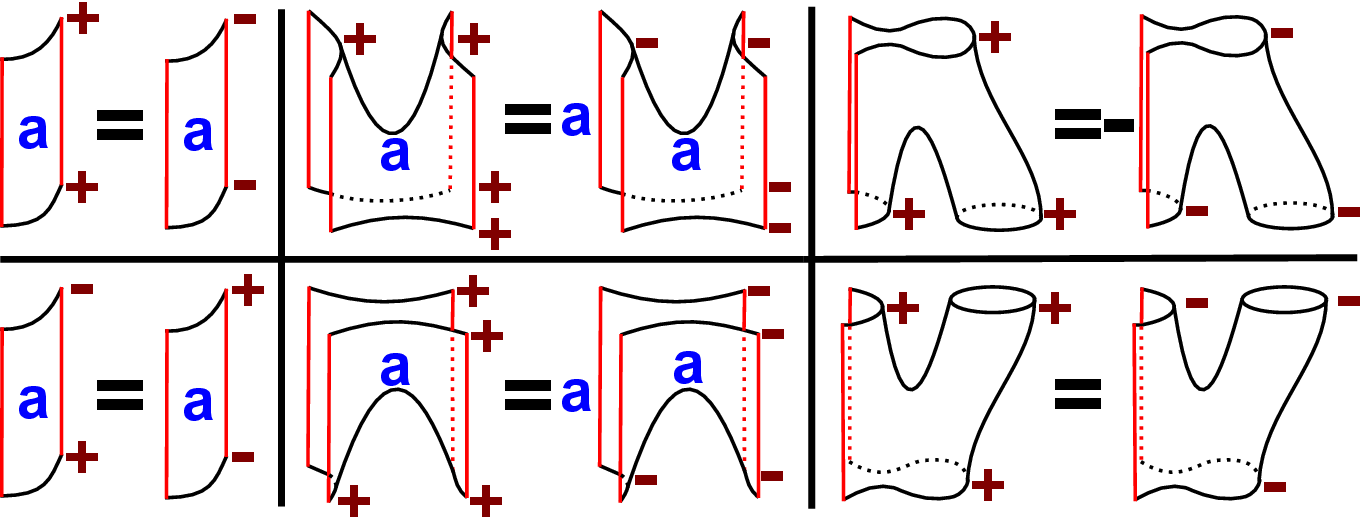}};\endxy\;\;\;\;\;\;\;\xy(0,0)*{\includegraphics[scale=0.25]{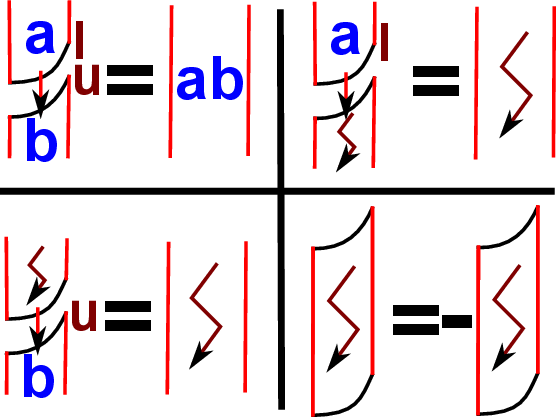}};\endxy\;\;\;\;\;\;\xy(0,0)*{\includegraphics[scale=0.25]{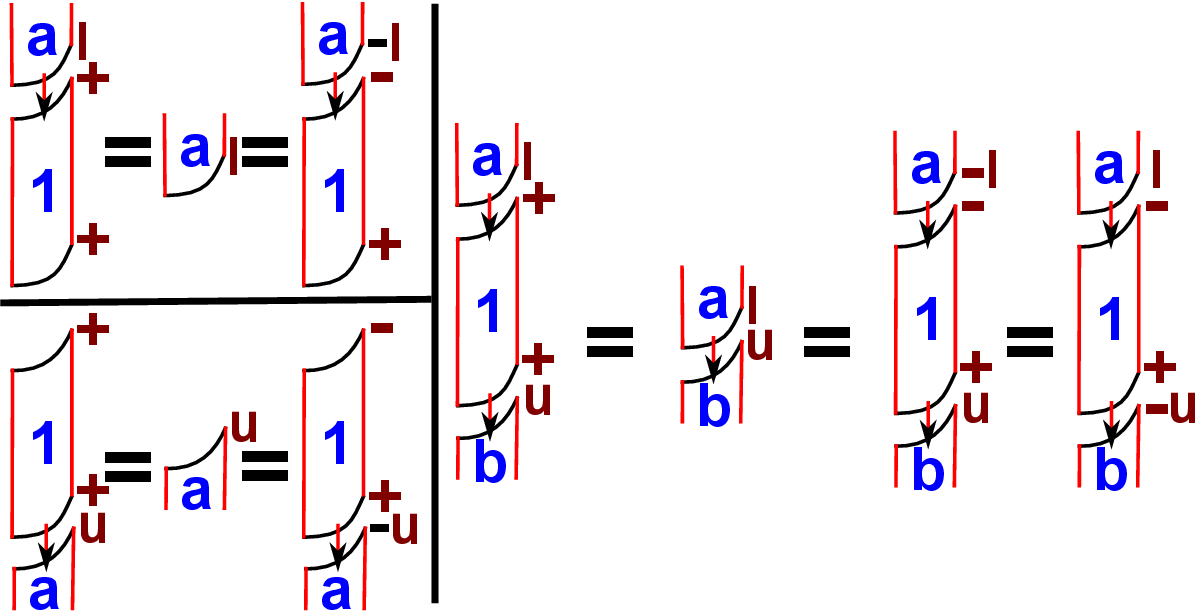}};\endxy
\end{align}
and the \textit{open} M\"obius relations (the glueing in these three cases is given by the glueing numbers, i.e. if there is an odd number of different glueing numbers, then the indicator is $0$ and just the product otherwise).
\begin{align}\label{eq-tanmoe}
\xy(0,0)*{\includegraphics[scale=0.325]{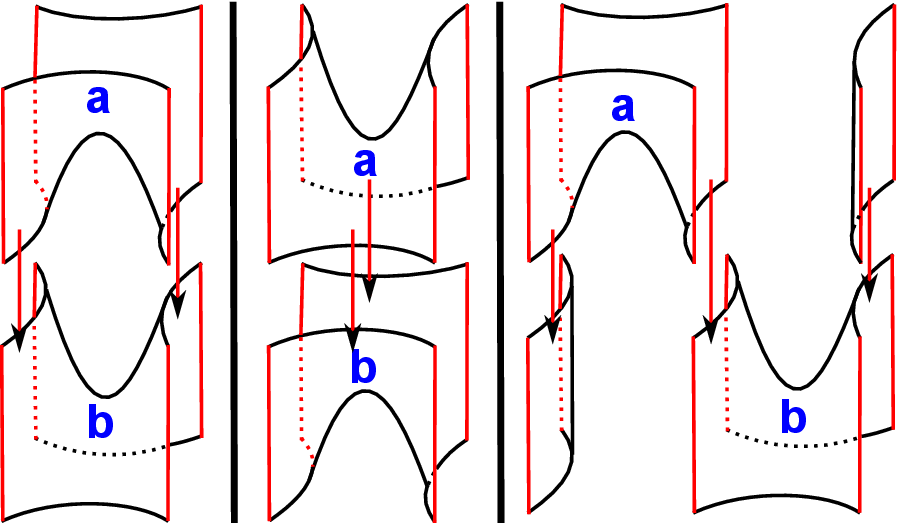}};\endxy
\end{align}
We define the category $\ucob_R(\omega)$ to be the category whose objects are $\bigcup_{k\in\bN}\Ob(\ucob_R(k))$ and whose morphisms are $\bigcup_{k\in\bN}\Mor(\ucob_R(k))$. Moreover, it should be clear how to convert Definition~\ref{defn-category2} to the open case. Note that this category is also graded, but the degree function has to be a little bit more complicated (since glueing with boundary behaves differently), that is the degree of a cobordism $C\colon\mathcal O_1\to\mathcal O_2$ is given by
\[
\mathrm{deg}(C)=\chi(C)-\frac{b}{2},\;\;\text{ where }b\text{ equals the number of vertical boundary components}.
\]
The reader should check that this definition makes the category graded, that is the degree of a composition is the sum of the degrees of its factors.
\end{defn}
\vspace*{0.5cm}

Note the following collection of formulas that follow from the relations. Recall that $\Phi^-_+$ and $\Phi(1)^-_+$ change the decorations and that $\theta$ and $\mathrm{id}(-1)^+_+$ change the indicators. With a slight abuse of notation, we omit to write $\amalg$ if it is not necessary, i.e. for the indicator changes. Moreover, since $\Phi^-_+$ and $\Phi(1)^-_+$ satisfy similar formulas, we only write down the equations for $\Phi^-_+$ and hope that it is clear how the others look.
\begin{lem}\label{lem-calc}
Let $\mathcal O,\mathcal O^{\prime}$ be two objects in $\ucob_R(k)$. Let $C\colon\mathcal O\to\mathcal O^{\prime}$ be a morphism that is connected, has $\mathrm{in}\in\{0,+1,-1\}$ as an indicator and $u$ and $l$ as decorated boundary strings. Then we have the following identities. We write $ C=C^u_l(\mathrm{in})$ as a shorthand notation if the indicators and glueing numbers do not matter. It is worth noting that the signs in (d) are important.
\begin{itemize}
\item[(a)] $C\circ\mathrm{id}(-1)^+_+=\mathrm{id}(-1)^+_+\circ C$ (indicator changes commute).
\item[(b)] $C\circ\theta=\theta\circ C$ ($\theta$ commutes).
\item[(c)] $C(0)\circ\Phi^-_+=\Phi^-_+\circ C(0)$ (first decoration commutation relation).
\item[(d)] Let $u^{\prime},l^{\prime}$ denote the decoration change at the corresponding positions of the words $u,l$. Then we have \begin{align*}
C(\pm 1)^u_l\circ(\mathrm{id}^+_+\amalg\dots\amalg\Phi^-_+\amalg\dots\amalg\mathrm{id}^+_+) &=C(\pm 1)^{u^{\prime}}_l=\pm C(\pm 1)^u_{l^{\prime}}\\ &=\pm(\Phi^-_+\amalg\dots\amalg\mathrm{id}^+_+\amalg\dots\amalg\Phi^-_+)\circ C(\pm 1)^u_l
\end{align*}(second decoration commutation relation).
\end{itemize}
\end{lem}
\begin{proof}
Everything follows by a straightforward (really!) usage of the relations in Definition~\ref{defn-category3}.
\end{proof}
\section{The topological complex for virtual links}\label{sec-vkhcom}
We note that the present section splits into three part, i.e. we define the virtual Khovanov complex first and we show that it is an invariant of v-links that agrees with the classical Khovanov complex for c-links. We have collected the more technical points, e.g. it is not clear why Definition~\ref{defn-topcomplex} gives a well-defined chain complex independent of all involved choices, in the last part. The last part is rather technical and the reader may skip it on the first reading.
\subsection{The definition of the complex}
In the present subsection we define the topological complex which we call the \textit{virtual Khovanov complex} $\bn{L_D}$ of an oriented v-link diagram $L_D$. This complex is an element of our category $\ukob_R$.

By Lemma~\ref{cor-saddleclass} we know that every saddle cobordism $S$ is homeomorphic to $\theta$, $m$ or $\Delta$ (disjoint union with cylinders for all v-cycles not affected by the saddle). We need extra information for the last two cases. We call these extra information the \textit{sign of the saddle} and the \textit{decoration of the saddle} (see Definitions~\ref{defn-sign} and~\ref{defn-deco}). 
\begin{defn}\label{defn-sign}(\textbf{The sign of a saddle}) We always want to read off signs or decorations for crossings that look like $\slashoverback$, but for a crossing $c$ in a general position there are two ways to rotate $c$ until it looks like $\slashoverback$ (which we call the \textit{standard position}). Since the sign depends on the two possibilities (see bottom row of Fig.~\ref{figure-xmarker}), we choose an \textit{x-marker} as in Fig.~\ref{figure-xmarker} for every crossing of $L_D$ and rotate the crossing in such a way that the markers match. 
\begin{figure}[ht]
     \centerline{\includegraphics[scale=0.8]{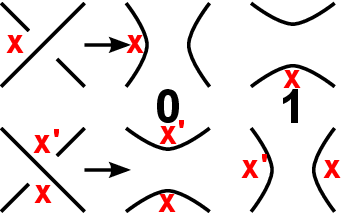}}
  \caption{Top: The x-marker for a crossing in the standard position. Bottom: Two possible choices (one denoted by $x^{\prime}$) for a crossing not in the standard position.}
  \label{figure-xmarker}
\end{figure}

We can say now that every orientable saddle $S$ can be viewed in a unique way as a formal symbol $S\colon\smoothing\rightarrow\hsmoothing$. Then the saddle $S$ carries an extra sign determined in the following way.
\begin{itemize}
\item Recall that the v-circles of any resolution are numbered. Moreover, the x-marker for the resolutions in the source and target of $S$ should be at the position indicated in the top row of Fig.~\ref{figure-xmarker}.
\item For a saddle $S\colon\gamma_a\to\gamma_{b}$ we denote the numbered v-circles of $\gamma_a,\gamma_{b}$ by $a_1,\dots, a_{k_a}$ and $b_1,\dots,b_{k_b}$ and the v-circles with the x-marker by $a^x_{i},b^x_{j}$.
\item Since the saddle $S$ is orientable, it either splits one v-circle or merges two v-circles. Hence, the two strings in the resolutions $\smoothing$ or $\hsmoothing$ are only different either in the target or in the source of $S$ and we denote the second affected v-circle by $b^y_{j^{\prime}}$ for a split and $a^y_{i^{\prime}}$ for a merge.
\item Then there exist two permutations $\sigma_1,\sigma_2$ from the fixed orderings for $S$, one for the source and one for the target, to other ordered sets $\{a_k\},\{b_k\}$ such that all $a_k\notin\{a^x_{i},a^y_{i^{\prime}}\}$ and all $b_{k^{\prime}}\notin\{b^x_{j},b^y_{j^{\prime}}\}$ are \textit{ordered ascending after} the (also ordered) $a^x_{i},a^y_{i^{\prime}}$ and $b^x_{j},b^y_{j^{\prime}}$.
\item Then we define the saddle sign $\mathrm{sgn}(S)$ by
\[
\mathrm{sgn}(S)=\mathrm{sgn}(\sigma_1)\cdot\mathrm{sgn}(\sigma_2).
\]
\end{itemize}
For completeness, we define the sign of a non-orientable saddle to be $0$. The \textit{sign $\mathrm{sgn}(F)$ of a face $F$} is then defined by the product of all the saddle signs of the saddles of $F$.
\end{defn}
\begin{ex}\label{ex-sign}
If we have a saddle $S$ between four v-circles numbered (fixed ordering) $u_1,u_2,u_3,u_4$ and three v-circles $l_1,l_2,l_3$ and the upper (in the target) x-marker is on the v-circle number $2$ and the lower (in the source) is on number $3$ and the second string of the upper part is number $1$, then the sign of $S$ is calculated by the product of the signs of the following two permutations.
\[
\sigma_1\colon(u_1,u_2,u_3,u_4)\mapsto (u_2,u_1,u_3,u_4)\;\;\text{and}\;\;\sigma_2\colon(l_1,l_2,l_3)\mapsto (l_3,l_1,l_2).
\]
Note that the saddle $S$ above ``multiplies $u_1,u_2$ to $l_3$'' and the x-marker is on $u_2$.
\end{ex}
Before we can define the virtual Khovanov complex we need to define the saddle decorations.
\begin{defn}\label{defn-deco}(\textbf{Saddle decorations})
By Lemma~\ref{cor-saddleclass} again, we only have to define the decorations in three different cases. First choose an x-marker as in Definition~\ref{defn-sign} for all crossings and choose orientations for the two resolutions $\gamma_a,\gamma_{a^{\prime}}$. We say the formal saddle of the form
\[
S^{++}_{++}\colon\du\to\ler
\]
is the \textit{standard oriented saddle}. Moreover, every saddle looks locally like the standard oriented saddle, but with possibly different orientations. Now we spread the decorations as follows.
\begin{itemize}
\item The non-orientable saddles do not get any extra decorations. It should be noted that locally non-alternating saddles, e.g. $S\colon \dd\to\lel$, are always non-orientable and vice versa.
\item The orientable saddles get a $+$ decoration at strings where the orientations agree and a $-$ where they disagree (after rotating it to the standard position defined above).
\item All cylinders of $S$ are $\mathrm{id}^+_+$ iff the corresponding unchanged v-circles of $\gamma_a$ and $\gamma_{a^{\prime}}$ have the same orientation and a $\Phi^-_+$ otherwise. 
\end{itemize}
To summarise we give the following table (we also give a way to denote the decorations for the saddles). We suppress the cylinders in the Table 1, but we note that the last point of the list above, i.e. the decorations of the cylinders, is important and can not be avoided in our context.

In the Table 1 below we write $m,\Delta$ for the corresponding saddles $S$.
\begin{table}[ht]
\begin{center}
\begin{tabular}{|c|c||c|c|}
\hline String & Comultiplication & String & Multiplication \\ 
\hline $\du\to\ler$ & $\Delta^{+}_{++}$ & $\du\to\ler$ & $m_{+}^{++}$ \\ 
\hline $\du\to\rir$ & $\Delta^+_{-+}=\Phi_1\circ\Delta^+_{++}$ & $\uu\to\ler$ & $m_{+}^{-+}=m_{+}^{++}\circ\Phi_1$ \\ 
\hline $\du\to\lel$ & $\Delta^+_{+-}=\Phi_2\circ\Delta^+_{++}$ & $\dd\to\ler$ & $m_{+}^{+-}\circ\Phi_2$ \\ 
\hline $\du\to\ril$ & $\Delta^+_{--}=\Phi_{12}\circ\Delta^+_{++}$ & $\ud\to\ler$ & $m_{+}^{--}\circ\Phi_{12}$ \\ 
\hline $\ud\to\ler$ & $\Delta^-_{++}=\Delta^+_{++}\circ\Phi^-_+$ & $\du\to\ril$ & $\Phi^-_+\circ m_{-}^{++}$ \\ 
\hline $\ud\to\rir$ & $\Delta^-_{-+}=\Phi_1\circ\Delta^+_{++}\circ\Phi^-_+$ & $\uu\to\ril$ & $\Phi^-_+\circ m_{-}^{-+}\circ\Phi_1$ \\ 
\hline $\ud\to\lel$ & $\Delta^-_{+-}=\Phi_2\circ\Delta^+_{++}\circ\Phi^-_+$ & $\dd\to\ril$ & $\Phi^-_+\circ m_{-}^{+-}\circ\Phi_2$ \\ 
\hline $\ud\to\ril$ & $\Delta^-_{--}=\Phi_{12}\circ\Delta^+_{++}\circ\Phi^-_+$ & $\ud\to\ril$ & $\Phi^-_+\circ m_{-}^{--}\circ\Phi_{12}$ \\ 
\hline 
\end{tabular}
\caption{The decorations are spread based on the local orientations.}\label{tab-deco}
\end{center}
\end{table}
\end{defn}
At this point we are finally able to define the \textit{virtual Khovanov complex}. We call this complex \textit{the topological complex}.
\begin{defn}\label{defn-topcomplex}(\textbf{The topological complex}) For a v-link diagram $L_D$ with $n$ ordered crossings we define \textit{the topological complex} $\bn{L_D}$ as follows. We choose an x-marker for every crossing.
\begin{itemize}
\item For $i=0,\dots,n$ the $i-n_-$ \textit{chain module} is the formal direct sum of all $\gamma_a$ of length $i$. We consider the resolutions as elements of $\Ob(\ucob_R(\emptyset))$.
\item There are only morphisms between the chain modules of length $i$ and $i+1$.
\item If two words $a,a^{\prime}$ differ only in exactly the $r$-th letter and $a_r=0$ and $a'_r=1$, then there is a morphism between $\gamma_a$ and $\gamma_{a'}$. Otherwise all morphisms between components of length $i$ and $i+1$ are zero.
\item This morphism $S$ is a \textit{saddle} between $\gamma_a$ and $\gamma_{a'}$.
\item We consider \textit{numbered} and \textit{oriented resolutions} (we choose a numbering/orientation) and the saddles carry the \textit{saddle sings} and \textit{decorations} from the Definitions~\ref{defn-sign} and~\ref{defn-deco}.
\item We consider the saddles $S$ as elements of $\Mor(\ucob_R(\emptyset))$ where we interpret the saddle signs as scalars in $R$ and the saddle decorations as the corresponding boundary decorations.
\end{itemize}
\end{defn}
\begin{rem}\label{rem-topcom}
At this point it is not clear why we can choose the numbering of the crossings, the numbering of the v-circles, the x-markers and the orientation of the resolutions. Furthermore, it is not clear why this complex is a well-defined chain complex.

But we show in Lemma~\ref{lem-commutativeindependence} that the complex is independent of these choices, i.e. if $\bn{L_D}_1$ and $\bn{L_D}_2$ are well-defined chain complexes with different choices, then they are equal up to chain isomorphisms. Moreover, we show in Theorem~\ref{theo-facescommute} and Corollary~\ref{cor-chaincomplex} that the complex is indeed a well-defined chain complex. Hence, we see that
\[
\bn{L_D}\in\Ob(\ukob_R).
\]
For an example see Fig.~\ref{figure0-big}. This figure shows the virtual Khovanov complex of a v-diagram of the unknot.
\end{rem}
\subsection{The invariance}
There is a way to represent the topological complex of a v-link diagram $L_D$ as a cone of two v-link diagrams $L^0_D,L^1_D$. Here one fixed crossing of $L_D$ is resolved $0$ in $L^0_D$ and $1$ in $L^1_D$. Note that the cone construction, as explained in Definition~\ref{defn-cone}, works in our setting.

It should be noted that there is a saddle between any two resolutions that are resolved equal at all the other crossings of $L^0_D$ and $L^1_D$. This saddle induces a chain map (as explained in the proof below) between the topological complex of $L^0_D$ and $L^1_D$. We denote this chain map by $\varphi\colon\bn{L^0_D}\to \bn{L^1_D}$.
\begin{lem}\label{lem-cone}
Let $L_D$ be a v-link diagram and let $c$ be a crossing of $L_D$. Let $L^0_D$ be the v-link where the crossing $c$ is resolved 0 and let $L^1_D$ be the v-link where the crossing $c$ is resolved 1. Then we have
\[
\bn{L_D}=\Gamma(\bn{L^0_D}\xrightarrow{\varphi}\bn{L^1_D}).
\]
\end{lem}
\begin{proof}
The proof is analogous to the proof for the classical Khovanov complex. The only new thing to prove is the fact that the map $\varphi$, which resolves the crossing, induces a chain map. This is true because we can take the induced orientation (from the orientations of the resolutions of $L^0_D$ and $L^1_D$) of the strings of $\varphi$. This gives us the glueing numbers for the morphisms of $\varphi$. Here we need the Lemma~\ref{lem-commutativeindependence} to ensure that all faces anti-commute.
\end{proof}
\begin{ex}\label{ex-cone}
Let $L_D$ be the v-diagram of the unknot from Fig.~\ref{figure0-big}. Then we have
\[
\bn{L_D}=\Gamma(\varphi\colon\bn{\jpg{12mm}{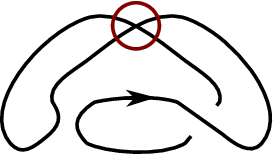}}\to\bn{\jpg{12mm}{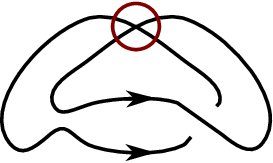}})=\Gamma(\varphi\colon L^0_D\to L^1_D).
\]
If we choose the orientation for the resolutions for the chain complexes $L^0_D,L^1_D$ to be the ones from Fig.~\ref{figure0-big}, then the map $\varphi$ is of the form $\varphi=(\theta,m^{--}_+)$.
\end{ex}
As a shorthand notation we only picture a certain part of a v-link diagram. The rest of the diagram can be arbitrary. Now we state the main theorem of this section.
\begin{thm}\label{thm-geoinvarianz}(\textbf{The topological complex is an invariant})
Let $L_D,L^{\prime}_D$ be two v-link diagrams which differ only through a finite sequence of isotopies and generalised Reidemeister moves. Then the complexes $\bn{L_D}$ and $\bn{L^{\prime}_D}$ are equal in $\ukob^{hl}_R$.
\end{thm}
\begin{proof}
We have to check invariance under the generalised Reidemeister moves from Fig.~\ref{figureintroa-1}. We follow the original proof of Bar-Natan in~\cite{bn2} (that is, his proof in Sec. 4 of~\cite{bn2}) with some differences. The main differences are the following.
\begin{itemize}
\item[(1)] We have to ensure that our cobordisms have the adequate decorations. For this we number the v-circles in a way such that the pictured v-circles have the lowest numbers and we use the orientations given below. It should be noted that Lemma~\ref{lem-commutativeindependence} ensures that we can use this numbering and these orientations without problems. We mention that we do not care about the saddle signs to maintain readability because they only affect the anti-commutativity of the faces. Hence, after adding some extra signs, the entire arguments work analogously.
\item[(2)] We have to check that the glueing of the cobordisms we give below works out correctly. This is a straightforward calculation using the relations in Lemma~\ref{lem-basiscalculations}.
\item[(3)] The proof of Bar-Natan uses the local properties of his construction. This is not so easy in our case. To avoid it we use some of the technical tools from homological algebra, i.e. Proposition~\ref{prop-sdr}.
\item[(4)] We have to check extra moves, i.e. the virtual Reidemeister moves vRM1, vRM2 and vRM3 and the mixed one mRM.
\end{itemize}
Recall that we have to use the Bar-Natan relations from Fig.~\ref{figureintroa-4} here. Note that the Bar-Natan relations do not contain any boundary components. Therefore, we do not need extra decorations for them. Because of this we can take the same chain maps as Bar-Natan (the cobordisms are the identity outside of the pictures). Furthermore, the whole construction is in $\ukob_R$.

The outline of the proof is as follows. For the RM1 and RM2 moves one has to show that the given maps induce chain homotopies, using the rules from Definition~\ref{defn-category} and Lemma~\ref{lem-basiscalculations} and the cone construction from Definition~\ref{defn-cone}. We note that we have to use Proposition~\ref{prop-sdr} to get the required statement for the RM1 and RM2 moves. Then the statement for the RM3 move follows with the cone construction from the RM2 move. The vRM1, vRM2 and vRM3 moves follow from their properties explained in Example~\ref{ex-vrmremoval}. Finally, the invariance under the mRM move can be obtained as an instance of Proposition~\ref{prop-sdr}.

We consider oriented v-link diagrams. Thus, there are a lot of cases to check. But all cases for the RM1 and RM2 moves are analogous to the cases shown below, i.e. one case for the RM1 move and three cases for the RM2 move. Note that the mirror images work similar.

The case for the RM1 move is pictured below. For the RM2 move we show that the virtual Khovanov complexes of
\[
\bn{\,\jpg{8mm}{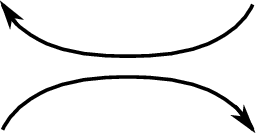}\,}\text{ and }\bn{\,\jpg{8mm}{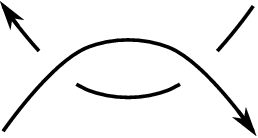}\,}\;\;\;\;\;\;\bn{\,\jpg{8mm}{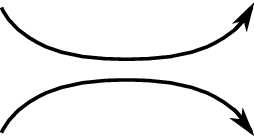}\,}\text{ and }\bn{\,\jpg{8mm}{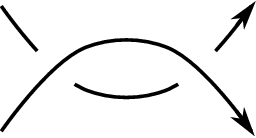}\,}
\]
are chain homotopic. Here both cases contain two different subcases. For the left case the upper left string can be connected to the upper right or to the lower left. For the other case the upper left string can be connected to the lower right or to the upper right. But the last case is analogous to the first. So we only consider the first three cases.

For the RM1 move (recall that the moves are pictured in Fig.~\ref{figureintroa-1}) we only have to resolve one crossing in the left picture and no crossing in the right. We choose the orientation in such a way that the saddle is a multiplication of the form $\ud\rightarrow\ril$. Thus, it is the multiplication $m^{--}_-=m^{++}_+$.

For the RM2 move we have to resolve two crossings in the left picture and no crossing in the right. For the first two cases we choose the orientation in such a way that the corresponding saddles are of the form $\du\rightarrow\ler$ for the left crossing and of the form $\ler\rightarrow\du$ for the right crossing. Hence, we only have $\Delta^+_{++}=-\Delta^-_{--}$ and $m^{--}_-=m^{++}_+$ saddles in the possible complexes.

For the third case we choose the orientation in such a way that the corresponding saddles are of the form $\uu\rightarrow\ril$ or $\uu\rightarrow\rir$ for the left crossing and of the form $\ler\rightarrow\du$ or $\rir\rightarrow\uu$ for the right crossing. Hence, we only have $m^{-+}_-$, $\theta$, $\Delta^-_{--}$ and $\theta$ saddles in the possible complexes.

We give the required chain maps $F,G$ and the homotopy $h$. Note our abuse of notation, that is, we denote the chain maps and homotopies and their parts with the same symbols. Moreover, the degree zero (we mean the homological degree) components are the leftmost non-trivial in the RM1 case and the middle non-trivial in the RM2 case.

One can prove that these maps are chain maps and that $F\circ G$ and $G\circ F$ are chain homotopic to the identity using the same arguments as Bar-Natan in Sec. 4 in~\cite{bn2} and the relations from Lemma~\ref{lem-basiscalculations}. We suppress the notation $\Gamma(\cdot)$ in the following. For the RM1 move we have
\[
\begin{xy}
  \xymatrix{
  \bn{\,\jpg{8mm}{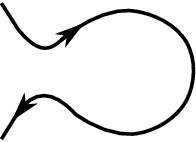}\,}: &  \bn{\,\jpg{8mm}{RM1-1}\,} \ar[rr]^0\ar@<2pt>[dd]^{F=\,\jpg{11mm}{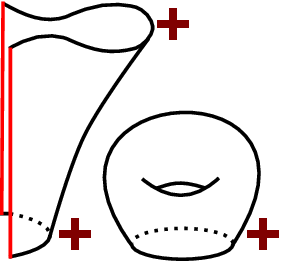}-\,\jpg{11mm}{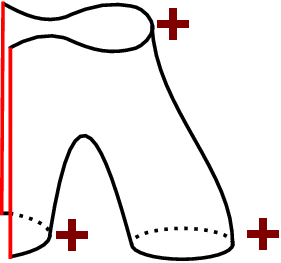}}   &  &   0\ar@<2pt>[dd]^0  \\
  & & & \\
  \bn{\,\jpg{8mm}{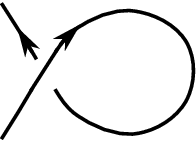}\,}: &  \bn{\,\jpg{8mm}{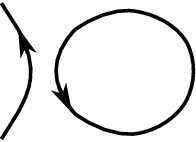}\,} \ar[rr]_{m^{++}_{+}}\ar@<2pt>[uu]^{G=\jpg{11mm}{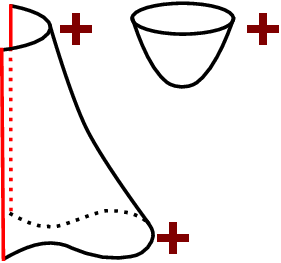}}           &  &   \bn{\,\jpg{8mm}{RM1-1}\,}\ar@<2pt>[uu]^0.   
  }
\end{xy}
\]
We also need to give an extra chain homotopy $h$. It is the one from below.
\[
h\colon\bn{\,\jpg{8mm}{RM1-1}\,}\to\bn{\,\jpg{8mm}{RM1-0}\,},\;\;h=-\,\jpg{11mm}{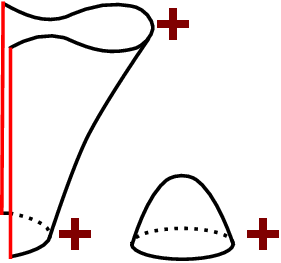}.
\]
An important observation is now that $G\circ F=\mathrm{id}$ and $h\circ F=0$. Beware our abuse of notation here, i.e. the parts of the homotopy $h$ and the chain map $F$ that can be composed are $0$. Thus, we are in the situation of Definition~\ref{defn-sdr} and can use Proposition~\ref{prop-sdr} to get
\[
\Gamma(\bn{\,\jpg{8mm}{RM1-1}\,})\simeq_h\Gamma(\bn{\,\jpg{8mm}{RM1}\,}).
\]
For the RM2 move the first two cases are
\[
\begin{xy}
  \xymatrix{
  \bn{\,\jpg{8mm}{RM2-1-a}\,}: &  0 \ar[rr]^0\ar@<2pt>[dd]^{0}   &  &   \bn{\,\jpg{8mm}{RM2-1-a}\,}\ar@<2pt>[dd]^{F=\begin{pmatrix}-\,\jpg{12mm}{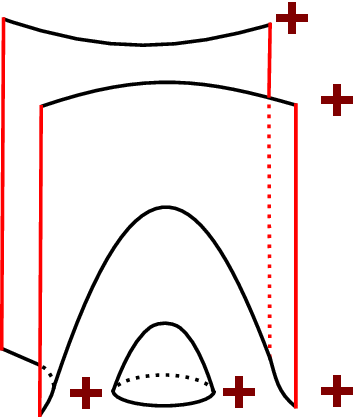} & \Phi^-_+\end{pmatrix}}\ar[rr]^0 & & 0\ar@<2pt>[dd]^{0} \\
  & & & & & \\
  \bn{\,\jpg{8mm}{RM2-a}\,}: &  \bn{\,\jpg{8mm}{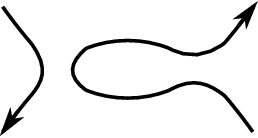}\,} \ar[rr]_/-1.0em/{d^{-1}}\ar@<2pt>[uu]^{0}           &  & \bn{\,\jpg{8mm}{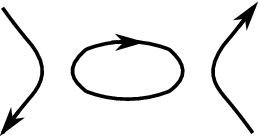}\,}\oplus\bn{\,\jpg{8mm}{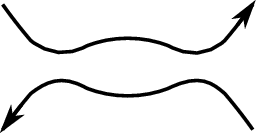}\,}\ar[rr]_/1.0em/{d^{0}}\ar@<2pt>[uu]^{G=\begin{pmatrix}\jpg{12mm}{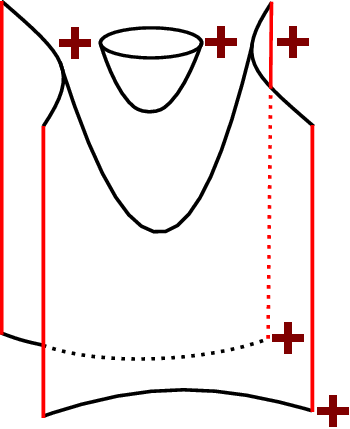} & \Phi^-_+ \end{pmatrix}^T} & & \bn{\,\jpg{8mm}{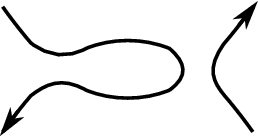}\,}\ar@<2pt>[uu]^0.   
  }
\end{xy}
\]
Here the differentials are either $d^{-1}=\begin{pmatrix}\Delta^-_{--} & \Delta^+_{++}\end{pmatrix}^T$ and $d^0=\begin{pmatrix}m^{++}_+ & m^{++}_+\end{pmatrix}$ in the second case or $d^{-1}=\begin{pmatrix}\mathrm{id}^+_+\amalg\,\Delta^-_{--} & m^{++}_+\end{pmatrix}^T$ and $d^0=\begin{pmatrix}m^{++}_+\amalg\,\mathrm{id}^+_+ & \Delta^-_{--}\end{pmatrix}$ in the first case. We can follow the proof of Bar-Natan again. Therefore, we need to give a chain homotopy. This chain homotopy is
\[
\begin{matrix}
h^{-1}\colon\bn{\,\jpg{8mm}{RM2-01-a}\,}\oplus\bn{\,\jpg{8mm}{RM2-10-a}\,}\to\bn{\,\jpg{8mm}{RM2-00-a}\,},\;h^{-1}=\begin{pmatrix}-\,\jpg{11mm}{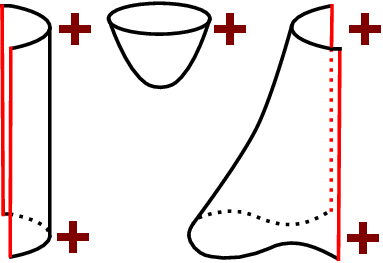} & 0 \end{pmatrix},\\h^0\colon\bn{\,\jpg{8mm}{RM2-11-a}\,}\to\bn{\,\jpg{8mm}{RM2-01-a}\,}\oplus\bn{\,\jpg{8mm}{RM2-10-a}\,},\;h^0=\begin{pmatrix}-\,\jpg{11mm}{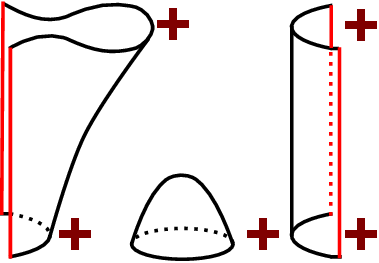} & 0 \end{pmatrix}^T.
\end{matrix}
\]
For the RM2 move the last case is
\[
\begin{xy}
  \xymatrix{
  \bn{\,\jpg{8mm}{RM2-1-b}\,}: &  0 \ar[rr]^0\ar@<2pt>[dd]^{0}   &  &   \bn{\,\jpg{8mm}{RM2-1-b}\,}\ar@<2pt>[dd]^{F=\begin{pmatrix}-\,\jpg{12mm}{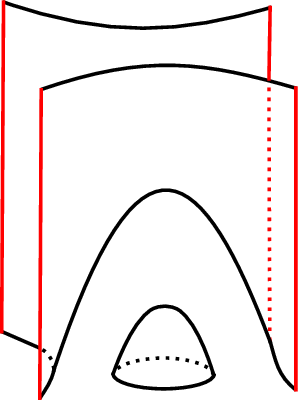} & \mathrm{id}^+_+\end{pmatrix}}\ar[rr]^0 & & 0\ar@<2pt>[dd]^{0} \\
  & & & & & \\
  \bn{\,\jpg{8mm}{RM2-b}\,}: &  \bn{\,\jpg{8mm}{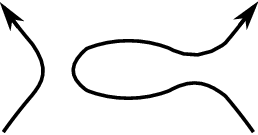}\,} \ar[rr]_/-1.0em/{d^{-1}}\ar@<2pt>[uu]^{0}           &  & \bn{\,\jpg{8mm}{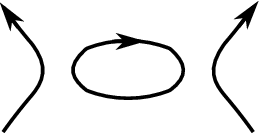}\,}\oplus\bn{\,\jpg{8mm}{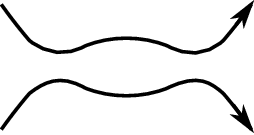}\,}\ar[rr]_/1.0em/{d^{0}}\ar@<2pt>[uu]^{G=\begin{pmatrix}\jpg{12mm}{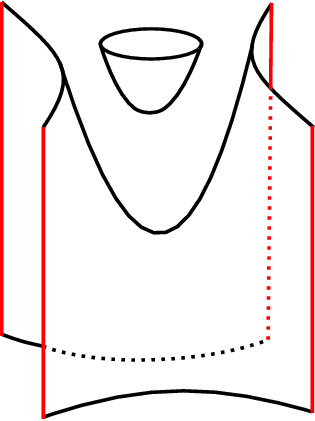} & \mathrm{id}^+_+ \end{pmatrix}^T} & & \bn{\,\jpg{8mm}{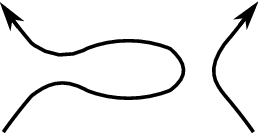}\,}\ar@<2pt>[uu]^0.   
  }
\end{xy}
\]
Here the differentials are either $d^{-1}=\begin{pmatrix}\Delta^-_{--} & \theta\end{pmatrix}^T$ and $d^0=\begin{pmatrix}m^{-+}_- & -\theta\end{pmatrix}$. Furthermore, saddles of the maps $F,G$ are also $\theta$ saddles. Hence, we do not need any decorations for them. The chain homotopy is defined by
\[
\begin{matrix}
h^{-1}\colon\bn{\,\jpg{8mm}{RM2-01-b}\,}\oplus\bn{\,\jpg{8mm}{RM2-10-b}\,}\to\bn{\,\jpg{8mm}{RM2-00-b}\,},\;h^{-1}=\begin{pmatrix}-\,\jpg{11mm}{RM2-h1} & 0 \end{pmatrix},\\h^0\colon\bn{\,\jpg{8mm}{RM2-11-b}\,}\to\bn{\,\jpg{8mm}{RM2-01-b}\,}\oplus\bn{\,\jpg{8mm}{RM2-10-b}\,},\;h^0=\begin{pmatrix}-\,\jpg{11mm}{RM2-h2} & 0 \end{pmatrix}^T.
\end{matrix}
\]
In all the cases it is easy to check that the given maps $F,G$ are chain homotopies. Furthermore, $G$ satisfies the conditions of a strong deformation retract, i.e. $G\circ F=\mathrm{id}$, $F\circ G=h^0\circ d^0+d^{-1}\circ h^{-1}$ and $h\circ F=0$. With the help of Proposition~\ref{prop-sdr} we get
\[
\bn{\,\jpg{8mm}{RM2-1-a}\,}\simeq_h\bn{\,\jpg{8mm}{RM2-a}\,}\;\;\;\text{and}\;\;\;\bn{\,\jpg{8mm}{RM2-1-b}\,}\simeq_h\bn{\,\jpg{8mm}{RM2-b}\,}.
\]
Because of this we can follow the proof of Bar-Natan again to show the invariance under the RM3 move. We skip this because this time it is completely analogously to the proof of Bar-Natan (with the maps from above).

The invariance under the virtual Reidemeister moves vRM1, vRM2 and vRM3 follows from Lemma~\ref{lem-virtualisation}. Therefore, the only move left is the mixed Reidemeister move mRM. We have
\[
\bn{\,\jpg{8mm}{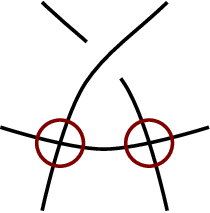}\,}=\Gamma(\bn{\,\jpg{8mm}{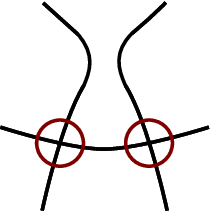}\,}\xrightarrow{\varphi}\bn{\,\jpg{8mm}{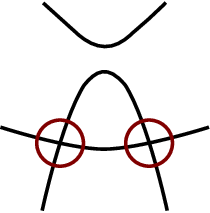}\,})\phantom{.}
\]
and
\[
\bn{\,\jpg{8mm}{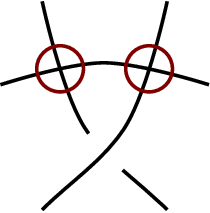}\,}=\Gamma(\bn{\,\jpg{8mm}{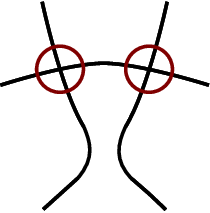}\,}\xrightarrow{\varphi^{\prime}}\bn{\,\jpg{8mm}{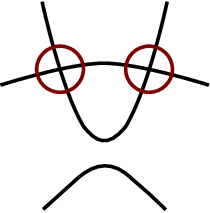}\,}).
\]
There is a vRM2 move in the rightmost parts of both cones. This move can be resolved. Hence, the complex changes only up to an isomorphism (see Lemma~\ref{lem-virtualisation}). Therefore, we have
\[
\bn{\,\jpg{8mm}{mRM}\,}\simeq\Gamma(\bn{\,\jpg{8mm}{mRM-0}\,}\xrightarrow{\varphi}\bn{\,\jpg{8mm}{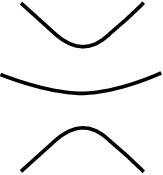}\,})\phantom{.}
\]
and
\[
\bn{\,\jpg{8mm}{mRM2}\,}\simeq\Gamma(\bn{\,\jpg{8mm}{mRM2-0}\,}\xrightarrow{\varphi^{\prime}}\bn{\,\jpg{8mm}{mRM-12}\,}).
\]
Thus, we see that the left and right parts of the cones are equal complexes. Hence, the complexes of two v-links diagrams which differ only through a mRM move are isomorphic. This finishes the proof, because with the obvious chain homotopy $h=0$, isomorphisms induced by the virtual Reidemeister cobordisms and Proposition~\ref{prop-sdr} again gives the desired
\[
\bn{\,\jpg{8mm}{mRM}\,}\simeq_h\bn{\,\jpg{8mm}{mRM2}\,}.
\]
This finishes the proof.
\end{proof}
A question which arises from Theorem~\ref{thm-geoinvarianz} is if the topological complex yields any new information for c-links (compared to the classical Khovanov complex). The following theorem answers this question negatively, i.e. the complex from Definition~\ref{defn-topcomplex} is the classical complex up to chain isomorphisms. It should be noted that Theorem~\ref{thm-geoinvarianz} and Theorem~\ref{thm-classic} imply that our construction can be seen as an extension of Bar-Natan's cobordism based complex to v-links. 

To see this we mention that the cobordisms $m^{++}_+,\Delta^+_{++}$ have the same behaviour as the classical (co)multiplications. Therefore, let $\bn{L_D}_c$ denote the classical Khovanov complex, i.e. every pantsup- or pantsdown-cobordisms are of the form $m^{++}_+,\Delta^+_{++}$ and we add the usual extra signs (e.g. see for example Sec. 3.2 in~\cite{bn1}). Beware that this complex is in general \textit{not} a chain complex for an arbitrary v-link diagram $L_D$. But it is indeed a chain complex for any c-link diagram, i.e. a diagram without v-crossings.
\begin{thm}\label{thm-classic}
Let $L_D$ be a c-link diagram. Then $\bn{L_D}$ and $\bn{L_D}_c$ are chain isomorphic.
\end{thm}
\begin{proof}
Because $L_D$ does not contain any v-crossing, the complex has no $\theta$-saddles. Moreover, every circle is a c-circle. Hence, we can orient them $+$ or $-$, i.e. counterclockwise or clockwise. We choose any numbering for the circles.

Because every circle is oriented clockwise or counterclockwise, every saddle $S$ is of the form $\ud\rightarrow\ril$ or $\du\rightarrow\ler$. Hence, every saddle is of the form $m^{++}_+=m^{--}_-$, $\Delta^+_{++}$ or $\Delta^-_{--}$. Thus, these maps are the classical maps (up to a sign).

We prove the theorem by a spanning tree argument, i.e. choose such a spanning tree. Start at the rightmost leaves and reorient the circles in such a way that the maps which belong to the edges in the tree are the classical maps $m^{++}_+$ or $\Delta^+_{++}$. This is possible because we can use $m^{++}_+=m^{--}_-$ here. We do this until we reach the end.

We repeat the process rearranging the numbering in such a way that the corresponding maps have the same sign as in the classical Khovanov complex. This is possible because every face has an odd number of minus signs (if we count the sign from the relation $\Delta^-_{--}=-\Delta^+_{++}$).

Note that such rearranging does not affect the anti-commutativity because of Lemma~\ref{lem-commutativeindependence}. Hence, after we reach the end every saddle is the classical saddle together with the classical sign. The change of orientations/numberings does not change the complex because of Lemma~\ref{lem-commutativeindependence}. This finishes the proof. 
\end{proof}
\begin{rem}\label{rem-gradings}
We could use the Euler characteristic to introduce the structure of a graded category on $\ucob_R(\emptyset)$ (and hence on $\ukob_R$).

The differentials in the topological complex from Definition~\ref{defn-topcomplex} have all $\deg=0$ (after a degree shift), because their Euler characteristic is -1 (see Lemma~\ref{cor-euler2}). Then it is easy to prove that the topological complex is a v-link invariant under graded homotopy.
\end{rem}
\begin{rem}\label{rem-tutu}
If one does the same construction as above in the category $\ucob_R(\emptyset)^*$, then the whole construction becomes easier in the following sense. First one does not need to work with the saddle signs anymore, i.e. the complex will be a well-defined chain complex if one uses the same signs as in the classical case. Furthermore, most of the constructions and arguments to ensure that everything is a well-defined chain complex are not necessary or trivial, e.g. most parts of the next subsection are ``obviously'' true, and the rest of this subsection can be proven completely analogously. This construction leads us to an equivalent of the construction of Turaev and Turner (see Sec. 3 in~\cite{tutu}). Note that this version does not generalise the classical Khovanov homology. In order to get a bi-graded complex one seems to need a construction related to $\wedge$-products.
\end{rem}
\subsection{The technical points of the construction}
In this subsection we give the arguments for why the topological complex is well-defined and independent of all choices involved.

The following lemma ensures that we can choose the x-marker and the order of the v-circles in the resolutions without changing the total number of signs mod two for all faces.
\begin{lem}\label{lem-xmarker}
Let $F$ be a face of the v-link diagram $L_D$ for a fixed choice of x-markers and numberings for the v-circles of the resolutions of $L_D$. Let $F^{\prime},F^{\prime\prime}$ denote the same face, but $F^{\prime}$ with a different choice of x-markers and $F^{\prime\prime}$ with a different choice for the numberings. Then
\[
\mathrm{sgn}(F)=\mathrm{sgn}(F^{\prime})=\mathrm{sgn}(F^{\prime\prime}).
\]
\end{lem}
\begin{proof}
It is sufficient to show that the statement holds if we only change one x-marker of one crossing $c$ or the numbers of only two v-circles with consecutive numbers in a fixed resolution $\gamma_a$. Moreover, the statement is clear if $c$ or $\gamma_a$ does not affect $F$ at all or one of the saddles is non-orientable.

Note that a change of the x-marker of $c$ affects exactly two saddles $S,S^{\prime}$ of $F$ and for both the number $\mathrm{sgn}(S),\mathrm{sgn}(S^{\prime})$ changes since, by definition, we demand that in the definition of the permutations $\sigma_1,\sigma_2$ from Definition~\ref{defn-sign} the two corresponding v-circles are ordered. Hence, the total change for the face is zero mod two.

If the numbering of the two v-circles changes in $\gamma_{00}$, then the sign of the permutation $\sigma_1$ changes for both saddles $S_{0*},S_{*0}$ but there are no changes for $S_{1*},S_{*1}$. Analogously for the $\gamma_{11}$ case.

In contrast, if the numbering of the two v-circles changes in $\gamma_{01}$, then the sign of the permutations $\sigma_1$ and $\sigma_2$ changes for $S_{*1}$ and $S_{0*}$ respectively, but no changes for $S_{*0},S_{1*}$. Analogously for the $\gamma_{10}$ case. Hence, no change for $F$ mod two.
\end{proof}
\begin{lem}\label{lem-commutativeindependence}
Let $L_D$ be a v-link diagram and let $\bn{L_D}_1$ be its topological complex from Definition~\ref{defn-topcomplex} with arbitrary orientations for the resolutions. Let $\bn{L_D}_2$ be the complex with the same orientations for the resolutions except for one circle $c$ in one resolution $\gamma_a$. If a face $F_1$ from $\bn{L_D}_1$ is anti-commutative, then the corresponding face $F_2$ from $\bn{L_D}_2$ is also anti-commutative.

Moreover, if $\bn{L_D}_1$ is a well-defined chain complex, then it is isomorphic to $\bn{L_D}_2$, which is also a well-defined chain complex.

The same statement is true if the difference between the two complexes is the numbering of the crossings, the choice of the x-marker, rotations/isotopies of the v-link diagram or the fixed numbering of the v-circles in the resolutions.
\end{lem}
\begin{proof}
Assume that the face $F_1$ is anti-commutative. Then the different orientations of the circle $c$ correspond to a composition of all morphisms of the face $F_2$ with this circle as a boundary component with $\Phi^-_+$.

Hence, the face $F_2$ is also anti-commutative, because both outgoing (or incoming) morphisms of $F_2$ are composed with an extra $\Phi^-_+$ if the circle is in the first (or last) resolution of the faces. If it is in one of the middle resolutions, then we have to use the relation $\Phi^-_+\circ\Phi^-_+=\mathrm{id}^+_+$ from Lemma~\ref{lem-basiscalculations}. Note that it is important for this argument to work that cylinders between differently oriented v-circles are $\Phi^-_+$, as in Definition~\ref{defn-topcomplex}.

Thus, if the first complex is a well-defined chain complex, then the same is true for the second. The isomorphism is induced (using a spanning tree construction) by the isomorphism $\Phi^-_+$.

The second statement is true because the numbering of the crossings does not affect the cobordisms at all. Hence, the argument can be shown analogously to the classical case (see for example Theorem 1 in~\cite{bn2}), but it should be noted that our way of spreading signs does not depend on the numbering of the crossings (in contrast to the classical case).

On the third point: That anti-commutative faces stay anti-commutative, if one changes between the two possible choices in Definition~\ref{defn-sign}, is part of Lemma~\ref{lem-xmarker}. The chain isomorphism is induced (using a spanning tree construction) by a sign permutation.

The penultimate statement follows directly from the definition of the saddle sign and decorations, while the last statement is also part of Lemma~\ref{lem-xmarker} with the isomorphisms again induced (using a spanning tree construction) by a sign permutation.
\end{proof}
\begin{lem}\label{lem-virtualisation}
Let $L_D,L^{\prime}_D$ be v-link diagrams which differs only by a virtualisation of one crossing $c$. If a face $F$ is anti-commutative in $\bn{L_D}$, then the corresponding face $F^{\prime}$ is anti-commutative in $\bn{L^{\prime}_D}$. Moreover, if $\bn{L_D}$ is a well-defined chain complex, then it is isomorphic to $\bn{L^{\prime}_D}$, which is also a well-defined chain complex.

The same statement is true if $L_D$ and $L^{\prime}_D$ differ only by a vRM1, vRM2, vRM3 or mRM move.
\end{lem}
\begin{proof}
The statement about anti-commutativity is clear, if one of the saddles which belongs to the crossing $c$ is non-orientable. This is true because of the relations from Equation~\ref{eq-combrel3} and Proposition~\ref{prop-nonorientablefaces}. Thus, we can assume that both saddles are orientable. Furthermore, it is clear that the two compositions of the saddles are boundary-preservingly homeomorphic after the virtualisation. Hence, the only thing we have to ensure is that the decorations and signs work out correctly.

We use the Lemma~\ref{lem-commutativeindependence} here, i.e. we can choose the orientations and the numberings in such a way that the saddles which do not belong to the crossing $c$ have the same local orientations and numberings. We observe the following. The sign and the local orientations of a saddle can only change if the saddle belongs to the crossing $c$, i.e. the local orientations always change (see Fig.~\ref{figure-virt2}) and the sign changes precisely if the two strings in the bottom picture of Fig.~\ref{figure-virt2} are part of two different v-circles.
\begin{figure}[ht]
     \centerline{\includegraphics[scale=0.75]{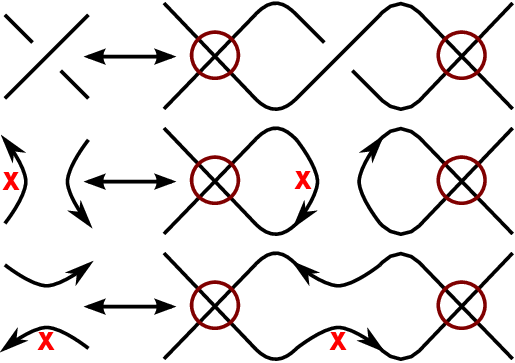}}
  \caption{The behaviour of the x-marker and orientations under virtualisation.}
  \label{figure-virt2}
\end{figure}

A change of the local orientations multiplies an extra sign for comultiplication, but no such extra sign for multiplication. This follows from Table 1 and the relations from Equation~\ref{eq-combrel12}. Hence, the anti-commutativity still holds if the two saddles which belong to the crossing $c$ are both multiplications or comultiplications, because their decorations and signs change in the same way.

If one is a multiplication and one is a comultiplication, then we have two cases: Either the multiplication gets an extra sign or it does not. The comultiplication always gets an extra sign because the local orientations change. But the multiplication will change its saddle sign iff the comultiplication does not change its saddle sign. Hence, the number of extra signs does not change mod two. This ensures that the faces stay anti-commutative.

That the face $F^{\prime}$ stays anti-commutative after a vRM1, vRM2, vRM3 or mRM move follows because neither the local orientations nor the signs of any cobordism change. Thus, all decorations and signs are the same.

The chain isomorphisms are induced by the vRM-cobordisms shown in Fig.~\ref{figure0-reide}, morphisms of type $\Phi^-_+$ and identities. Recall that all these cobordisms are isomorphisms in our category.
\end{proof}
For the proof of the next lemma we refer the reader to the paper~\cite{ma2}. We call faces of the following type the \textit{basic (non-)orientable faces}.
\begin{figure}[ht]
  \centerline{\begin{minipage}[c]{6,9cm}
	\includegraphics[scale=0.35]{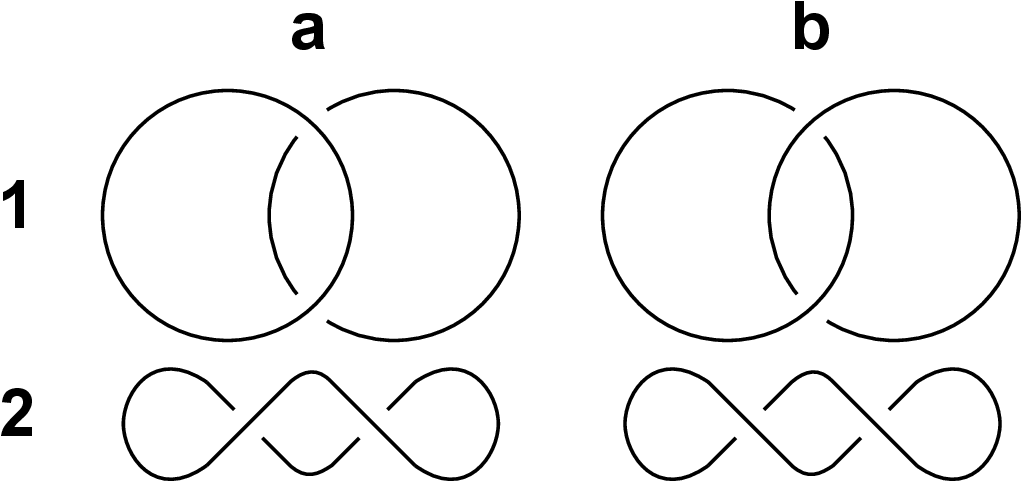}
\end{minipage}
\begin{minipage}[c]{6,9cm}
	\includegraphics[scale=0.30]{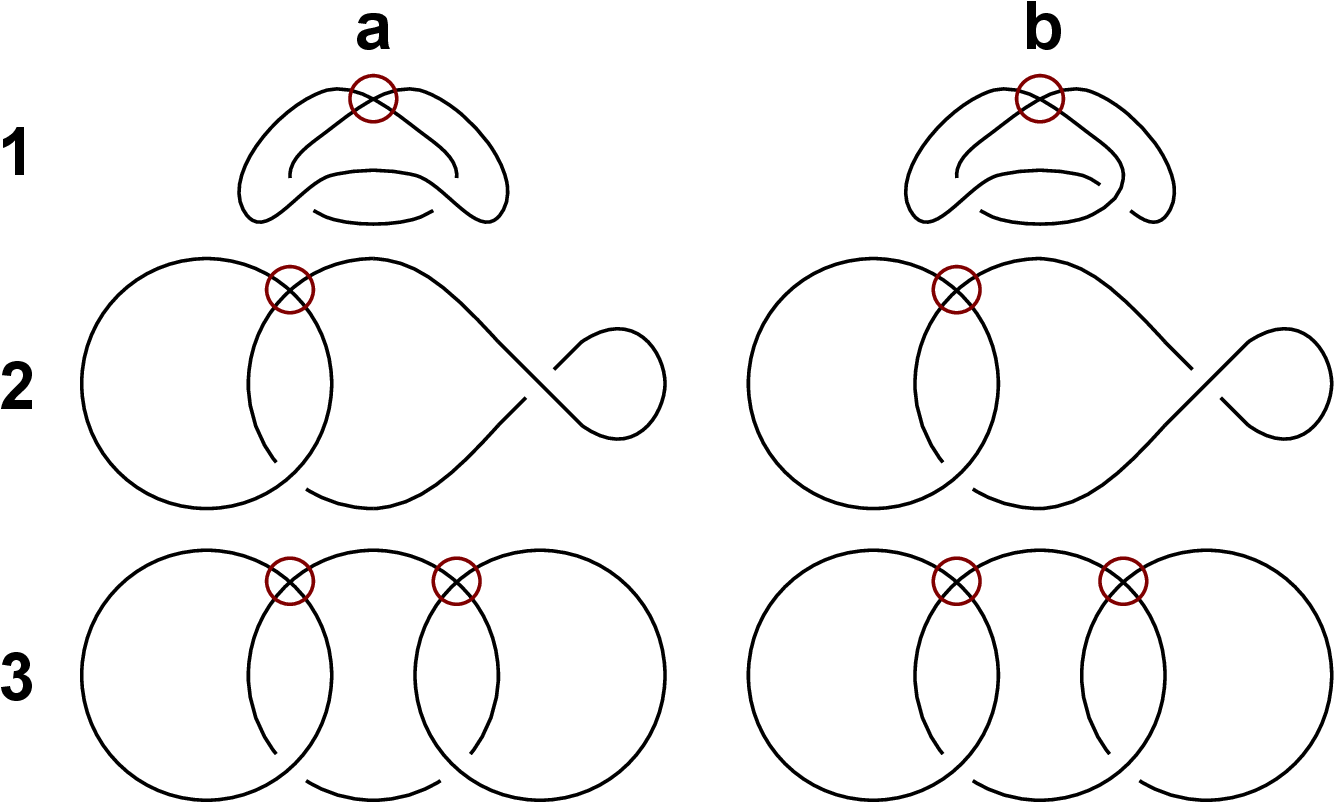}
\end{minipage}}
\caption{Left: The basic orientable faces. Right: The basic non-orientable faces.}\label{figure-basic}
\end{figure} 
\begin{lem}\label{lem-basicfaces}
Let $L_D$ be a v-link diagram. Then $L_D$ can be reduced by a finite sequence of isotopies, vRM1, vRM2, vRM3, mRM moves and virtualisations to a v-link diagram $L^{\prime}_D$ in such a way that a fixed connected face of $L^{\prime}_D$ is isotopic to one of the basic faces from Fig.~\ref{figure-basic} (or to one of their mirror images) up to vRM1, vRM2, vRM3 moves on the face itself.
\end{lem}
\begin{proof}
See Manturov's paper, i.e. the discussion after Lemma 2 in~\cite{ma2}.
\end{proof}
Note that these lemmata allow us to check arbitrary orientations on the basic faces with arbitrary numbering of crossings and components.
\begin{prop}\label{prop-standardfaces}
Let $L_D$ be a v-link diagram with a diagram which is isotopic to one of the projections from Fig.~\ref{figure-basic}. Then $\bn{L_D}$ is a chain complex, i.e. the basic faces are anti-commutative. Moreover, disjoint faces, i.e. faces such that the corresponding four-valent graph is unconnected, are always anti-commutative.
\end{prop}
\begin{proof}
Because of Lemma~\ref{lem-commutativeindependence}, we only need to check that the faces are anti-commutative for orientations of the resolutions of our choice with an arbitrary numbering. Then we are left with three different cases: Either the v-link diagram of $L_D$ is orientable, i.e. all saddles are orientable, or the face is non-orientable, i.e. two or four of the saddles are non-orientable, or the face is disjoint.

For the first case we see that every resolution contains only c-circles. We prove the anti-commutativity of the corresponding face for the following orientations of the resolutions. All appearing circles are numbered in ascending order from left to right or outside to inside. Moreover, the position of the x-marker does not affect our argument and we suppress it in this proof.

Because every resolution contains only c-circles, we choose a positive orientation for the circles except for the two nested circles that appear in two resolution of a face of type 1a or 1b. This is a clockwise orientation for all the non-nested circles and a counterclockwise orientation for the two nested circles. Hence, all appearing cylinders are identities.

It follows from this convention that every $0$-resolution (or $1$-resolution) $\smoothing$ of a crossing $\slashoverback$ (or a crossing $\backoverslash$) is of the form $\du$ and every $1$-resolution (or $0$-resolution) $\smoothing$ of a crossing $\slashoverback$ (or a crossing $\backoverslash$) is of the form $\ler$. Moreover, the only face with an even number of saddle signs is of type 1a.

All we need to do is compare these local orientations with the ones from Table 1. We see that we have to check (indicated by the $\overset{!}{=}$) the following equations.
\begin{itemize}
\item $\Delta^{++}_+\circ m^-_{--}\overset{!}{=}-\Delta^{--}_-\circ m^+_{++}$ (face of type 1a).
\item $m^{++}_+\circ\Delta^+_{++}\overset{!}{=}m^{++}_+\circ\Delta^+_{++}$ (face of type 1b).
\item $(\Delta^+_{++}\amalg\,\mathrm{id}^+_+)\circ(\mathrm{id}^+_+\amalg\,m^{++}_+)\overset{!}{=}m^{++}_+\circ\Delta^{++}_+$ (face of type 2a).
\item $m^{++}_{+}\circ(m^{++}_-\amalg\,\mathrm{id}^+_+) =m^{++}_{+}\circ(\mathrm{id}^+_+\amalg\,m^{++}_+)$ (face of type 2b).
\end{itemize}
Most of these equations are easy to calculate. The reader should check that the cobordisms on the left and the right side of every equation are homeomorphic (using Proposition~\ref{prop-nonorientablefaces} and Lemma~\ref{lem-basiscalculations}).

Furthermore, the second equation is clear and the other three follows easily using the result of Lemma~\ref{lem-basiscalculations}. Hence, they are all anti-commutative because only the first face has an even number of saddle signs.

The non-orientable faces of type 1b, 2a, 2b, 3a and 3b are easy to check. One can use the Euler characteristic here and the relations in Eq.~\ref{eq-combrel3}.

The non-orientable face of type 1a is the face from~\ref{probcube}. Here we have to use Proposition~\ref{prop-nonorientablefaces}. We get two $\theta$-cobordisms and a $\Delta$- and a $m$-cobordism. Because of the relations in Eq.~\ref{eq-combrel3} we can ignore the saddle signs.

Again we can choose an orientation for the resolutions. We can do this for example in the following way (compare to Fig.~\ref{figure0-big}).
\begin{itemize}
 \item The first M\"obius strips are $\theta\colon\dd\rightarrow\lel$ and $\theta\colon\dd\rightarrow\rir$.
 \item The pantsdown is $\Delta^+_{-+}\colon\du\rightarrow\lel$ and the pantsup is $m^{--}_+\colon\ud\rightarrow\ler$.
\end{itemize}
We use Proposition~\ref{prop-nonorientablefaces} to see that this face is anti-commutative.

The reader should check that all disjoint faces with only orientable saddles have an odd number of saddle signs. The disjoint faces with two or four non-orientable saddles anti-commute because of the relations in Eq.~\ref{eq-combrel3} and (k) of Lemma~\ref{lem-basiscalculations}.
\end{proof}
This proposition leads us to an important theorem and an easy corollary.
\begin{thm}\label{theo-facescommute}(\textbf{Faces commute}) Let $L_D$ be a v-link diagram. Let $\bn{L_D}$ be the complex from Definition~\ref{defn-topcomplex} with arbitrary possible choices. Then every face of the complex $\bn{L_D}$ is anti-commutative.
\end{thm}
\begin{proof}
This is a direct consequence of the Proposition~\ref{prop-standardfaces} and the three Lemmata~\ref{lem-commutativeindependence},~\ref{lem-virtualisation} and~\ref{lem-basicfaces}.
\end{proof}
\begin{cor}\label{cor-chaincomplex}
The complex $\bn{L_D}$ is a chain complex. Thus, it is an object in the category $\ukob_R$. 
\end{cor}
\begin{proof}
We can directly use Theorem~\ref{theo-facescommute}.
\end{proof}
\section{Skew-extended Frobenius algebras}\label{sec-vkhfa}
We note that this section has three subsections. We construct the ``algebraic'' complex of a v-link diagram $L_D$ in the first subsection. Its homotopy type is an invariant of virtual links $L$, i.e. invariant modulo the generalised Reidemeister moves from Fig.~\ref{figure0-reide}. It should be noted that, in contrast to the topological complex, the notion of ``homology'' makes sense for the algebraic complex.
\vspace*{0.5cm}

We describe the relation between uTQFTs and skew-extended Frobenius algebras in the second subsection. A relation of this kind was discovered by Turaev and Turner for extended Frobenius algebras and the functors they use (see Proposition 2.9 in~\cite{tutu}). Even though our construction is different, their ideas can be used in our context too. This is the main part of Theorem~\ref{thm-tututheo}. But our uTQFT correspond to \textit{skew-extended} Frobenius algebras, i.e. the map $\Phi$ is a skew-involution rather than an involution.
\vspace*{0.5cm}

In the last subsection we are able to classify all aspherical uTQFTs which can be used to define v-link invariants, see Theorem~\ref{thm-class}. It is worth noting that we get an invariant for v-links which is an extension of the Khovanov complex for $R=\bZ$ or $R=\bQ$, see Corollary~\ref{cor-khovhomo}. Note that this includes that our construction can be seen as a \textit{categorification} of the virtual Jones polynomial, see Corollary~\ref{cor-khovhomocat}. Moreover, we also get extension for other classical link homologies, see e.g. the Corollaries~\ref{cor-leekhovhomo} and~\ref{cor-bnhomo}.
\subsection{The algebraic complex}
We denote any v-link diagram of the unknot by the symbol $\bigcirc$. Furthermore, we view v-circles, i.e. v-links without classical crossings, as disjoint circles immersed into $\bR^2$. Recall that $R$ is always an unital, commutative ring of arbitrary characteristic.
\begin{defn}\label{defn-utqft}(\textbf{uTQFT}) A \textit{(1+1)-dimensional unoriented TQFT} $\mathcal F$ (we call this a \textit{uTQFT}) is a strict, symmetric, covariant, $R$-pre-additive functor
\[
\mathcal F\colon\ucob_R(\emptyset)\rightarrow\RMOD.
\]
Here $\mathcal F(\bigcirc)$ is a finitely generated, free $R$-module. Let $\mathcal O,\mathcal O^{\prime}$ be two homeomorphic objects from $\ucob_R(\emptyset)$. Then $\mathcal F(\mathcal O)=\mathcal F(\mathcal O^{\prime})$ should hold. The functor $\mathcal F$ should also satisfy the following axioms.
\begin{itemize}
\item[(1)] Let $\mathcal O,\mathcal O^{\prime}$ be two disjoint objects in $\Ob(\ucob_R(\emptyset))$. Then there exists a natural (with respect to homeomorphisms) isomorphism between $\mathcal F(\mathcal O\amalg \mathcal O^{\prime})$ and $\mathcal F(\mathcal O)\otimes\mathcal F(\mathcal O^{\prime})$.
\item[(2)] The functor satisfies $\mathcal F(\emptyset)=R$.
\item[(3)] For a cobordism $C \colon \mathcal O \to \mathcal O^{\prime}\in\Mor(\ucob_R(\emptyset))$ the homomorphism $\mathcal F(C)$ is natural with respect to homeomorphisms of cobordisms.
\item[(4)] Let the cobordism $C\colon \mathcal O \to \mathcal O^{\prime}\in\Mor(\ucob_R(\emptyset))$ be a disjoint union of the two cobordisms $C_{1,2}$. Then $\mathcal F(C)=\mathcal F(C_1)\otimes\mathcal F(C_2)$ under the identification from axiom (1).
\end{itemize}
Two uTQFTs $\mathcal{F,F^{\prime}}$ are called \textit{isomorphic} if for each object of $\mathcal O\in\Ob(\ucob_R(\emptyset))$ there is an isomorphism $\mathcal F(\mathcal O)\rightarrow\mathcal F^{\prime}(\mathcal O)$, natural with respect to homeomorphisms of the objects and homeomorphisms of cobordisms, multiplicative with respect to disjoint union and the isomorphism assigned to $\emptyset$ is the identity morphism. 
\end{defn}
\begin{rem}
There are several things to note about the definition.
\begin{itemize}
\item Recall that our category is $R$-pre-additive. An uTQFT is a $R$-pre-additive functor. So we can extend this to a functor
\[
\mathcal F\colon\ukob_R\rightarrow\Kom_b(\mat(\RMOD)),
\]
i.e. for every formal chain complex $(C_*,d_*)$ of objects of $\ucob_R(\emptyset)$, i.e. v-circles, the object $\mathcal F((C_*,d_*))$ is a chain complex of $R$-modules and for every formal chain map $f\colon (C_*,d_*)\to (C^{\prime}_*,d^{\prime}_*)$ of possibly non-orientable, decorated cobordisms the morphism $\mathcal F(f)$ is a chain map of $R$-module homomorphisms.
\item An uTQFT $\mathcal F$ is a covariant functor. Hence, we see that $\mathcal F(\mathrm{id}^+_+)=\mathrm{id}$. Furthermore it is symmetric and hence $\mathcal F(\tau^{++}_{++})=\tau$. Here $\tau$ denotes the canonical permutation.
\item The permutation $\tau^{++}_{++}$ is natural. So we can assume that $A\otimes B$ and $B\otimes A$ are equal and not merely isomorphic.
\item For the definition of natural (converted to our setting) we refer the reader to Chapter III in~\cite{tur2}.
\end{itemize}
\end{rem}
\begin{defn}\label{defn-alcomplex}(\textbf{Algebraic complex}) Let $L_D$ be a v-link diagram. Then the \textit{algebraic complex of $L_D$} induced by the uTQFT $\mathcal F$ is the complex $\mathcal F(\llbracket L_D\rrbracket)$. 
\end{defn}
We have the following important result. Here $L_D,L^{\prime}_D$ are a v-link diagrams. The proof is a direct consequence of Theorem~\ref{thm-geoinvarianz}.
\begin{thm}\label{thm-algcomplex}(\textbf{The algebraic complex is an invariant}) Let $\mathcal F$ an uTQFT which satisfies the Bar-Natan-relations of Fig.~\ref{figureintroa-4}. Then the algebraic complex $\mathcal F(\llbracket L_D\rrbracket)$ is a v-link invariant in the following sense.

For two equivalent (up to the generalised Reidemeister moves) v-link diagrams $L_D,L^{\prime}_D$ the two chain complexes $\mathcal F(\llbracket L_D\rrbracket)$ and $\mathcal F(\llbracket L^{\prime}_D\rrbracket)$ are equal up to chain homotopy.
\end{thm}
\begin{proof}
Directly from Theorem~\ref{thm-geoinvarianz}.
\end{proof}
This theorem allows us to speak of \textit{the algebraic complex} $\mathcal F(\llbracket L\rrbracket)$ of any oriented v-link $L$. Furthermore, the category $\RMOD$ is abelian. Hence, the category $\Kom_b(\mat(\RMOD))$ is also an abelian category. So unlike in the category $\ukob_R$, we have the notion of homology. We denote the homology of the algebraic chain complex by $H(\mathcal F(\bn L))$.
\subsection{Skew-extended Frobenius algebras and uTQFTs}
We continue with the definition of an algebra that we call a \textit{skew-extended Frobenius algebra}. For a $R$-algebra $A$ with comultiplication $\Delta$ and counit $\varepsilon$ we call an $R$-algebra homomorphism $\Phi\colon A\to A$ a \textit{skew-involution} if it satisfies the following.
\begin{itemize}
\item[(a)] $\Phi^2=\mathrm{id}$ (involution).
\item[(b)] $(\Phi\otimes\Phi)\circ\Delta\circ\Phi=-\Delta$ and $\varepsilon\circ\Phi=-\varepsilon$ (skew-property).
\end{itemize}
\begin{defn}\label{defn-exfrob}(\textbf{Skew-extended Frobenius algebras}) A \textit{Frobenius algebra $A$ over $R$} is a unital, commutative\footnote{For us a Frobenius algebra is necessary commutative. Beware that some authors distinguish between commutative and non-commutative Frobenius algebras.} algebra over $R$ which is projective and of finite type (as an $R$-module), together with a module homomorphism $\varepsilon\colon A\to R$, such that the bilinear form $\left\langle \cdot,\cdot\right\rangle$ defined by $\left\langle a,b\right\rangle=\varepsilon(m(a,b))$ for all $a,b\in A$ is non-degenerate. Here $m=m(\cdot,\cdot)\colon A\otimes A\to A$ denotes the multiplication in $A$.

An \textit{skew-extended Frobenius algebra $A$ over $R$} is a Frobenius algebra together with a skew-involution of Frobenius algebras $\Phi\colon A\to A$ and an element $\theta\in A$ which satisfies the two equations below. Note that one can use $\varepsilon$ together with the multiplication $m$ to define a comultiplication $\Delta\colon A\to A\otimes A$. In fact, there are around a dozen equivalent definitions of a Frobenius algebra and some include the comultiplication in the definition. Street list a few of them in~\cite{street}.
\vspace*{0.5cm}

The two defining equations for a skew-extended Frobenius algebra $A$ over $R$ are the following.
\begin{itemize}
\item[(1)] $\Phi(\theta a)=\theta a=\theta\Phi(a)$ for all $a\in A$.
\item[(2)] $(m\circ(\Phi\otimes\mathrm{id})\circ\Delta)(1)=\theta^2$. 
\end{itemize}
\end{defn}
\begin{nota}
Because of $1\in A$, we can define $\iota\colon R\to A$ by $1\mapsto 1$. We can write a Frobenius algebra uniquely as $\mathcal F=(R,A,\varepsilon,\Delta)$. Moreover, we can write a skew-extended Frobenius algebra $\mathcal F$ uniquely as $\mathcal F=(R,A,\varepsilon,\Delta,\Phi,\theta)$.
\end{nota}
\begin{defn}\label{defn-isofrobenius}
Two skew-extended Frobenius algebras, denoted by $\mathcal F_1=(R,A,\varepsilon,\Delta,\Phi,\theta)$ and $\mathcal F_2=(R,A^{\prime},\varepsilon^{\prime},\Delta^{\prime},\Phi^{\prime},\theta^{\prime})$, are called \textit{isomorphic} if there exists an isomorphism of Frobenius algebras $f\colon A\rightarrow A^{\prime}$, which satisfies $f(\theta)=\theta^{\prime}$ and $f\circ\Phi=\Phi^{\prime}\circ f$.

We call a Frobenius algebra \textit{aspherical} if $\varepsilon(\iota(1))=0$. Furthermore, we say it is a \textit{rank2-Frobenius algebra} if $A\cong R\cdot 1\oplus R\cdot X$ for some elements $1,X\in A$ as $R$-modules.
\end{defn}
\begin{rem}\label{rem-categoryfrob}
As above, we will call $m\colon A\otimes A\rightarrow A$ the \textit{multiplication} of $A$. The map $\varepsilon$ is called the \textit{counit} of $A$. It, together with the multiplication $m$, can be used to define a \textit{comultiplication} $\Delta\colon A\to A\otimes A$. The coproduct and the product make the two diagrams
\[
\begin{xy}
  \xymatrix{
 A\otimes A\ar[r]^{\mathrm{id}\otimes\Delta}\ar[d]_{\Delta\otimes\mathrm{id}}\ar[dr]|/.3em/{\Delta\circ m} & A\otimes A\otimes A\ar[d]^{m\otimes\mathrm{id}} & & A\ar[dr]|{\mathrm{id}}\ar[r]^{\Delta}\ar[d]_{\Delta} & A\otimes A\ar[d]^{\mathrm{id}\otimes\varepsilon}\\
 A\otimes A\otimes A\ar[r]_{\mathrm{id}\otimes m} & A\otimes A  & & A\otimes A\ar[r]_{\varepsilon\otimes\mathrm{id}} & A}
\end{xy}
\]
commutative. In a skew-extended Frobenius algebra the skew-involution $\Phi$ and the element $\theta$ make the two diagrams
\[
\begin{xy}
  \xymatrix{
& A\ar[rd]^{\Phi} & & & & A\otimes A\ar[rd]^{m^{\prime}} & \\
A\ar[ru]^{\cdot\theta}\ar[rr]_{\cdot\theta} & & A & & A\ar[ru]^{\Delta}\ar[rd]_{\cdot\theta} &  & A\\
 & & & & & A\ar[ru]_{\cdot\theta} & }
\end{xy}
\]
commutative (it is easy to check that the two equations from Definition~\ref{defn-exfrob} already imply the equation $(m\circ(\phi\otimes\mathrm{id})\circ\Delta)(a)=\theta^2 a$ for all $a\in A$). Here the map $\cdot\theta\colon A\to A$ is the multiplication with $\theta$ and the map $m^{\prime}\colon A\otimes A\to A$ is the map $m\circ(\Phi\otimes\mathrm{id})$.

We recognise that the lower right diagram is the problematic face from~\ref{probcube}. So the second equation from Definition~\ref{defn-exfrob} is a key point in the definition. 
\end{rem}
The following theorem is inspired by Proposition 2.9 in~\cite{tutu} and the ``classical'' correspondence between commutative Frobenius algebras and TQFTs, see for example Theorem 3.3.2 in~\cite{ko}.
\begin{thm}\label{thm-tututheo}
The isomorphism classes of (1+1)-dimensional uTQFTs over $R$ are in bijective correspondence with the isomorphism classes of skew-extended Frobenius algebras over $R$.
\end{thm}
\begin{proof}
First let us consider an uTQFT $\mathcal F$ over $R$. We describe a way to get a skew-extended Frobenius algebra from it. Let us denote this algebra by $(R,A,\varepsilon,\Delta,\Phi,\theta)$.

We take $A=\mathcal F(\bigcirc)$ as our underlying $R$-module. Next we need a skew-involution $\Phi\colon A\to A$. We take the second cylinder from Fig.~\ref{figure1-1}: Set $\Phi=\mathcal F(\Phi^-_+)$.

The unit $\iota$ should be $\mathcal F(\iota_+)$. There is no further choice because $\iota_+=\iota_-$. The counit should be $\mathcal F(\varepsilon^+)$. Here we have a choice because $\varepsilon^+\neq\varepsilon^-$. But because of $\varepsilon^+=-\varepsilon^-$, both choices lead to isomorphic algebras.

Now we need a multiplication $m$ and a comultiplication $\Delta$. One may suspect, that we have different choices for either of them, namely the eight $m^{\pm\pm}_{\pm},\Delta^{\pm}_{\pm\pm}$. But the relations of a Frobenius algebra only allow one option. We discuss this now. It should be noted that the computations below can be done using Lemma~\ref{lem-basiscalculations}.
\begin{itemize}
\item The lower boundary components of $\Delta^u_{l_1l_2}$ must have the same glueing numbers as the boundary component of $\varepsilon^+$ because $\mathcal F(\varepsilon^+)$ should be the counit.
\item Because of the relation $\varepsilon\circ m\circ(\mathrm{id}\otimes\iota)=\varepsilon=\varepsilon\circ m\circ(\iota\otimes\mathrm{id})$, the lower boundary of $m^{u_1u_2}_l$ must have the same glueing number as the boundary component of $\varepsilon^+$. The same is true for the upper boundary (this means we need $m^{++}_+=m^{--}_-$).
\item Because of the relation $(\mathrm{id}\otimes\,m)\circ(\Delta\otimes\mathrm{id})=\Delta\circ m=(m\otimes\mathrm{id})\circ(\mathrm{id}\otimes\Delta)$, the $m^{u_1u_2}_{l}$ must have the same glueing number on the lower boundary as the upper boundary of $\Delta^u_{l_1l_2}$ (the reader should check that this is the only possible choice for the glueing numbers for $m^{u_1u_2}_{l}$ and $\Delta^u_{l_1l_2}$).
\end{itemize}
Therefore, we have $\mathcal F(\iota_+)=\iota$, $\mathcal F(\varepsilon^+)=\varepsilon$, $\mathcal F(m^{++}_+)=m$ and $\mathcal F(\Delta^+_{++})=\Delta$.

The last piece missing is the element $\theta\in A$. Consider a two times punctured projective plane $\mathbb{RP}^2_2$ (a punctured M\"obius strip). This is $\theta$ in our notation.

Then $\theta\circ\iota_+\colon\emptyset\to\bigcirc$ is a punctured projective plane (hence a M\"obius strip). Set $\theta=\mathcal F(\theta\circ\iota_+)(1)$. Because of the definition, this is an element of $\mathcal F(\bigcirc)=A$.

We have to prove the equations needed for a skew-extended Frobenius algebra, i.e. that $\iota$ is a unit, $\varepsilon$ is a counit, $\Phi$ is a skew-involution, $m$ ($\Delta$) is a (co)multiplication and the commutativity of the faces from Remark~\ref{rem-categoryfrob}.

This is a straightforward verification based on the relations from Lemma~\ref{lem-basiscalculations} (we omit it here). This shows that every uTQFT has an underlying skew-extended Frobenius algebra.
\vspace*{0.5cm}

For the other direction, i.e. if we assume that we have a skew-extended Frobenius algebra, we note that this algebra has an underlying ``classical'' Frobenius algebra. Therefore we get a TQFT $\mathcal F^{\prime}$ from this underlying Frobenius algebra (by for example Theorem 3.3.2 in~\cite{ko}). We want to use this TQFT to define an uTQFT $\mathcal F$. The TQFT $\mathcal F^{\prime}$ is a covariant functor
\[
\mathcal F^{\prime}\colon\cob_R(\emptyset)\rightarrow\RMOD.
\]
Let $\mathcal O$ be an object in $\ucob_R(\emptyset)$. This object gives us (modulo homeomorphisms) a corresponding object $\mathcal O^{\prime}$ in $\cob_R(\emptyset)$. We set $\mathcal F(\mathcal O)=\mathcal F^{\prime}(\mathcal O^{\prime})$. This assignment clearly satisfies that $\mathcal F(\bigcirc)$ is a finitely generated, free $R$-module and $\mathcal F(\mathcal O_1)=\mathcal F(\mathcal O_2)$ for two homeomorphic objects $\mathcal O_1,\mathcal O_2$.

Moreover, because $\mathcal F^{\prime}$ is a TQFT, this satisfies the first two axioms from our Definition~\ref{defn-utqft}. Now we need to define $\mathcal F(\mathcal C)$ for morphisms from $\ucob_R(\emptyset)$.

First we assume that $\mathcal C\colon\mathcal O_1\to\mathcal O_2$ is orientable and connected. Then we have a corresponding morphism in $\cob_R(\emptyset)$, i.e. the same without the boundary decorations, which we will denote by $\mathcal C^{\prime}\colon\mathcal O^{\prime}_1\to\mathcal O^{\prime}_2$.

We denote the cap-, cup-, pantsup- and pantsdown-cobordisms in the category $\cob_R(\emptyset)$ by $\iota,\varepsilon,m$ and $\Delta$ respectively. Let us define
\[
\mathcal F(\iota_+)=\mathcal F^{\prime}(\iota),\hspace*{0.25cm}\mathcal F(\varepsilon^+)=\mathcal F^{\prime}(\varepsilon),\hspace*{0.25cm}\mathcal F(\Phi^-_+)=\Phi
\]
(the leftmost assignment is the important new piece) and
\[
\mathcal F(m^{++}_+)=\mathcal F^{\prime}(m),\hspace*{0.25cm}\mathcal F(\Delta^+_{++})=\mathcal F^{\prime}(\Delta).
\]
The map $\Phi$ is the skew-involution in the skew-extended Frobenius algebra. Thus, we can define $\mathcal F(\mathcal C)$ in the following way. We decompose $\mathcal C^{\prime}$ into the basic pieces $\iota,\varepsilon,m,\Delta$. Then $\mathcal F^{\prime}(\mathcal C^{\prime})$ is independent of this decomposition because $\mathcal F^{\prime}$ is a TQFT, see Theorem 3.3.2 in~\cite{ko}. If we use the same decomposition for $\mathcal C$ (under the identification from above), we get a cobordism $\tilde{\mathcal C}$. For this cobordism we can define $\mathcal F(\tilde{\mathcal C})$. We see that we only have to change some of the boundary decorations of $\tilde{\mathcal C}$ to obtain $\mathcal C$. Hence, we have
\[
\mathcal C=\mathcal C_1\circ\tilde{\mathcal C}\circ\mathcal C_2,
\]
where $\mathcal C_1,\mathcal C_2$ are cylinders of the type $\mathrm{id}^+_+$ or $\Phi^-_+$. Hence, we can define
\[
\mathcal F(\mathcal C)=\mathcal F(\mathcal C_1)\circ\mathcal F(\tilde{\mathcal C})\circ\mathcal F(\mathcal C_2).
\]
That this is also independent of the decomposition follows from the fact that $\mathrm{id}^+_+,\Phi^-_+$ and the corresponding maps in the skew-extended Frobenius algebra are (skew-) involutions and a ``level-by-level''\footnote{One can for example verify the statement by induction on the number of generators in the decomposition.} change of decorations using the relations in Lemma~\ref{lem-basiscalculations}. Moreover, for a non-connected, orientable cobordism $\mathcal C$ we extend the definition from above multiplicatively.
\vspace*{0.5cm}

For a non-orientable, connected cobordism $\mathcal C$ we have to define $\mathcal F(\theta)=\cdot\theta$ first. Here the map $\cdot\theta\colon A\to A$ is the multiplication with the element $\theta$ in our skew-extended Frobenius algebra. Hence, if we decompose $\mathcal C=\mathcal C_{or}\#n\mathbb{RP}^2$ into a (non-decorated) orientable part $\mathcal C_{or}$ and $n$-times a projective plane we define
\[
\mathcal F(\mathcal C)=\theta^n\mathcal F^{\prime}(\mathcal C_{or}).
\]
This is again independent of the decomposition of $\mathcal C_{or}$, because of the first relation in a skew-extended Frobenius algebra, namely $\Phi(\theta a)=\theta a=\theta\Phi(a)$ for all $a\in A$. Furthermore, it is independent from the decomposition $\mathcal C=\mathcal C_{or}\#n\mathbb{RP}^2$, because if we replace a $2-\mathbb{RP}^2$ with a torus $\mathcal T$, we see that $\mathcal F(\mathcal C_{or})$ is multiplied by a factor $(m\circ(\Phi\otimes\mathrm{id})\circ\Delta)(1)\theta^{n-2}$. Hence, using the second relation of the skew-extended Frobenius algebra, we get
\begin{align*}
\mathcal F(\mathcal C_{or}\#n\mathbb{RP}^2)&=(\theta^n)\mathcal F^{\prime}(\mathcal C_{or})=\theta^{n-2}(m\circ(\Phi\otimes\mathrm{id})\circ\Delta)(1)\mathcal F^{\prime}(\mathcal C_{or})\\ &=\mathcal F(\mathcal C_{or}\#\mathcal T\#(n-2)\mathbb{RP}^2).
\end{align*}
For a non-connected, non-orientable cobordism $\mathcal C$ we extend the definition from above multiplicatively. Hence, we only have to show the remaining axioms from the Definition~\ref{defn-utqft}. The reader should check these axioms (one could follow the end of the proof of Lemma 2.11 in~\cite{tutu}).
\end{proof}
\subsection{Classification of v-link homologies}
From now on we use the notions uTQFT and skew-extended Frobenius algebra interchangeably.
\begin{prop}\label{cor-univalg}(\textbf{The universal skew-extended Frobenius algebra}) Every rank2-uTQFT that is aspherical comes from the rank2-uTQFT $\mathcal F_U=(R_U,A_U,\varepsilon_U,\Delta_U,\Phi_U,\theta_U)$ through base change. Here the ring $R_U$ is $R_U=\bZ[a,a^{-1},\alpha,\beta,\gamma,t]/\mathcal I$ with $\mathcal I$ is the ideal generated by the relations (we use the notation $h=a^{-1}\gamma-\alpha^2-\beta^2 t$ here)
\[
\alpha\gamma=\beta\gamma=2\alpha=2\beta=a^2\beta^2 h=0.
\]
Furthermore, the algebra is $A_U=R_U[X]/(X^2=t+ahX)$, the element $\theta_U\in R_U$ is given by $\theta_U=\alpha+\beta\cdot X$ and the maps will be the ones from Table 2. The table is the following.
\begin{table}[ht]
\begin{center}
\begin{tabular}{|c||c|}
\hline
$\iota_U\colon R\to A,\;1\mapsto 1.$ & $\Phi_U\colon A\to A,\;\begin{cases}1\mapsto 1,\\X\mapsto \gamma-X.\end{cases}$\\
\hline
$\varepsilon_U\colon A\to R,\;1\mapsto 0,\,X\mapsto a.$ & $\cdot\theta_U\colon A\to A,\;\begin{cases}1\mapsto \alpha+\beta\cdot X,\\X\mapsto \beta t+(\alpha+a\beta h) \cdot X.\end{cases}$\\
\hline
\multicolumn{2}{|c|}{$m_U\colon A\otimes A\to A,\;\begin{cases}1\otimes1\mapsto 1,\,1\otimes X\mapsto X,\\X\otimes 1\mapsto X,\,X\otimes X\mapsto t+h\cdot X.\end{cases}$}\\
\hline
\multicolumn{2}{|c|}{$\Delta_U\colon A\to A\otimes A,\;\begin{cases}1\mapsto -h\cdot1\otimes 1+a^{-1}(1\otimes X+X\otimes 1),\\X\mapsto a^{-1}t\cdot 1\otimes 1+a^{-1}\cdot X\otimes X.\end{cases}$}\\
\hline
\end{tabular}
\caption{The maps for the generators from Figure~\ref{figure1-1}.}\label{tabular-matrixtable}
\end{center}
\end{table}
\end{prop}
\begin{proof}
We start by showing that the data given above give rise to a skew-extended Frobenius algebra, i.e. the satisfy the axioms given in Definition~\ref{defn-exfrob}. Note that the algebra $A_U$ is certainly a rank2-algebra over $R_U$, $\theta_U=\alpha+\beta X\in A_U$ and $(\varepsilon_U\circ\iota_U(1))=0$.

Moreover, it satisfies the axioms of an aspherical Frobenius algebra, since it, forgetting the new structure, coincides with the classical one given before Proposition 4 in~\cite{kh2}.

A direct computation verifies that $\Phi_U$ is a skew-involution, i.e.
\[
\Phi_U\circ \Phi_U=\mathrm{id_{A_U}},\;(\Phi_U\otimes\Phi_U)\circ\Delta_U\circ\Phi_U=-\Delta_U \;\text{and}\;\varepsilon_U\circ\Phi_U=-\varepsilon_U.
\]
Furthermore, a direct computation shows that $\cdot\theta_U$ and $\Phi_U$ also satisfy the axioms (a) and (b) from Definition~\ref{defn-exfrob}, i.e. the whole data is an aspherical rank2-uTQFT. 
\vspace*{0.5cm}

Now assume that we have a given aspherical rank2-uTQFT $\mathcal F=(R,A,\varepsilon,\Delta,\Phi,\theta)$.

First we observe that a skew-extended Frobenius algebra $A$ has an underlying Frobenius algebra of rank two. Hence, $\iota$ has to be of the given form. Because it is also aspherical, i.e. $\varepsilon(\iota(1))=0$, we see that $\varepsilon(1)=0$ and $\varepsilon(X)=a\cdot 1$. The element $a\in R$ is invertible because of the relation
\[
(\varepsilon\otimes\mathrm{id})\circ\Delta=\mathrm{id}=(\mathrm{id}\otimes\text{ }\varepsilon)\circ\Delta.
\]
It is known (e.g. see Proposition 5~\cite{kh2}) that such an algebra is of the form $A=R[X]/(X^2=t+ahX)$ with multiplication $m$ and comultiplication $\Delta$ from the Table 2 above.

We look at the new structure now. Because $\theta$ is an element of $A\cong R\cdot 1\oplus R\cdot X$ we find $\alpha,\beta\in R$ such that $\theta=\alpha+\beta X$. Using the multiplication we see that $X^2=t+ah\cdot X$. So an easy calculation shows that $\theta\cdot X=\beta t+(\alpha+a\beta h)X$ which gives us the map $\cdot\theta$ as above.

Because the map $\Phi\colon A\to A$ is not only $R$-linear, but also a skew-involution, we get $\Phi(1)=1$ and with $\varepsilon\circ\Phi=-\varepsilon$ we get $\Phi(X)=\gamma-X$ for an element $\gamma\in R$. Using then the first relation of a skew-extended Frobenius algebra we get the relations $\alpha\gamma=\beta\gamma=2\beta=0$ and $2(\alpha+a\beta h)=2\alpha=0$.

Using the second relation of a skew-extended Frobenius algebra, namely
\[
m\circ(\Phi\otimes\mathrm{id_{A}})\circ\Delta=(\cdot\theta)^2,
\]
we get the last two relations $ah=\gamma-a\alpha^2-a\beta^2t$ and $a^2\beta^2 h=0$.

These are all relations we get from the axioms of an aspherical rank2-uTQFT, i.e. any other axiom will also lead to one of these relations.
\end{proof}
\begin{rem}\label{rem-tutuextension}
The reader familiar with the paper of Turaev and Turner will recognise that our \textit{universal skew-extended Frobenius algebra} $\mathcal F_U$ (given in and before Proposition 2.15 of~\cite{tutu}) is different from the one from Turaev and Turner. But this is an advantage (see Corollary~\ref{cor-khovhomo}).
\vspace*{0.5cm}

As mentioned before in Remark~\ref{rem-tutu}, the version of Turaev and Turner can be obtained from our approach too. The difference again are the relations $\varepsilon^+=-\varepsilon^-$ and $\Delta^+_{++}=-\Delta^-_{--}$. This forces $\Phi=F(\Phi^-_+)$ from the proof above to send $X\mapsto\gamma -X$ instead of $X\mapsto\gamma +X$ (but over $R$ with $\mathrm{char}(R)=2$ they coincide).
\end{rem}
\vspace*{0.5cm}
The next corollary allows us to characterise the uTQFTs which lead to v-link homology.
\begin{cor}
\label{cor-bnraretrue}
Every aspherical rank2-uTQFT $\mathcal F$ satisfies the local relations from Fig.~\ref{figureintroa-4}.
\end{cor}
\begin{proof}
View a sphere $S^2$ as a cobordism $S^2\colon\emptyset\rightarrow\emptyset$. Then $\mathcal F(S^2)=\mathcal F(\varepsilon^+)\circ\mathcal F(\iota_+)$. So we calculate $\mathcal F(S^2)=0$. Because of the axiom (4) from Definition~\ref{defn-utqft}, this is true for every cobordism with a sphere. Analogously view a torus $\mathcal T$ as a cobordism $\mathcal T\colon\emptyset\rightarrow\emptyset$. Thus, it is of the form $\mathcal F(\mathcal T)=\mathcal F(\varepsilon^+)\circ \mathcal F(m^{++}_+)\circ\mathcal F(\Delta^+_{++})\circ\mathcal F(\iota_+)$. An easy calculation with the maps of Table 2 shows, that $\mathcal F(\mathcal T)=2$. Because of the axiom (4), this is true for every cobordism with a torus.

The \textit{4Tu}-relation is algebraically just the formula
\[
\Delta_{12}\circ\iota+\Delta_{34}\circ\iota=\Delta_{13}\circ\iota+\Delta_{24}\circ\iota.
\]
Here $\Delta_{ij}:A\rightarrow A\otimes A\otimes A\otimes A$ is the map which sends an element $a\in A$ to an element $a_1\otimes a_2\otimes a_3\otimes a_4$ with $a_k=1$ for $k\neq i,j$ and $a_i,a_j$ the first respectively the second tensor factor of $\Delta(a)$ (see Fig.~\ref{figure-4tu}).
\begin{figure}[ht]
    \centerline{\includegraphics[scale=0.525]{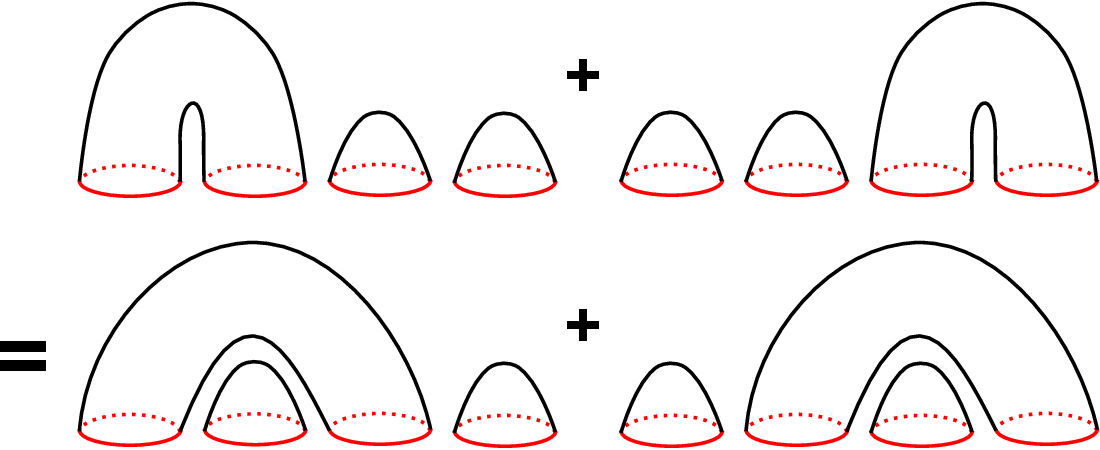}}
 \caption{The relation $\Delta_{12}\circ\iota+\Delta_{34}\circ\iota=\Delta_{13}\circ\iota+\Delta_{24}\circ\iota.$}
 \label{figure-4tu}
\end{figure}

That this relation is true is also an easy calculation. Again axiom (4) gives us the global statement. Because this is true for the universal skew-extended Frobenius algebra $\mathcal F_U$, we get the statement for all aspherical rank2-uTQFTs from Proposition~\ref{cor-univalg}.
\end{proof}
\vspace*{0.5cm}
Because with an aspherical, rank2 skew-extended Frobenius algebra we can define a corresponding rank2-uTQFT which satisfies the Bar-Natan relations, we note the following two corollaries.
\begin{cor}\label{cor-invariant}
Every aspherical, rank2 uTQFT can be used to define a v-link invariant.
\end{cor}
\begin{proof}
Directly from Theorem~\ref{thm-algcomplex}, Proposition~\ref{cor-univalg} and Corollary~\ref{cor-bnraretrue}.
\end{proof}
\vspace*{0.5cm}
\begin{cor}\label{cor-khovhomo}(\textbf{The virtual Khovanov complex}) The above construction enables one to extend the Khovanov complex ($R_{Kh}=\bZ,A_{Kh}=\bZ[X]/(X^2=0,t=h=0)$) from c-links to v-links by setting $\alpha=\beta=\gamma=0$ and $a=1$.
\end{cor}
\begin{proof}
Use Proposition~\ref{cor-univalg} and the substitution above.
\end{proof}
\vspace*{0.5cm}
From now on we denote by
\[
\Kh (L)=\mathcal F_{Kh}(\bn L)
\]
the \textit{virtual Khovanov complex} of a v-link $L$ and by $H(\Kh (L))$ its homology.
\begin{rem}
It is possible to introduce gradings (by setting $\deg 1=1$ and $\deg X=-1$) for the complex from Corollary~\ref{cor-khovhomo}. This is true because the map $\cdot\theta$ is equal to zero. In fact, this is the only possibility where we can introduce gradings, because all maps in the Khovanov complex must decrease the grading by one. And this is only possible if $\cdot\theta\colon A\rightarrow A$ is equal to zero.
\end{rem}
\vspace*{0.5cm}
\begin{cor}\label{cor-khovhomocat}(\textbf{Categorification of the virtual Jones polynomial}) The virtual Khovanov complex~\ref{cor-khovhomo} is a categorification of the virtual Jones polynomial in the sense that its graded Euler characteristic gives the polynomial.
\end{cor}
\begin{proof}
The classical Jones polynomial is uniquely determined by the skein-relations. The same is true for the virtual Jones polynomial, see for example Sec. 5 in~\cite{ka1} (although he uses the bracket polynomial $\langle\cdot\rangle$ instead of the skein relations, but this can be adopted straightforward as in the classical case). One can now easily check that the (graded) Euler characteristic of the virtual Khovanov complex $\Kh(\cdot)$ satisfies these relations.
\end{proof}
There is also an extension of the Khovanov-Lee complex (see Sec. 4.3 in~\cite{lee}) and two different extensions of Bar-Natan's variant ($R=\bZ/2,h=1,t=0$) (see Sec. 9.3 in~\cite{bn2}).
\vspace*{0.5cm}
\begin{cor}\label{cor-leekhovhomo}(\textbf{The virtual Khovanov-Lee complex}) The above construction enables us to extend the Khovanov-Lee complex ($R_{Lee}=\bZ$ and $A_{Lee}=\bZ[X]/(X^2=0,t=1,h=0)$) from c-links to v-links by setting $\alpha=\beta=\gamma=0$ and $a=1$.
\end{cor}
\begin{proof}
Directly from Proposition~\ref{cor-univalg}.
\end{proof}
\vspace*{0.5cm}
\begin{cor}\label{cor-bnhomo}(\textbf{The virtual Khovanov-Bar-Natan complex}) The above construction enables us to extend Bar Natan's variant of Khovanov homology (this is the Frobenius algebra over the field $R_{BN}=\bZ/2$ with $A_{BN}=\bZ/2[X]/(X^2=0,t=0,h=1)$) from c-links to v-links in two different ways by setting $\alpha=\beta=0$ and $\gamma=a=1$ or by setting $\beta=\gamma=0$ and $\alpha=a=1$. The two extensions are non-isomorphic skew-extended Frobenius algebras.
\end{cor}
\begin{proof}
That these two skew-extended Frobenius algebras can be used as v-link homologies follows from Proposition~\ref{cor-univalg} and Corollary~\ref{cor-invariant}. To see that they are non-isomorphic skew-extended Frobenius algebras we note that $\theta=0$ in the first case and $\theta=1$ in the second case. Because any isomorphism of skew-extended Frobenius algebras satisfies $f(1)=1$ and $f(\theta)=\theta^{\prime}$, they are not isomorphic.
\end{proof}
\vspace*{0.5cm}
We denote these three extensions by
\[
\mathcal F_{Lee} (L)=\mathcal F_{Lee}(\bn L),\hspace*{0.5cm}\mathcal F_{BN1}(L)=\mathcal F_{BN1}(\bn L)
\]
and
\[
\mathcal F_{BN2}(L)=\mathcal F_{BN2}(\bn L)
\]
respectively.
\begin{prop}\label{prop-classiccomplex}
Let $L_D$ be a c-link diagram and let $\mathcal F$ be an aspherical rank2-uTQFT. Then the complex $\mathcal F(\bn{L_D})$  is the classical Khovanov complex (up to chain isomorphisms) which is obtained by using the underlying TQFT $\mathcal F^{\prime}$ of $\mathcal F$.
\vspace*{0.5cm}

A similar statement is true for the Khovanov-Lee complex and the two different versions of the Khovanov-Bar-Natan complex.
\end{prop}
\begin{proof}
This is just the algebraic version of Theorem~\ref{thm-classic}.
\end{proof}
It is worth noting that, if $L$ is a c-link, then these three (and any other of the possibilities) are the classical complexes (up to chain homotopies) due to Theorem~\ref{thm-algcomplex}. Moreover, it should be noted that Corollary~\ref{cor-bnhomo} and Proposition~\ref{prop-classiccomplex} give a method to distinguish v-links that are not c-links. That is, a v-link diagram $L_D$ with
\[
H(\mathcal F_{BN1}(L_D))_*\not\cong H(\mathcal F_{BN2}(L_D))_*
\]
can not be a v-diagram of a c-link, since a c-link diagram does not need the map $\cdot\theta$.
\vspace*{0.5cm}

Because Khovanov showed (see the discussion after Corollary 2 in~\cite{kh2}) that every TQFT which respects the first Reidemeister move must have an underlying $R$-module $A\cong R\cdot 1\oplus R\cdot X$ for an element $X\in A$, we also get the following theorem.
\vspace*{0.5cm}
\begin{thm}\label{thm-class}(\textbf{Classification of aspherical uTQFTs}) The following statements are equivalent for an aspherical uTQFT.
\begin{itemize}
\item[(a)] It respects the first Reidemeister move RM1.
\item[(b)] It is a rank2-uTQFT.
\item[(c)] It can be obtained from the one of Proposition~\ref{cor-univalg}.
\item[(d)] It can be used as a v-link invariant.
\end{itemize} 
\end{thm}
\vspace*{0.5cm}
With the work already done the proof is simple.
\begin{proof}
(a)$\Rightarrow$(b): This was done by Khovanov and stays true.

(b)$\Rightarrow$(c): This is just the Proposition~\ref{cor-univalg}.

(c)$\Rightarrow$(d): This is the Corollary~\ref{cor-invariant}.

(d)$\Rightarrow$(a): This is clear.
\end{proof}
\vspace*{0.5cm}
\begin{rem}\label{rem-mant}
We conjecture that Manturov's $\bZ$-version~\cite{ma2} is a strictly weaker invariant than our extension of the Khovanov complex~\ref{cor-khovhomo} in the following sense. A v-link with ``lots'' of classical crossings is likely to have ``lots'' of faces of type 1b or its mirror image (see~\ref{figure-basic}). We call these faces the \textit{virtual trefoil faces}. In our construction the two multiplications (or comultiplications for the mirror image) are not the same, i.e. they have different boundary decorations as pictured for example in Fig.~\ref{figure-big2}, since we take extra information of this face in account. In contrast, in Manturov's version they are just the same maps. It is worth noting that we use the extra information explicitly in Sec.~\ref{sec-vkhapp}.
\end{rem}
\vspace*{0.5cm}
\begin{rem}\label{rem-otherrelations}
At this state it is a fair question to ask why we use the relations (1) from Equation~\ref{eq-combrel12} (or the one without the signs for the variant of Turaev and Turner) for our cobordisms, i.e. why do we assume that $\Delta^{+}_{++}$ changes its sign under conjugation with $\Phi^-_+$ and not $m^{++}_+$ (or neither of them changes its sign for the variant of Turaev and Turner).
\vspace*{0.5cm}

So what happens if we assume that $m^{++}_+$ changes its sign under conjugation with $\Phi^-_+$ (or both)? One can repeat the whole construction from Sec.~\ref{sec-vkhcat}, Sec.~\ref{sec-vkhcom} and this Sec.~\ref{sec-vkhfa} for these cases too. But this does not lead to anything new, i.e. if we assume that $m^{++}_+$ changes its sign, then we get an equivalent to the construction above and if we assume that both of them change their signs, then we get an equivalent to the variant of Turaev and Turner again.
\end{rem}
\vspace*{0.5cm}
\begin{rem}\label{rem-newthings}
Note that the classification of Theorem~\ref{thm-class} and Table 2 include non-classical invariants. To be more precise, if we work for example over $R=\bQ$, then the relations force us to set $\theta=0$. But if we work over $R=\bZ/2$, then we have different choices for $\theta$. It should be noted that, since c-links do not require the map $\cdot\theta$, these invariants can not appear in the classical setting. Note that in both of Manturov's versions~\cite{ma2} and~\cite{ma1} he sets $\theta=0$. 
\end{rem}
\section{The topological complex for virtual tangles}\label{sec-vkhtan}
We will define the \textit{topological complex of a v-tangle diagram} $T^k_D$ in this section. For this construction we use our notations for the \textit{saddle decorations} and \textit{saddle signs} of v-link diagrams $L_D$ from Sec.~\ref{sec-vkhcom}. Recall that a crucial ingredient for the construction of the topological complex were the \textit{decorations} of the saddles. Note that we work in a slightly different category now, i.e. the one from Definition~\ref{defn-category3}. Hence, we need signs, glueing numbers and indicators.
\vspace*{0.5cm}

It is worth noting that the idea how to find non-trivial solutions to the problems that come with the observation summarised in Fig.~\ref{figure0-order} (where non-trivial means that we do not define open saddles to be zero) is the following. Take the signs and decorations of a closure of the v-tangle diagram, since we already defined how to spread them for v-link diagrams in a ``good'' way. Note that this convention makes it easy to show analogous statements as in Sec.~\ref{sec-vkhcom}.
\vspace*{0.5cm}

We note that this section has two subsections. We define the topological complex of a v-tangle diagram with a *-marker and show that it is v-tangle invariant in the first part, i.e. in Definition~\ref{defn-geocomplex} and Theorem~\ref{thm-invarianz}. In the second part we discuss how the position of the *-marker has influence on the complex. We can show in Theorem~\ref{thm-marker} that in general the position of the *-marker gives rise to two different v-tangle invariants, but it agrees with the classical construction for c-tangles.
\subsection{The topological complex for virtual tangles}
We start by explaining how we are going to extend the important notions of saddle sign and decorations to v-tangle diagrams.
\vspace*{0.5cm}

Recall that $T^k_D$, as in Definition~\ref{defn-vtangle}, should denote a v-tangle diagram with $k\in\bN$ boundary points. Moreover, such diagrams should always have a *-marker on the boundary and let $\mathrm{Cl}(T^k_D)$ be the closure of the diagram. Recall that such diagrams come with x-markers.
\begin{defn}\label{defn-deco2}(\textbf{``Open'' saddle decorations})
Let $T^k_D$ be a v-tangle diagram with a *-marker on the boundary and let $\mathrm{Cl}(T^k_D)$ be the closure of the diagram. The \textit{saddle decorations} of the saddles of $T^k_D$ should be the ones induced by the saddle decorations of the closure. To be more precise.
\begin{itemize}
\item[(a)] The \textit{signs} of the saddles of $T^k_D$ should be the same as the signs of the corresponding saddles of $\mathrm{Cl}(T^k_D)$ as defined in Definition~\ref{defn-sign}.
\item[(b)] The \textit{indicators} of the saddles should be obtained from the corresponding saddles of $\mathrm{Cl}(T^k_D)$ as follows.
\begin{itemize}
\item Every orientable surface should carry an indicator $+1$ iff the number of upper boundary components of the saddle is two and a $-1$ iff the number is one.
\item Every non-orientable saddle gets a $0$ as an indicator.
\end{itemize}
\item[(c)] The \textit{glueing numbers} of the saddles of $T^k_D$ should be the same as the glueing numbers of the corresponding saddles of $\mathrm{Cl}(T^k_D)$ as defined in Definition~\ref{defn-deco}.
\end{itemize}
Note that saddles with a $0$-indicator do not have any boundary decorations. Everything together, i.e. boundary decorations, the saddle sign, and the indicator, is called the \textit{saddle decorations} of $S$.
\end{defn}
Beware again that many choices are involved. But they do not change the complex up to chain isomorphisms as we show in Lemma~\ref{lem-everythingfine} which is an analogon of Lemma~\ref{lem-commutativeindependence}.
\begin{defn}\label{defn-geocomplex}(Topological complex for v-tangles) For a v-tangle diagram $T^k_D$ with a *-marker on the boundary and with $n$ ordered crossings we define \textit{the topological complex} $\bn{T^k_D}$ as follows.
\begin{itemize}
\item For $i\in\{0,\dots,n\}$ the $i-n_-$ \textit{chain module} is the formal direct sum of all resolutions $\gamma_a$ of length $i$.
\item There are only morphisms between the chain modules of length $i$ and $i+1$.
\item If two words $a,a^{\prime}$ differ only in exactly one letter and $a_r=0$ and $a^{\prime}_r=1$, then there is a morphism between $\gamma_a$ and $\gamma_{a^{\prime}}$. Otherwise all morphisms between components of length $i$ and $i+1$ are zero.
\item This morphism is a \textit{saddle} between $\gamma_a$ and $\gamma_{a^{\prime}}$.
\item The saddles should carry the \textit{saddle decorations} from Definition~\ref{defn-deco2}.
\end{itemize}
\end{defn}
We note again that it is not clear at this point why we can choose the numbering of the crossings, the numbering of the v-circles, and the orientation of the resolutions of the closure. Furthermore, it is not clear why this complex is a well-defined chain complex. But we show in Lemma~\ref{lem-everythingfine} that the complex is independent of these choices, i.e. if $\bn{L_D}_1$ and $\bn{L_D}_2$ are well-defined chain complexes with different choices, then they are equal up to chain isomorphisms. The same lemma ensures that the complex is a well-defined chain complex.
\vspace*{0.5cm}

Another point that is worth mentioning is that the signs in our construction, in contrast to the classical Khovanov homology, do not depend on the order of the crossings of the diagram.

Beware that the position of the *-marker is important for v-tangle diagrams. But Theorem~\ref{thm-marker} ensures that the position is not important for c-tangles and v-links.

If it does not matter which of the two possibly different chain complexes is which, i.e. if it is just important that they could be different, then we denote them by $\bn{T^k_D}^*$ and $\bn{T^k_D}_*$ for a given v-tangle diagram $T^k_D$ without a chosen *-marker position. 
\vspace*{0.5cm}

\begin{figure}[ht]
     \centerline{\includegraphics[scale=0.425]{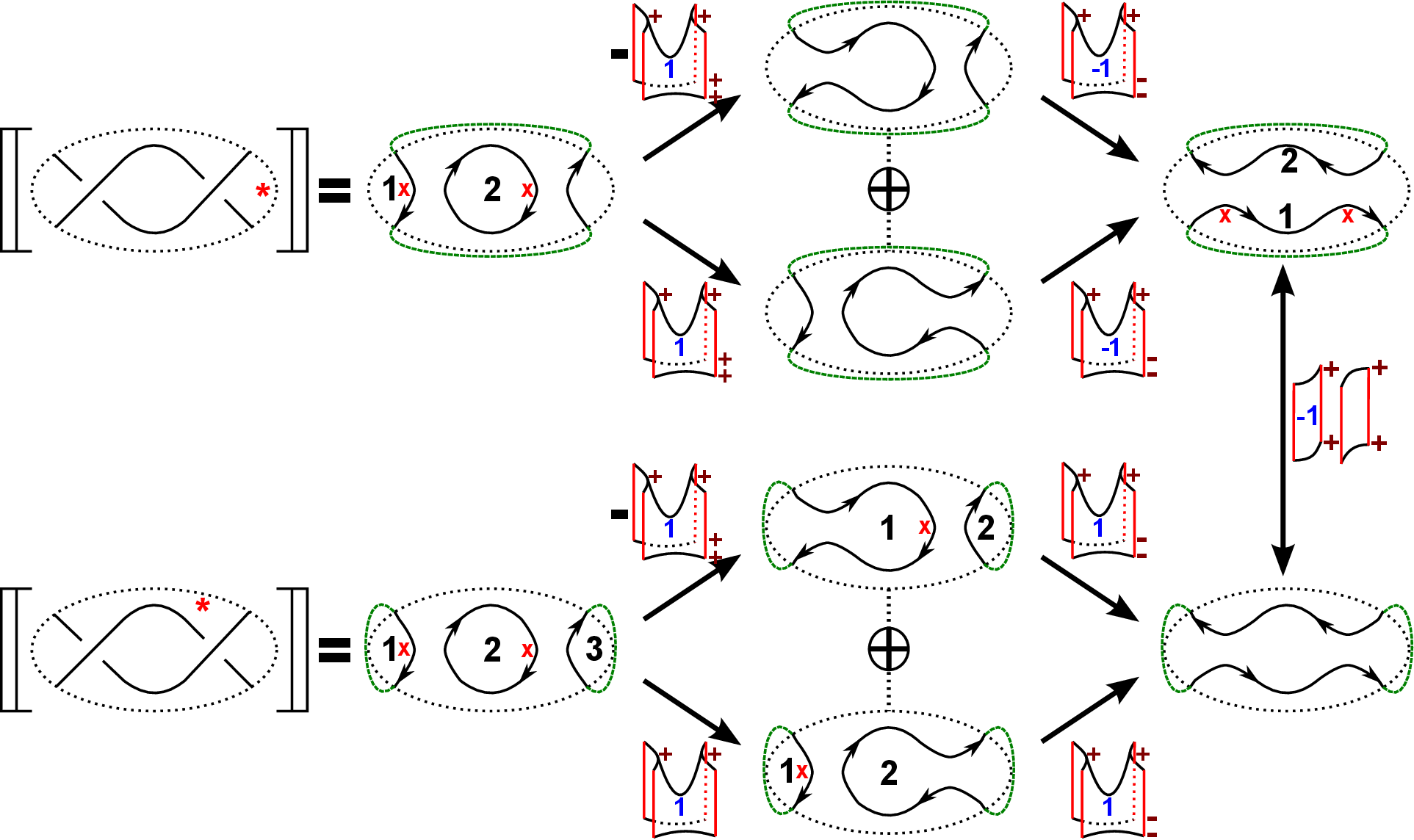}}
  \caption{The complexes of the same v-tangle with different *-marker positions. The two complexes are (in this case) isomorphic.}
  \label{figure-tanglecomplex}
\end{figure}

For an example see Fig.~\ref{figure-tanglecomplex}. This figure shows the virtual Khovanov complex of a v-tangle diagram with two different *-marker positions. The vertical arrow between them indicates that they are (in this case) chain isomorphic. It is worth noting at this point that, as we show in Theorem~\ref{thm-marker}, they are always isomorphic if the diagram is a c-tangle diagram (as the two diagrams in the figure below).
\begin{lem}\label{lem-everythingfine}
Let $T^k_D$ be a v-tangle diagram with a *-marker and let $\bn{T^k_D}_1$ be its topological complex from Definition~\ref{defn-geocomplex} with arbitrary orientations for the resolutions of the closure. Let $\bn{T^k_D}_2$ be the complex with the same orientations for the resolutions except for one circle $c$ in one resolution $\gamma_a$. If a face $F_1$ from $\bn{T^k_D}_1$ is anti-commutative, then the corresponding face $F_2$ from $\bn{T^k_D}_2$ is also anti-commutative.
\vspace*{0.5cm}

Moreover, if $\bn{T^k_D}_1$ is a well-defined chain complex, then it is isomorphic to $\bn{T^k_D}_2$, which is also a well-defined chain complex.

The same statement is true if the difference between the two complexes is the numbering of the crossings, the choice of the x-marker for the calculation of the saddle signs or the fixed numbering of the v-circles of the closure. Moreover, the same is true for any rotations/isotopies of the v-tangle diagram.
\end{lem}
\begin{proof}
For v-tangle diagrams $T^k_D$ with $k=0$ the statement is the same as the corresponding statements in Lemma~\ref{lem-commutativeindependence} and Corollary~\ref{cor-chaincomplex}. Recall that the trick is to reduce all faces through a finite sequence of vRM1, vRM2, vRM3 and mRM moves in Fig.~\ref{figure0-reide} and virtualisations from Fig.~\ref{figure0-virt} to a finite number of different possible faces. Then one does a case-by-case check.
\vspace*{0.5cm}

Because the saddles in the two chain complexes are topologically the same, we only have to worry about the decorations. But the decorations are spread based on the closure of the v-tangle diagram and the relations from Definition~\ref{defn-category3} are build in such a way that the open cases behave as the closed ones.

Hence, we can use the statement for $k=0$ to finish the proof, since the only possible differences for $k>0$ are the indicators, but they only depend on the *-marker.
\end{proof}
In the same vein as in Sec.~\ref{sec-vkhcom} we obtain the following Corollary.
\begin{cor}\label{cor-chaincomplex2}
The complex $\bn{T^k_D}$ is a chain complex. Thus, it is an object in the category $\ukobk_R$. 
\end{cor}
\begin{proof}
As in Sec.~\ref{sec-vkhcom}.
\end{proof}
Hence, we can speak of \textit{the} topological complex $\bn{T^k_D}$ of the v-tangle diagram with a *-marker. The complex is by Corollary~\ref{cor-chaincomplex2} a well-defined chain complex.
\vspace*{0.5cm}

The next theorem is very important but the proof itself is almost equal to the proof of Theorem~\ref{thm-geoinvarianz}. Therefore, we skip the details.
\begin{thm}\label{thm-invarianz}
Let $T^k_D,\tilde T^k_D$ be two v-tangle diagrams with the same *-marker position which differ only through a finite sequence of isotopies and generalised Reidemeister moves. Then the complexes $\bn{T^k_D}$ and $\bn{\tilde T^k_D}$ are equal in $\ukobk^{hl}_R$.
\end{thm}
\begin{proof}
We can copy the arguments of Theorem~\ref{thm-geoinvarianz}. Lemma~\ref{lem-everythingfine} guarantees that we can choose the numbering and orientations without changing anything up to chain isomorphisms.
\vspace*{0.5cm}

Beware that the chain homotopies in~\ref{thm-geoinvarianz} should all carry $+1$ as an indicator. Again, one can check that the involved chain homotopies satisfy the condition of a strong deformation retract.
\end{proof}
\subsection{The *-marker and the classical complex}
We need some notions now. Note that they seem to be ad-hoc, but the main motivation is that in general the position of the *-marker is important. But to recover at least some local properties, as discussed in Sec.~\ref{sec-vkhca}, we need to identify basic parts of v-tangle diagrams such that the two complexes are isomorphic.

\begin{figure}[ht]
     \centerline{\includegraphics[scale=0.56]{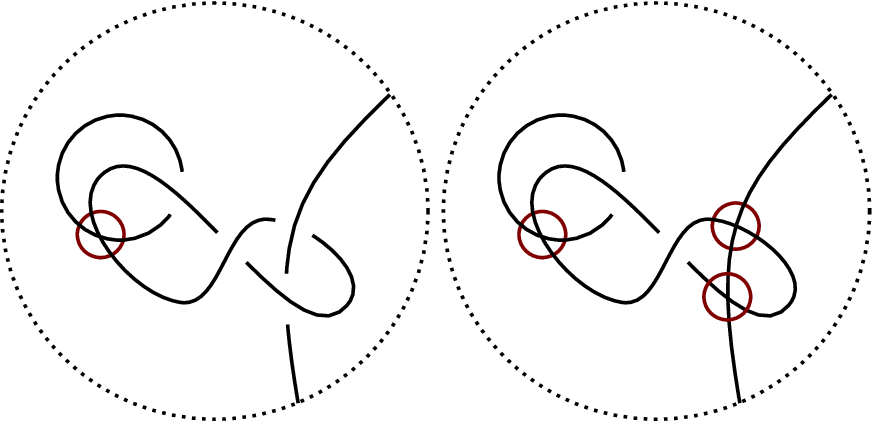}}
  \caption{The left v-tangle diagram is not fully internal, but the right diagram has a fully internal component (the two internal v-circles).}
  \label{figure-internal}
\end{figure}

Let $T^k_D$ denote a v-tangle diagram. We call a part of $T^k_D$ a \textit{connected part} if it is connected as the four-valent graph obtained by ignoring the v-crossings and treating c-crossings as vertices. We call a connected part of a v-tangle diagram \textit{fully internal} if it is not adjacent to the boundary, see Fig.~\ref{figure-internal}. The left v-tangle diagram has one connected part, which is not fully internal, and the right v-tangle diagram has two connected parts, one fully internal and one not fully internal.

A v-crossing is called \textit{negligible} if it is part of a fully internal component, e.g. all v-crossings of the right v-tangle diagram in Fig.~\ref{figure-internal} are negligible. Note that, by convention, negligible v-crossings are never part (for all resolutions) of any string that touches the boundary.

We call a v-tangle diagram $T^k_D$ \textit{nice} if there is a finite sequence of vRM1, vRM2, vRM3 and mRM moves and virtualisations which transforms $T^k_D$ into a diagram for which every v-crossing is negligible, e.g. every v-link diagram is nice and every c-tangle diagram is nice.

\begin{figure}[ht]
     \centerline{\includegraphics[scale=0.5]{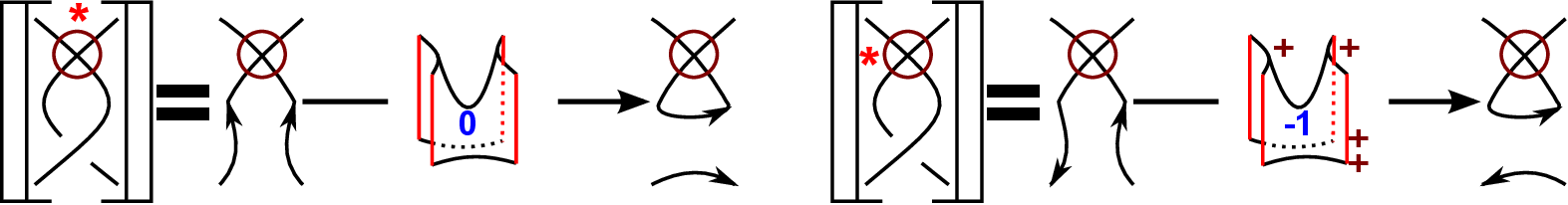}}
  \caption{A counterexample. The diagram is not a nice v-tangle diagram.}
  \label{figure-notnice}
\end{figure}

An example of a non-nice v-tangle diagram is shown in Fig.~\ref{figure-notnice}. Note that the complexes are not chain homotopic.

We note that for a v-tangle diagram $T^k_D$ the chain complexes $\bn{T^k_D}^*$ and $\bn{T^k_D}_*$ are ``almost'' the same, i.e. they have the same vertices, but possibly different edges (which are still in the same positions). The next lemma makes this observation precise.
\vspace*{0.5cm}

It is worth noting that the Khovanov cube of a v-tangle diagram with $n$ crossings has $2^{n-1}n$ saddles. We number these saddles and the numbering in the lemma below should be the same for the two complexes.
\begin{lem}\label{lem-factorbyisos}
Let $T^k_D$ be a v-tangle diagram with $n$ crossings. Let ${}_iS(\mathrm{in})^*$ and ${}_iS(\mathrm{in})_*$ denote the numbered saddles of $\bn{T^k_D}^*$ and of $\bn{T^k_D}_*$. If $T^k_D$ is a nice v-tangle diagram, then we have for all $i=1,\dots,2^{n-1}n$ a factorisation of the form ${}_iS(\mathrm{in})^*=\alpha\circ {}_iS(\mathrm{in})_*\circ\beta$ for two invertible cobordisms $\alpha,\beta$.
\end{lem}
\begin{proof}
It is clear that the saddles are topologically equivalent. So we only need to consider the decorations. The main point is the following observation. Of the four outer (the two leftmost and the two rightmost) cobordisms in the bottom row of Fig.~\ref{figure1-3}, i.e. $\mathrm{id}(1)^+_+$, $\Phi(1)^-_+$, $\mathrm{id}(0)$ and $\mathrm{id}(-1)^+_+$, only the third is not invertible. The first is the identity, the second and fourth are their own inverses. The third is not invertible because the $0$-indicator can not be changed to a $\pm 1$-indicator.

Note that neither the vRM1, vRM2, vRM3 and mRM moves nor a virtualisation change the indicator of a saddle cobordism. Hence, it is sufficient to show the statement for a v-tangle diagram with only negligible v-crossings. From the observation above it is enough to show that every saddle gets a $0$-indicator in one closure iff it gets a $0$-indicator in the other closure.

The only possible way that a saddle gets an indicator from the set $\{+1,-1\}$ for one closure and a $0$-indicator for the other closure is the rightmost case in Fig.~\ref{figure0-order}. But for this case the existence of a non-negligible v-crossing is necessary. Hence, we get the statement.    
\end{proof}
\begin{prop}\label{prop-nicechainiso}
Let $T^k_D$ be a v-tangle diagram. If $T^k_D$ is nice, then $\bn{T^k_D}^*$ and $\bn{T^k_D}_*$ are chain isomorphic.
\end{prop}
\begin{proof}
Let $T^k_D$ be a nice v-tangle diagram. Then Lemma~\ref{lem-factorbyisos} ensures that every saddle is the same, up to isomorphisms, in $\bn{T^k_D}^*$ and $\bn{T^k_D}_*$. Furthermore, Lemma~\ref{lem-everythingfine} ensures that both complexes are well-defined chain complexes. Hence, the number of signs of every face is odd (also counting the ones from the decorations).
\vspace*{0.5cm}

Thus, we can use a spanning tree argument to construct the chain isomorphism explicitly, i.e. start at the rightmost leaves of a spanning tree of the Khovanov cube and change the orientations of the resolutions at the corresponding vertices such that the unique outgoing edges of the tree have the same decorations in both cases (Lemma~\ref{lem-everythingfine} ensures that nothing changes modulo chain isomorphisms). Continue along the vertices of the spanning tree, but remove already visited leaves. This construction generates a chain isomorphism.
\vspace*{0.5cm}

Next repeat the whole process, but change the indicators and afterwards the signs. It is worth noting that Lemma~\ref{lem-everythingfine} ensures that the two processes will never run into ambiguities or problems and Lemma~\ref{lem-factorbyisos} ensures that they will generate chain isomorphisms.

The chain isomorphism that we need is the composition of the three isomorphisms constructed before. See for example Fig.~\ref{figure-tanglecomplex}. 
\end{proof}
\begin{thm}\label{thm-marker}(\textbf{Two different chain complexes}) Let $T^k_D$ be a v-tangle diagram with two different *-marker positions. Let $\bn{T^k_D}^*$ and $\bn{T^k_D}_*$ be the topological complexes for the two positions. Then the two complexes are equal in $\ukobk^{hl}_R$ if the v-tangle has $k=0$ or is a c-tangle.
\end{thm}
\begin{proof}
We can use the Proposition~\ref{prop-nicechainiso} above for a v-tangle diagram with $k=0$. Moreover, we can choose a diagram without virtual crossing for a c-tangle without changing anything up to chain homotopies, because of Theorem~\ref{thm-geoinvarianz}. Then we can use the Proposition~\ref{prop-nicechainiso} again.
\end{proof}
\begin{rem}\label{rem-disums}
Note that the whole construction can be done with an \textit{arbitrary} closure of a v-tangle diagram, i.e. cap off in any possible way without creating new c- or v-crossings. The direct sum of all possibilities is then a v-tangle invariant. Or one can even allow v-crossings and take direct sums over all possibilities again. But since both is inconvenient for our purpose, we do not discuss it in detail here. 
\end{rem}
\begin{rem}\label{rem-gradings2}
Again, we could use the Euler characteristic to introduce the structure of a grading on $\ucob_R(k)$ (and hence on $\ukobk_R$). The differentials in the topological complex from Definition~\ref{defn-geocomplex} have all $\deg=0$ (after a degree shift), because their Euler characteristic is $-1$. Then it is easy to prove that the topological complex is a v-tangle invariant under graded homotopy.
\end{rem}
\section{Circuit algebras}\label{sec-vkhca}
In the present section we describe the notion of a \textit{circuit algebra}. A circuit algebra is almost the same as a planar algebra, but we allow virtual crossings.

Planar algebras were introduced by Jones~\cite{vj} to study $II_1$-subfactors. In our setting, they were for example studied by Bar-Natan in the case of classical Khovanov homology, see e.g. Sec.5 in~\cite{bn2}. Hence, we can use most of his constructions in our context, too. A crucial difference is that we need to \textit{decorate} our \textit{circuit diagrams}. This is necessary because our cobordisms are also decorated.

We start the section with the definition of a (decorated) circuit diagram. In the whole section every v-tangle diagram should have a *-marker. We call a v-tangle diagram \textit{decorated} if it has an orientation, a number (same numbers are allowed), one coloured (green and red) dot  for each of its v-circle/v-string and we call a cobordism \textit{decorated} if it has gluing numbers and an indicator. In the following we use the notation $\omega^*$ to illustrate that we consider all possibilities for $k\in\bN$ together.

\begin{defn}\label{defn-circuit diagram}
Let $D^2_o$ denote a disk embedded into $\bR^2$, the so-called \textit{outside disk}. Let $\mathrm{I}_1,\dots,\mathrm{I}_{m-1}$ denote disks $D^2$ embedded into $\bR^2$ such that for all $k\in\{0,\dots,m-1\}$ the disk $\mathrm{I}_k$ is also embedded into $D^2_o$ without touching the boundary of $D^2_o$, i.e. $\mathrm{I}_k\subset D^2_o\subset \bR^2$, $\mathrm{I}_k\cap\mathrm{I}_{k^{\prime}}=\emptyset$ for $k\neq k^{\prime}$ and $\mathrm{I}_k\cap\partial D^2_o=\emptyset$. We denote $\mathcal D_m=D^2_o-(\mathrm{I}_0\cup\cdots\cup\mathrm{I}_{m-1})$. The $\mathrm{I}_k$ are called \textit{input disks}.

A \textit{circuit diagram with $m$ input disks} $\mathcal{CD}_m$ is a planar graph embedded into $\mathcal D_m$ with only vertices of valency one and four and such that every vertex of valency one is in $\partial\mathcal D_m$ and every vertex of valency four is in $\mathrm{Int}(\mathcal D_m)$. All vertices of valency four are marked with a v-crossing. Again we allow circles, i.e. closed edges without any vertices. A \textit{*-marked circuit diagram} is the same, but with $m+1$ extra *-markers for every boundary component of $\mathcal D_m$. Moreover, we call the vertices at $\partial D^2_o$ the \textit{outer boundary points}.

See for example Fig.~\ref{figure-deccircuit}, i.e. the figure shows a *-marked (decorated) circuit diagram with three input disks.
\begin{figure}[ht]
     \centerline{\includegraphics[scale=0.45]{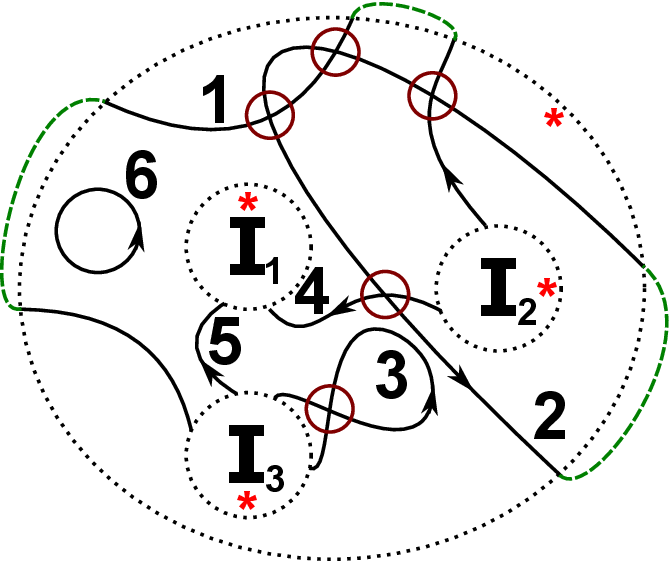}}
  \caption{A decorated circuit diagram with three input disks.}
  \label{figure-deccircuit}
\end{figure}
\end{defn}

A \textit{closure} of a *-marked circuit diagram with $m$ input disks $\mathrm{Cl}(\mathcal{CD}_m)$ is a circuit diagram with $m$ input disks and without any outer boundary points which is constructed from $\mathcal{CD}_m$ by capping off neighbouring strings starting from the \textit{outer} *-marker and proceeding counterclockwise. Note that we only cap off the outside disk and not the small inside disks.

A \textit{decoration} for a *-marked circuit diagram is a tuple of a numbering and an orientation of the strings of the diagram in such a way that it is also a numbering and orientation of the closure. We call a circuit diagram together with a decoration a \textit{decorated circuit diagram}. See for example Fig.~\ref{figure-deccircuit}. The decoration of the circuit diagram in this figure is also a decoration for the closure (the diagram together with the green lines).

\begin{figure}[ht]
     \centerline{\includegraphics[scale=0.45]{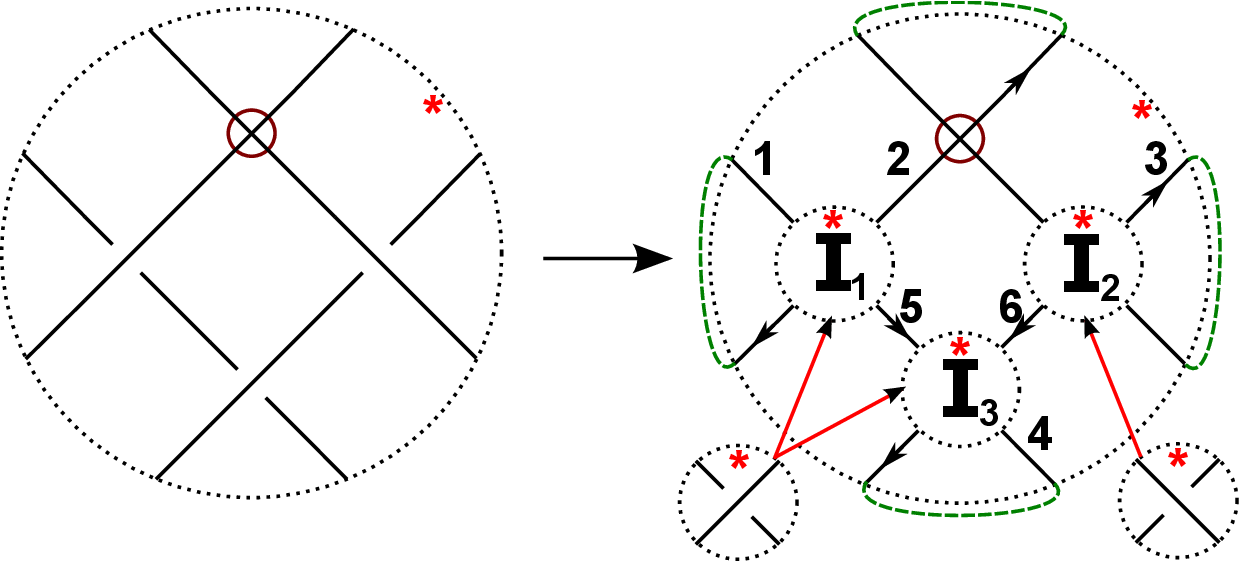}}
  \caption{A decorated circuit diagram induced by a v-tangle.}
  \label{figure-tanglealg}
\end{figure}

We can state now the definition of a \textit{(decorated) circuit algebra} with these notions. Recall that our v-tangle diagrams should always be oriented with the usual orientations but we suppress these again to maintain readability.
\begin{defn}\label{defn-circuit algebra}(\textbf{Circuit algebra}) Let $\mathrm{T^{\prime}}(k)$ be the set of (decorated) v-tangle diagrams with $k$ boundary points and a *-marker and let $\mathrm{T}(k)$ denote the quotient by boundary preserving isotopies and generalised Reidemeister moves.

Furthermore, let $\mathcal{CD}_m$ denote a (decorated) circuit diagram with $m$ input disks and $k^{\prime}$ outer boundary points such that the $j$-th input disk has $k_j$ numbered boundary points. 

Because $\mathcal{CD}_m$ has no c-crossings, this induces operations
\[
\mathcal{CD}_m\colon \mathrm{T^{\prime}}(k_0)\times\cdots\times \mathrm{T^{\prime}}(k_{m-1})\to \mathrm{T^{\prime}}(k^{\prime})\text{ and }\mathcal{CD}_m\colon \mathrm{T}(k_0)\times\cdots\times \mathrm{T}(k_{m-1})\to \mathrm{T}(k^{\prime})
\]
by placing the $i$-th v-tangle diagram from $\mathrm{T^{(\prime)}}(k_i)$ in the $i$-th boundary component of $\mathcal{CD}_m$, i.e. glue the v-tangle inside in such a way that the *-markers match. See the right side of Fig.~\ref{figure-tanglealg}. 

There is an identity operation on $\mathrm{T^{(\prime)}}(k)$ (it is of the form $\jpg{5mm}{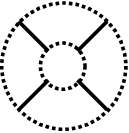}$) and the operations are compatible in a natural way (``associative''). We call a set of sets $\mathfrak C(\omega^*)$ with operations $\mathcal{CD}_m$ as above a \textit{circuit algebra}, provided that the identity and associativity from above hold.

If the operators and elements are decorated, first with numbers and orientations and latter with any kind of suitable decorations, then we call a set of sets $\mathfrak C(\omega^*)$ as before a \textit{decorated circuit algebra}. Note that in this case we have to define how the decorations change after glueing, that is, we can run into ambiguities if the glueing is defined in a non-compatible way.
\end{defn}
We should note that Fig.~\ref{figure-tanglealg} below also illustrates how a v-tangle diagram induces a decorated circuit diagram (choices for the decorations are involved).
Here are some examples. The reader may also check the ``classical'' examples in Sec. 5 of~\cite{bn2}.
\begin{ex}\label{ex-circuit}
The first example is the set $\Ob(\ucob(\omega^{*}))$ from Definition~\ref{defn-category}, i.e. v-tangle diagrams with $k\in\bN$ boundary points, an extra *-marker, but without c-crossings. This is a sub-circuit algebra of the circuit algebra that allows c-crossings.

But we want to view it as a decorated circuit algebra, denoted by $\Ob_d(\ucob(\omega^{*}))$, i.e. the elements are \textit{decorated} v-tangle diagram (all possible decorations). We have to define the operations in more detail now, since we can run into ambiguities, see top row of Fig.~\ref{figure-gluedeco}.
\begin{figure}[ht]
     \centerline{\includegraphics[scale=0.6]{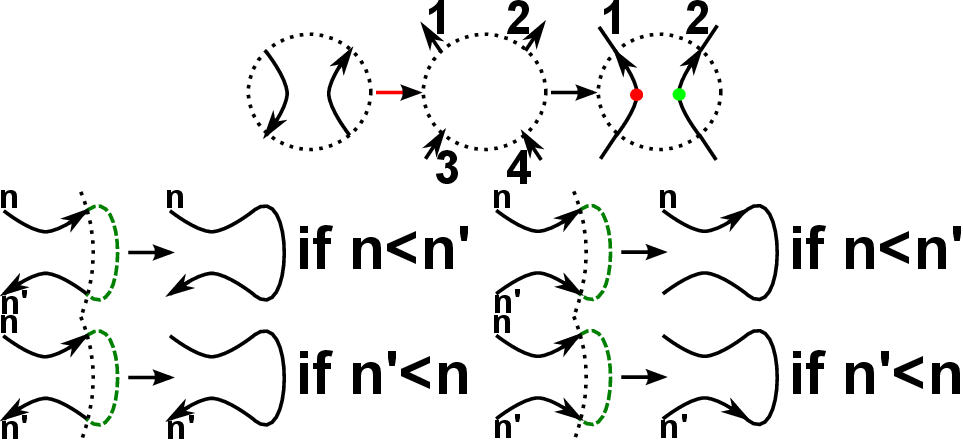}}
  \caption{The operation in the decorated circuit algebra.}
  \label{figure-gluedeco}
\end{figure}

First, we can run into ambiguities if the decorations of the operators that are glued together do not match. In this case we define the new decoration based on the rule \textit{``lower (number) first''}, i.e. the new number is the lower and the new orientation is the one from the lower numbered string. See the two lower rows of Fig.~\ref{figure-gluedeco}. Not all four cases are pictured, but we hope that it should be clear how the other two work. Moreover, it is worth noting that the order of these local steps does not affect the end result, since, by construction, the lowest number of all strings that are connected and its orientation will always determine the output. 

Furthermore, if we glue a decorated v-tangle diagram in an input disk, then we run into ambiguities if the shared decorations, i.e. the orientations and numbers, do not match. In this case we change the decorations of the v-tangle diagram (as above). We add in a \textit{red dot $r$} if we have to change the orientation and a \textit{green dot $g$} otherwise. This is pictured in the top row of Fig.~\ref{figure-gluedeco}. As above, in order to make it well-defined, one has to allow the dots to change stepwise, i.e. one uses the ``multiplication'' rules $g\cdot g=g=r\cdot r$ and $g\cdot r=r=r\cdot g$ for dots on the same string.

The reader should check that this gives rise to a (well-defined!) decorated circuit algebra.

Another important example is the whole collection $\Mor_d(\ucob(\omega^{*}))$ from Definition~\ref{defn-category}, i.e. decorated cobordisms (all possible decorations) with $k\in\bN$ vertical boundary lines and an extra *-marker. We want to view this example as a decorated circuit algebra again.

Hence, we have to define the operations. The most important point is the question how to handle the decorations again, because it should be clear how to glue a cobordism with $m$ vertical boundary lines into $\mathcal{CD}_m\times [-1,1]$. This time we have to define the behaviour of two decorations, i.e. the glueing numbers and indicators. The glueing numbers are treated as the orientations before, i.e. by using the ``lower first'' rule. The indicators (recall that they are just numbers of $\{0,+1,-1\}$) are multiplied. Recall that a cobordism with a $0$-indicator does not get any glueing numbers. We simply remove them in this case. To be more precise, we make the following definition of the operation of $\mathcal{CD}_m$ on $\Ob_d(\ucob(\omega^{*}))$ and $\Mor_d(\ucob(\omega^{*}))$ (compare to Fig.~\ref{figure-gluedeco}).
\begin{itemize}
\item A cobordism with a $+$ glueing number (or $-$) is composed with $\Phi^-_+$ iff the decorated v-tangle diagram (short: diagram) gets a red dot (or green) at the corresponding position.
\item A cobordism is composed with a $0$-indicator surface iff the strings of the diagram get identified at the bottom and top resolution at the corresponding position.
\item A cobordism with a $1$-indicator is composed with a $1$/$-1$-indicator surface iff the strings of the diagram get identified at the bottom/top resolution at the corresponding position.
\item A cobordism with a $-1$-indicator is composed with a $1$/$-1$-indicator surface iff the strings of the diagram get identified at the top/bottom resolution at the corresponding position.
\end{itemize} 
These rules define a new decoration for the new cobordism. The reader should check again that this gives rise to a (well-defined!) decorated circuit algebra (to see this we note that everything behaves multiplicative as the elements of $\{+1,-1\}\cong\bZ/2\bZ$ or has a $0$-indicator).
\end{ex}
We summarise the notions in a definition. Recall that v-tangle diagrams are decorated with orientations, numbers and coloured dots and cobordisms have glueing numbers and an indicator.
\begin{defn}\label{defn-dotcalculus}(\textbf{Dot-calculus}) Let $\mathcal{CD}_m$ denote a decorated circuit diagram with $m$ input disks and $k^{\prime}$ outer boundary points such that the $j$-th input disk has $k_j$ numbered boundary points. Then $\mathcal{CD}_m$ induces an associative and unital (as above) operation on decorated v-tangle diagrams (with a corresponding number of boundary points) by the ``lower first''-rule, i.e. if the orientation does not match, then the lower number induces the new orientation. Put a red dot $r$ on every string that has its orientation changed and a green dot $g$ otherwise (two dots on the same v-string are multiplied by the convention $g=1,r=-1$). We call this the \textit{v-tangle dot-calculus}.

Moreover, $\mathcal{CD}_m$ induces an associative and unital (as above) operation on decorated cobordisms (with a corresponding number of boundary lines) by the ``lower first''-rule, i.e. if the orientation does not match, then the lower number induces the new orientation. Put a red dot on every string that has its orientation changed and a green dot otherwise and compose the corresponding boundary with a $+$ glueing number (or $-$) with $\Phi^-_+$ iff the string has a red dot (or a green dot), multiply indicators via identity surfaces with corresponding indicators $+1$/$-1$ iff the v-tangle numbers get identified at the bottom/top (or vice versa for surfaces with a $-1$-indicator) resolution at the corresponding position, multiply with an $0$-identity iff in both resolutions the strings are identified (do everything repeatedly using the rules as explained above). We call this the \textit{dot-calculus}. The reader should compare the notions above with Fig.~\ref{figure0-main}.
\end{defn}
Recall that we assume that all v-tangle diagrams have already a fixed *-marker. A v-tangle diagram $T^k_D$ gives rise to a decorated circuit diagram $\mathcal{CD}_{T^k_D}$ as already illustrated in Fig.~\ref{figure-tanglealg}. If the diagram has $m$ crossings, denoted by $cr_1,\dots,cr_m$, then we choose a neighbourhood of the $cr_i$ without any other crossings and a *-marker for all $cr_i$. We obtain by this procedure $m$ v-tangle diagrams with one crossings and four boundary components, denoted by an abuse of notation by $cr_1,\dots,cr_m$, and we call these \textit{crossing diagrams associated to $T^k_D$}.
\begin{defn}\label{defn-assocomp}
Let $T^k_D$ be a v-tangle diagram with $m$ crossings and let $\mathcal{CD}_{T^k_D}$ and $cr_1,\dots,cr_m$ be its associated decorated diagram and crossings. Then the \textit{tensored complex}
\[
\mathcal{CD}_{T^k_D}(cr_1,\dots,cr_m)=(C_*,c_*)
\]
is defined as follows. Let $(C_j,c_j)$ with $j\in\{1,\dots,m\}$ be the topological complex of the $cr_j$ and defined in Definition~\ref{defn-geocomplex} such that the unique saddle is of the form $c_j\colon\smoothing\to\hsmoothing$ for any suitable orientation (without a sign). Let $\alpha_i,\beta_i$ denote the compositions of the morphisms that we compose after applying the circuit diagram on cobordisms (see Example~\ref{ex-circuit} and Definition~\ref{defn-dotcalculus} above), i.e. the red dots induce a composition with $\Phi^-_+$ (or with a $0$-identity surface in the degenerated case) and a change in the numbering induces a composition with a cobordisms that changes indicators.

Therefore, we denote the operation of $\mathcal{CD}_{T^k_D}$ on cobordisms, i.e. the dot-calculus, by $\alpha\circ\mathcal{CD}_{T^k_D}\circ\beta$ to illustrate the difference to the classical case. We do not use this notation for the objects to maintain readability. The $i$-th chain module is
\[
C^i=\bigoplus_{i=j_0+\dots+j_{m-1}}\mathcal{CD}_{T^k_D}(C_0^{j_0},\dots,C_{m-1}^{j_{m-1}})
\]
and the differentials are
\[
c|_{\mathcal{CD}_{T^k_D}(C_0^{j_0},\dots,C_{m-1}^{j_{m-1}})}=\sum_{i=0}^{m-1}\alpha\circ\mathcal{CD}_{T^k_D}(\mathrm{Id}_{C_0^{j_0}},\dots,c_i,\dots,\mathrm{Id}_{C_{m-1}^{j_{m-1}}})\circ\beta.
\]
Note that this complex does not have any extra signs and will in general not be a chain complex.
\end{defn}
We call a Khovanov cube \textit{of type p}, if all its faces are commutative up to a unit of $R$, and a projectivisation of such a cube is given by identifying morphisms up to units. We denote the latter usually with a superscript $P$. Details are in Sec.~\ref{sec-techcube}.

It should be noted that the choice of the *-markers in the definition of $\mathcal{CD}_{T^k_D}(cr_1,\dots,cr_m)$ or the choice of the decorations for $\mathcal{CD}_{T^k_D}$ is not important for our purpose (and we will suppress the difference). To be more precise, we give the following lemma. We should note that it is not clear at this point why the complexes are $m$ cubes of type p. But we show it in Theorem~\ref{thm-pcomplex} below.
\begin{lem}\label{lem-pcomplex}
Let $\mathcal{CD}_{T^k_D}(cr_1,\dots,cr_m)$ and $\mathcal{CD}^{\prime}_{T^k_D}(cr_1,\dots,cr_m)$ denote two different choices for the *-markers of the $cr_j$. Then the two complexes are equal.

Moreover, if the difference between $\mathcal{CD}_{T^k_D}(cr_1,\dots,cr_m)$ and $\mathcal{CD}^{\prime}_{T^k_D}(cr_1,\dots,cr_m)$ is the choice of decorations, either for the circuit diagram or for the saddles of the complexes of $cr_j$, then the two complexes are isomorphic as $m$ cubes of type p.
\end{lem}
\begin{proof}
This is the case because the $c_j$ has always an indicator $+1,-1$ in the definition of the complex $(C_j,c_j)$ and never a $0$-indicator. Moreover, the result depends only on the position of the *-marker for $T^k_D$, since the involved operations only depend on how strings are connected.

The second statement can be verified analogously as in Lemma~\ref{lem-everythingfine}.
\end{proof}
By a slight abuse of notation, we denote the topological complex by $\bn{T^k_D}$, although some choices are involved (but they do ``not matter'', see Lemma~\ref{lem-everythingfine}).
\begin{thm}\label{thm-pcomplex}(\textbf{Semi-locality I}) Let $T^k_D$ be a v-tangle diagram with $m$ crossings. Let $\bn{T^k_D}$ be (one of) its topological complex(es) from Definition~\ref{defn-geocomplex} and let $\mathcal{CD}_{T^k_D}(cr_1,\dots,cr_m)$ be its tensored complex from Definition~\ref{defn-assocomp}. Then $\mathcal{CD}_{T^k_D}(cr_1,\dots,cr_m)$ is a $m$-cube of type p and
\[
\mathcal{CD}_{T^k_D}(cr_1,\dots,cr_m)=\bn{T^k_D}^P
\]
for a suitable choice of orientations for the resolutions of $\bn{T^k_D}$.
\end{thm}
\begin{proof}
This is true because the dot-calculus is exactly built in such a way that the resulting saddles have some glueing numbers induced by a suitable choice of orientations of the resolutions. To be more precise, it is clear that the construction from Definition~\ref{defn-assocomp} gives rise to a $m$-cube as explained in Sec.~\ref{sec-techcube}.

Moreover, since we do not spread any formal signs in the construction from Definition~\ref{defn-assocomp}, the only thing we can expect is that the corresponding cube will be of type p, i.e. faces commute up to a sign. So we only have to care that the glueing numbers and indicators work out as claimed.

That the glueing numbers work out follows from the definition of the dot-calculus, since the orientation of the lowest numbered string will always determine the result and the decorations of the circuit diagram are also decorations of the closure, i.e. we can use Theorem~\ref{theo-facescommute} to see that the glueing numbers work out as claimed (up to a formal sign).

Moreover, the indicators of the saddles are spread based on a topological information, namely how certain strings are connected in the closure of the diagram $T^k_D$. Hence, since we have fixed the *-marker positions, these indicators are the same for $\mathcal{CD}_{T^k_D}(cr_1,\dots,cr_m)$ and any of the $\bn{T^k_D}$. Note that it is important that the indicators at the beginning are all $+1,-1$, since we can not change a $0$ using the conventions above.

This proves the statement, since there is a choice of orientations of the resolutions such that all saddles of $\mathcal{CD}_{T^k_D}(cr_1,\dots,cr_m)$ and $\bn{T^k_D}$ are equal up to a sign. 
\end{proof}
Given a Khovanov cube, then an \textit{edges assignment (with signs)} of this cube is a choice of extra signs for some of the saddles. We denote such an assignment using $\epsilon$ as a superscript, see Definition~\ref{defn-edge}.
\begin{cor}\label{cor-pcomplex}
There is an edge assignment such that $\mathcal{CD}^{\epsilon}_{T^k_D}(cr_1,\dots,cr_m)$ is a chain complex. Moreover, there is a chain isomorphism between $\mathcal{CD}^{\epsilon}_{T^k_D}(cr_1,\dots,cr_m)$ and $\bn{T^k_D}$ (for all possible choices involved in the definition of latter).  
\end{cor}
\begin{proof}
The first statement follows from Theorem~\ref{thm-pcomplex}, Theorem~\ref{theo-facescommute} and Lemma~\ref{lem-edge}. The second from Lemma~\ref{lem-everythingfine}.
\end{proof}
We note that Theorem~\ref{thm-pcomplex} and Corollary~\ref{cor-pcomplex} allow us to be ``sloppy'' when it comes to signs.

It is a natural question if one can generalise the statement of Theorem~\ref{thm-pcomplex}, since in the classical case one can allow arbitrary c-tangle diagrams as inputs. In fact, we do not know the answer in general. The main problem is that ``non-orientablity'' is not a local property.

We can make an analogous definition as in Definition~\ref{defn-assocomp}, but we allow the $cr_1,\dots,cr_m$ to be \textit{non-nice} v-tangle diagrams with one crossing. We denote them by $cr^{\prime}_1,\dots,cr^{\prime}_m$ to illustrate the difference and we call the corresponding complex \textit{generalised tensored complex}. An example is shown in Fig.~\ref{figure-counter}. Even this slight generalisation has unsatisfying properties.
\begin{thm}\label{thm-semiloc}(\textbf{Semi-locality II}) Let $T^k_D$ be a v-tangle diagram with $m$ crossings. And let
\[
\mathcal{CD}_{T^k_D}(cr^{\prime}_1,\dots,cr^{\prime}_m)=(C_*,c_*)
\]
be its generalised tensored complex.
\begin{itemize}
\item[(a)] The complex $(C_*,c_*)$ is a complex of type p, i.e. faces commute up to a unit of $R$.
\item[(b)] Let $\bn{T^k_D}$ denote (one of) its topological complex. Then we do not have a suitable choice for $\bn{T^k_D}$ in general such that
\[
\mathcal{CD}_{T^k_D}(cr^{\prime}_1,\dots,cr^{\prime}_m)=\bn{T^k_D}^P.
\]
\item[(c)] The complexes $(C_*,c_*)$ and $\bn{T^k_D}$ are not p-homotopic (see Definition~\ref{defn-chainhomotopy2}) in general.
\end{itemize}
\end{thm}
\begin{proof}
(a) This statement can be verified analogously to Theorem~\ref{thm-pcomplex}, since, if the corresponding saddles have an $+1,-1$ indicator, as in Theorem~\ref{thm-pcomplex}, then one can copy the arguments from before. If it has a $0$-indicator, then the arguments are even easier to verify, since we do not need any decorations in this case.

(b)+(c) This is true, because a surfaces with a $0$-indicator can not be changed to a surfaces with a $\pm 1$-indicator, since indicators behave multiplicatively. For an explicit example see Fig.~\ref{figure-counter}, i.e. the two complexes are not p-homotopy equivalent, since we can not change the $0$-indicator.
\begin{figure}[ht]
     \centerline{\includegraphics[scale=0.41]{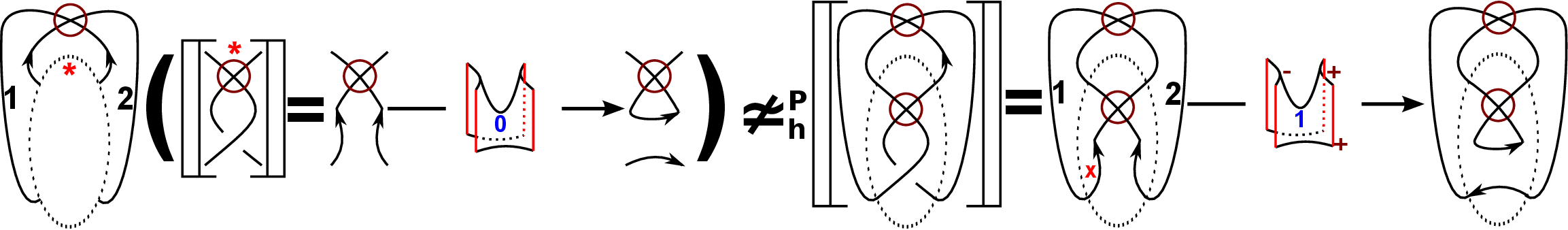}}
  \caption{A counterexample. The diagram is not a nice v-tangle diagram.}
  \label{figure-counter}
\end{figure}

Note that this includes that no choice will make them equal as complexes of type p.
\end{proof}
It should be noted again that the whole discussion in this section could be done with oriented (in the usual sense) v-tangle diagrams and oriented (decorated) circuit algebras. But to maintain readability we only refer to Sec. 5 in~\cite{bn2}, i.e. the reader can adopt the notions there and use them in our context without any difficulties.
\begin{rem}\label{rem-mod2}
It should be noted that the constructions presented in this section can be extended relatively easily to work in an even better way if one works over rings of characteristic two, e.g. over the ring $R=\bZ/2\bZ$.

This is the case because all appearing problems are in some sense ``sign problems''. If one works over $R=\bZ/2\bZ$, then, for example, the indicators are not necessary and most constructions will work analogously to the classical case (see e.g. Sec. 5 in~\cite{bn2}).
\end{rem}
\begin{rem}\label{rem-fastcalc}
One application of the local construction in the classical case is a way to calculate the classical Khovanov homology of a c-link with $n$ crossings in approximately $2^{\sqrt{n}}$ operations instead of $2^n$ operations of the ``brute force method'', see~\cite{bn3}. But in view of Theorem~\ref{thm-semiloc}, one has to be very careful if one tries to copy the method given in Sec.4 of~\cite{bn3}. 
\end{rem}
\section{An application: Degeneration of Lee's variant}\label{sec-vkhapp}
This section splits into three subsections. We explain the main motivations in the first and we are going to show that some facts about the classical Khovanov-Lee link homology are still true in the context of v-links (e.g. see Theorem~\ref{thm-leedeg}) in the last subsection. In order to do so, we identify the two generators with so-called non-alternating resolution~\ref{thm-nonalternating} in the second subsection. We note that these correspond to colourings in the c-case.
\vspace*{0.5cm}

The approach (we follow the idea of~\cite{bnsm}) to show that the degeneration of Khovanov-Lee's variant (see Subsection~\ref{sec-khlee-main} below) is still true is the following. First we define two orthogonal idempotents in our category, which we call \textit{down and up}. Then we can go to the \textit{Karoubi envelope} of our category, denoted by $\KAR(\ukobk_R)$.

The idea of the Karoubi envelope is to find a ``completion'' of a category such that every idempotent splits. It is named after the french mathematician Karoubi, but it already appears in an earlier work by Freyd~\cite{pf}. Note that it is sometimes called \textit{idempotent completion}. Then we show that the topological complex of a single crossing (as a v-tangle), if considered in $\KAR(\ukob_R(k))$, is homotopy equivalent to a very simple complex with only $0$-morphisms. After that we use the semi-local constructions from Sec.~\ref{sec-vkhca} to prove in Theorem~\ref{thm-leedeg} the virtual analogon of Lee's theorem.
\vspace*{0.5cm}

In the whole section let $R$ denote a commutative and unital ring such that $2$ is invertible, e.g. $R=\bZ\left[\frac{1}{2}\right]$. Moreover, throughout the whole section, we denote the topological complex by $\bn{\cdot}$ and its algebraic version by $\mathcal F(\bn{\cdot})$ or short by $\mathcal F(\cdot)$, e.g. we denote Khovanov-Lee's version by
\[
\mathcal F_{\mathrm{Lee}}(\cdot)=\mathcal F(\bn{\cdot}_{\mathrm{Lee}}).
\]
But, in order for the signs to work out correctly, we have to fix x-marker positions. In the whole section we, by convention, say that the x-marker for $\slashoverback$ is at the left side and for $\backoverslash$ the x-marker should be on the top.
\vspace*{0.5cm}

\begin{figure}[ht]
      \centerline{\includegraphics[width=0.55\linewidth]{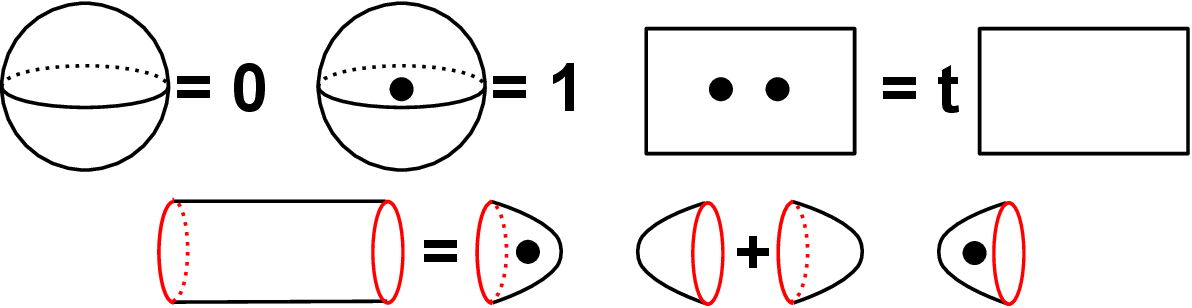}}
  \caption{The dot-relations. A dot is a shorthand notation for $\frac{1}{2}$-times a handle.}
  \label{figure-dotrel}
\end{figure}

Moreover, recall that the topological picture of Khovanov-Lee's variant is given by the \textit{dot-relations} in Fig.~\ref{figure-dotrel} with $t=1$, while the graded case of the Khovanov complex is $t=0$. Recall (see Proposition 11.1 in~\cite{bn1}) that $\frac{1}{2}\in R$ allows us to use the \textit{dot-relation} in Fig.~\ref{figure-dotrel} instead of the local relations of Fig.~\ref{figureintroa-4}. We give an example of the Khovanov-Lee complex of a v-knot in Example~\ref{ex-leedeg}.

\subsection{Main observations}\label{sec-khlee-main}
Recall (see Sec. 4.3 in~\cite{lee} for the classical and~\ref{cor-khovhomo} for the virtual case) that \textit{Khovanov-Lee's variant} for v-links is given by the filtered algebra $A=A_{\mathrm{Lee}}=R[X]/(X^2=1)$ and the following maps.
\[
m^{++}_+\colon A\otimes A\to A,\;\begin{cases}1\otimes 1\phantom{.}\mapsto 1, & X\otimes X\mapsto 1, \\ 1\otimes X\mapsto X, & X\otimes 1\phantom{.}\mapsto X\end{cases}
\]
for the multiplication and
\[
\Delta_{++}^+\colon A\to A\otimes A,\;\begin{cases}1\phantom{.}\mapsto 1\otimes X+X\otimes 1,\\ X\mapsto 1\otimes 1+X\otimes X\end{cases},\;\theta\colon A\to A,\; \begin{cases}1\phantom{.}\mapsto 0,\\ X\mapsto 0\end{cases}
\]
for the comultiplication and $\cdot\theta$. Furthermore, the very important map $\Phi^-_+$ given by
\[
\Phi^-_+\colon A\to A,\;\begin{cases}1\phantom{.}\mapsto 1,\\ X\mapsto -X.\end{cases}
\]
Khovanov-Lee's variant has a remarkable property in the classical case, i.e. Lee showed that her variant just ``counts'' the number of components of the c-link, i.e. she showed that (for $R=\bQ$) the homology of a $n$-component link $L$ is
\[
H(\mathcal F_{\mathrm{Lee}}(L))\cong\bigoplus_{2^n}\bQ.
\]
So at the first glance this seems to be a ``boring'' invariant. But Rasmussen used in~\cite{ra} this degeneration in a masterful way to define the \textit{Rasmussen invariant} of a c-knot and his invariant has lots of nice properties.

Therefore, a natural question is if this degeneration of Khovanov-Lee's variant is still true for v-links. In this section we show that this is indeed the case. It is worth noting that this is an unexpected result, since $\theta=0$ for $2^{-1}\in R$ (see the relations in Definition~\ref{defn-category2}). Hence, there are ``tons'' of $0$-morphisms in the complex. But these $0$-morphisms also come with isomorphisms ``in a lot of'' cases.
\vspace*{0.5cm}

The following example for the Khovanov-Lee complex of a v-knot is a blueprint of this effect. It is very important, as indicated in Example~\ref{ex-leedeg} below, that our construction keeps track of the \textit{extra information} how the cobordisms are glued together depending on the orientations of the v-circle diagrams in the resolutions. We note that, even though the orientations can be read off locally, this information has some ``global character''.
\vspace*{0.5cm}
\begin{ex}\label{ex-leedeg}
Consider the diagram of the virtual trefoil $L_D$ given in Fig.~\ref{figure-big2}. In this example the number of negative crossings is zero, i.e. the leftmost object is the $0$-degree part (we mean homological degree).
\begin{figure}[ht]
      \centerline{\includegraphics[width=0.575\linewidth]{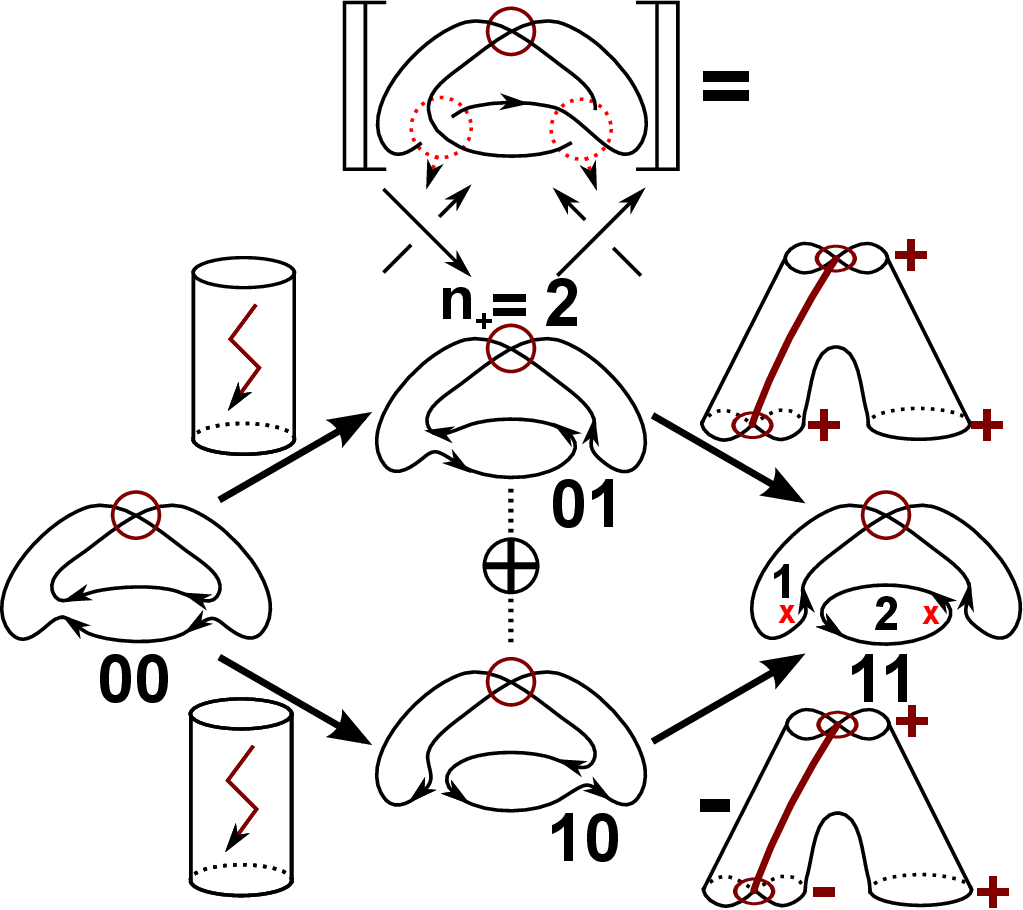}}
  \caption{The Khovanov-Lee complex of the v-trefoil. We note that the first map is a $0$-morphism, but the second is an isomorphism.}
  \label{figure-big2}
\end{figure}
Let us consider $R=\bQ$. Then $\theta=0$ and therefore the first two maps are $0$-morphisms. But note that the two right morphisms are not the same, i.e. one is $\Delta^+_{++}$ and the other is $\Delta^+_{-+}$. So on the algebraic level we get, using the maps from before, the following complex, if we fix $B_1=\{1,X\}$ as a basis for $A$ and $B_2=\{1\otimes 1,1\otimes X,X\otimes 1,X\otimes X\}$ for $A\otimes A$.
\[
\begin{xy}
  \xymatrix{
      A\ar[rrrr]|{\begin{pmatrix}
      0 & 0 \\
      0 & 0 \\
      0 & 0 \\
      0 & 0 
      \end{pmatrix}} &  & & & A\oplus A\ar[rrrr]|{\begin{pmatrix}
      0 & 1 & 0 & -1\\
      1 & 0 & 1 & 0\\
      1 & 0 & -1 & 0\\
      0 & 1 & 0 & 1
      \end{pmatrix}} & & & & A\otimes A.
      }
\end{xy}
\]
An easy calculations shows that the second matrix is an isomorphism. Hence, the homology of the virtual trefoil is only non-trivial for $k=0$, i.e.
\[
H_k(\mathcal F_{\mathrm{Lee}}(L_D))=\begin{cases}\bQ\oplus\bQ, & \text{if}\;k=0,\\ 0, &\text{else}.\end{cases}
\]
Another example is the v-knot in Fig.~\ref{figure-rasexample}, e.g. with the pictured orientation and numbering of the circles from left to right, the three outgoing morphisms from resolution $000$ to $001$, $010$ and $100$ are (up to, in this case, not important signs) the morphisms $m^{+-}_-$, $m^{++}_+$ and $m^{--}_-$, i.e. one alternating and two non-alternating. Hence, the kernel is trivial. The reader should check that the rest also works out in the same fashion as before.
\end{ex}
\subsection{Non-alternating resolutions}
We prove the following interesting result about the number of decorations of v-link resolutions with the ``colours'' down and up. Note that we call an oriented resolution $\mathrm{Re}$ of a v-link diagram \textit{non-alternating} if it is of the form $\uu$ or $\dd$ at the corresponding positions of the saddles. Recall that all the v-link diagrams should be oriented and that such a diagram with $n\in\bN_{>0}$ components has $2^n$ different orientations $\mathrm{Or}_1,\dots,\mathrm{Or}_{2^n}$.
\vspace*{0.5cm}

We note that one can also colour the resolutions with ``honest'' colours, say red and green, in such a way that the colour changes at every v-crossing. We call this a \textit{colouring of a v-link resolution} if at the corresponding saddle-position the colours are different, i.e. (red,green) or (green,red). The reader should compare this with the coloured dots in Fig.~\ref{figure0-main}.
\begin{thm}\label{thm-nonalternating}(\textbf{Non-alternating resolutions}) Let $L_D$ denote a v-link diagram with $n\in\bN_{>0}$ components. There are bijections of sets
\begin{align*}
\{\mathrm{Or}\mid \mathrm{Or}\;\text{is an orientation of}\;L_D\}&\simeq\{\mathrm{Re}\mid \mathrm{Re}\;\text{is a non-alternating resolution of}\;L_D\}\\
&\simeq\{\mathrm{Co}\mid \mathrm{Co}\,\text{is a coloured resolution of}\;L_D\}.
\end{align*}
If $L_D$ is a v-knot diagram, i.e. $n=1$, then the two non-alternating resolutions are in homology degree $0$. A similar statement holds for the coloured resolutions.
\end{thm}
\begin{proof}
With a slight abuse of notation let us denote the first two sets by $\mathrm{Or}$ and $\mathrm{Re}$. To show the existence of a bijection we construct an explicit map $f\colon\mathrm{Or}\to\mathrm{Re}$ and its inverse.

Given an orientation $\mathrm{Or}$ of the v-link diagram $L_D$, the map $f$ should assign the resolution $\mathrm{Re}$ which is obtained by replacing every oriented crossing of the form $\overcrossing$ and $\undercrossing$ with $\uu$ (and the same for rotations). This is clearly an injection.
\vspace*{0.5cm}

Now, given a non-alternating resolution $\mathrm{Re}$, we assign to it an orientation of $L_D$ in the following way. At any non-alternating part of the form $\uu$ and $\dd$ replace the non-alternating part with the corresponding oriented crossing $\overcrossing$ and $\undercrossing$ (or a rotation in the $\dd$ case).

Note that both maps are well-defined and that these two maps are clearly inverses for a v-knot diagram. Moreover, the corresponding non-alternating resolutions are in homology degree $0$, since all $n_+$-crossings are resolved $0$ and all $n_-$-crossings are resolved $1$ in this procedure.   
\vspace*{0.5cm}

To see the second bijection use a checker-board colouring of the v-link diagram. Then start at any point of the non-alternating resolution and use the right-hand rule, i.e. the index finger follows the orientation and the string should get the colour of the face on the side of the thumb. As above, one checks that all $n_+$-crossings are resolved $0$ and all $n_-$-crossings are resolved $1$.
\end{proof}
\begin{cor}\label{cor-numbernonalt}
Let $L_D$ be a v-link diagram with $n$ components. Then it has $2^n$ non-alternating resolutions.
\end{cor}
\begin{proof}
Such a diagram has $2^n$ possible orientations. Then the bijection of Theorem~\ref{thm-nonalternating} finishes the proof.
\end{proof}
\begin{ex}\label{ex-nonalt}
Let $L_D$ be the v-knot diagram in Fig.~\ref{figure-rasexample}.
\begin{figure}[ht]
      \centerline{\includegraphics[width=0.55\linewidth]{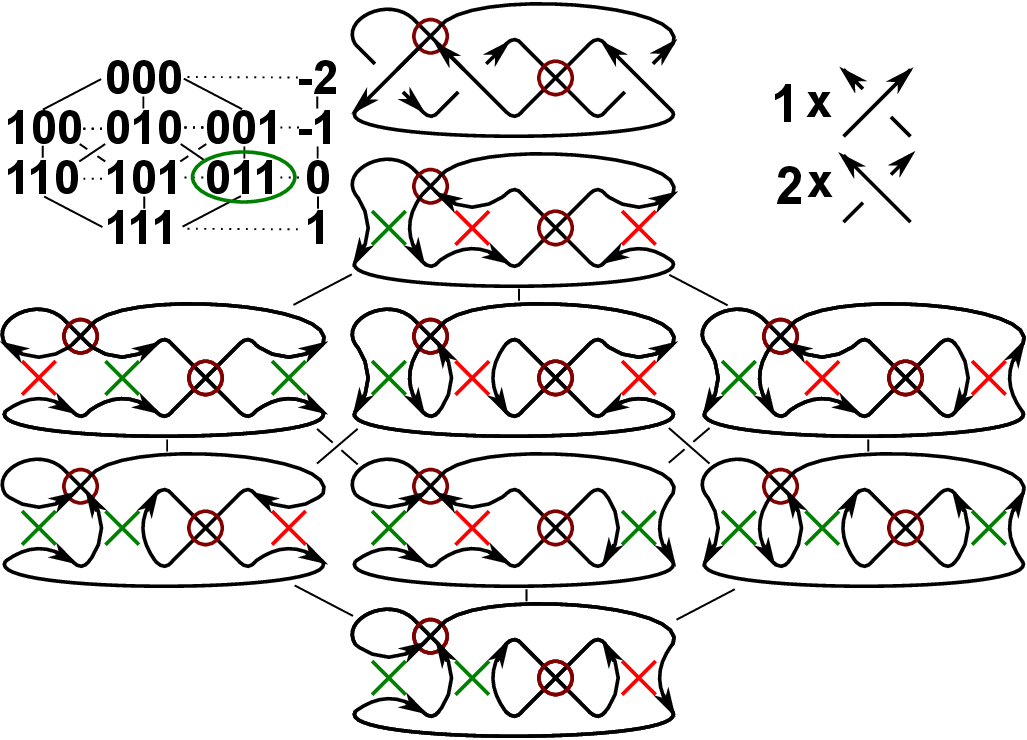}}
  \caption{There are exactly two non-alternating resolutions, i.e. the one pictured and the same resolution, but with all orientations reversed.}
  \label{figure-rasexample}
\end{figure}

Then only the $011$ resolution of the v-knot diagram allows a non-alternating resolution. Moreover, the orientation of the diagram induces this non-alternating resolution by replacing the three crossings with $\dd,\uu$ and $\uu$. The other orientation induces the non-alternating resolution $\uu,\dd$ and $\dd$. Note that, by construction, these resolutions are in homology degree $0$. A computation as in Example~\ref{ex-leedeg} shows that these two non-alternating resolutions give the only two generators of the homology, i.e.
\[
H_k(\mathcal F_{\mathrm{Lee}}(L_D))=\begin{cases}\bQ\oplus\bQ, & \text{if}\;k=0,\\ 0, &\text{else}.\end{cases}
\]
\end{ex}
\subsection{Degeneration}
We start by recalling the motivation, definition and some basic properties of the \textit{Karoubi envelope} of a pre-additive category $\mathcal C$. We denote the envelope as before by $\KAR(\mathcal C)$.
\vspace*{0.5cm}

For any category the notion of an idempotent morphism, i.e. a morphism with $e\circ e=e$, makes sense. Moreover, in a pre-additive category the notion $\mathrm{id}-e$ also makes sense. A classical trick in modern algebra is to use an idempotent, e.g. in $\mathrm{End}_K(V)$ for a given $K$-vector space $V$ (where $K$ is any field), to split the algebra into
\[
\mathrm{End}_K(V)\cong \mathrm{im}(e)\oplus\mathrm{im}(\mathrm{id}-e).
\]
Hence, it is a natural question to ask if on can ``split'', given an idempotent $e$, an object of a category $\mathcal O$ in the same way, i.e.
\[
\mathcal O\cong \mathrm{im}(e)\oplus \mathrm{im}(\mathrm{id}-e).
\]
The main problem is that the notion of an ``image'' of a morphism could possibly not exist in an arbitrary category. The Karoubi envelope is an extension of a category such that for a given idempotent $e$ the notion $\mathrm{im}(e)$ makes sense. Therefore, one can ``split'' a given object in the Karoubi envelope that could be indecomposable in the category itself.
\vspace*{0.5cm}
\begin{defn}\label{defn-karoubi}
Let $\mathcal C$ be a category and let $e,e^{\prime}\colon\mathcal O\to\mathcal O$ denote idempotents in $\mathrm{Mor}(\mathcal C)$. The \textit{Karoubi envelope of $\mathcal C$}, denoted by $\KAR(\mathcal C)$, is the following category.
\begin{itemize}
\item Objects are ordered pairs $(\mathcal O,e)$ of an object $\mathcal O$ and an idempotent $e$ of $\mathcal C$.
\item Morphisms $f\colon(\mathcal O,e)\to(\mathcal O^{\prime},e^{\prime})$ are all arrows $f\colon\mathcal O\to\mathcal O^{\prime}$ of $\mathcal C$ such that the equation $f=f\circ e=e^{\prime}\circ f$ holds.
\item Compositions are defined in the obvious (in this case it is really ``obvious'') way. The identity of an object $(\mathcal O,e)$ is $e$ itself.
\end{itemize}
It is straightforward to check that this is indeed a category. We denote an object $(\mathcal O,e)$ by $\mathrm{im}(e)$, the \textit{image} of the idempotent $e$. Moreover, we identify the objects of $\mathcal C$ with their image via the embedding functor
\[
\iota\colon\mathcal C\to\KAR(\mathcal C),\;\mathcal O\mapsto (\mathcal O,\mathrm{id}).
\] 
\end{defn}
Note that, if $\mathcal C$ is pre-additive, then $\mathrm{id}-e$ is also an idempotent and, under the identification above, we can finally write
\[
\mathcal O\cong \mathrm{im}(e)\oplus \mathrm{im}(\mathrm{id}-e).
\]
The following proposition is ``well-known''\footnote{``Well-known'' means ``strictly more than one person knows its true''.}. The proposition allows us to shift the problem if two chain complexes are homotopy equivalent to the Karoubi envelope. Recall that $\Kom(\mathcal C)$ denotes the category of formal chain complexes.
\vspace*{0.5cm}
\begin{prop}\label{prop-karoubi}
Let $(C,c),(D,d)$ be two objects, i.e. formal chain complexes, of $\Kom(\mathcal C)$. If the two objects are homotopy equivalent in $\Kom(\KAR(\mathcal C))$, then the two objects are also homotopy equivalent in $\Kom(\mathcal C)$.
\end{prop}
\begin{proof}
See e.g. Proposition 3.3 in~\cite{bnsm}.
\end{proof}
We define the two orthogonal idempotents $\mathrm{d},\mathrm{u}$ now and show some basic, but very important, properties afterwards.
\vspace*{0.5cm}

We call the idempotents \textit{``down and up''}. The reader should be careful not to confuse them with the orientations on the resolutions or the colourings of Theorem~\ref{thm-nonalternating}, i.e. latter colours change at v-crossings, but ``down and up'' do not change.
\vspace*{0.5cm}

\begin{defn}\label{defn-idempotent}
We call the two cobordisms in Fig.~\ref{figure-idem} the \textit{``down and up'' idempotents}. We denote them by $\mathrm{d}$ and $\mathrm{u}$.
\begin{figure}[ht]
      \centerline{\includegraphics[width=0.6\linewidth]{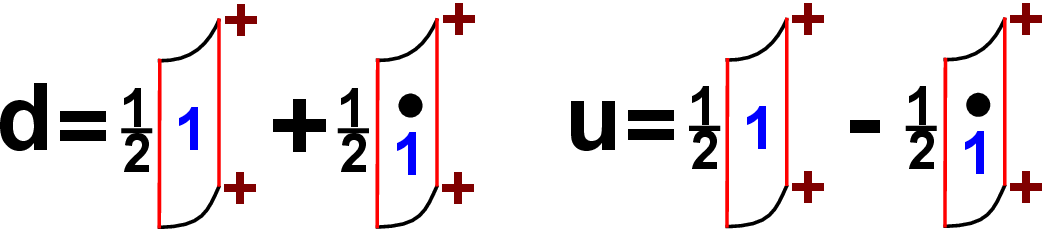}}
  \caption{The two idempotents up and down.}
  \label{figure-idem}
\end{figure}

Recall that the dot represents $\frac{1}{2}$-times a handle.
\end{defn}
It is worth noting that (e) is very important. Moreover, we write $\Phi^-_+$ instead of $\Phi^-_+(1)$.
\begin{lem}\label{lem-idem}
The cobordisms $\mathrm{d},\mathrm{u}$ satisfy the following identities.
\begin{enumerate}
\item[(a)] $\mathrm{d}^2=\mathrm{d}$ and $\mathrm{u}^2=\mathrm{u}$ (idempotent).
\item[(b)] $\mathrm{d}\circ\mathrm{u}=0=\mathrm{u}\circ\mathrm{d}$ (orthogonal).
\item[(c)] $\mathrm{d}+\mathrm{u}=\mathrm{id}$ (complete).
\item[(d)] $\mathrm{id}_{\mathrm{dot}}\circ\mathrm{d}=\mathrm{d}$ and $\mathrm{id}_{\mathrm{dot}}\circ\mathrm{u}=-\mathrm{u}$ (Eigenvalues).
\item[(e)] $\Phi^-_+\circ\mathrm{d}=\mathrm{u}\circ\Phi^-_+$ and $\mathrm{d}\circ\Phi^-_+=\Phi^-_+\circ\mathrm{u}$ (change of orientations).
\item[(f)] $[\mathrm{d},\Phi^-_+]=\mathrm{id}(1)_{\mathrm{dot}}=-[\mathrm{u},\Phi^-_+]$ (commutator relation).
\end{enumerate}
\end{lem}
\begin{proof}
All equations are straightforward to prove. One has to use the dot-relations from Fig.~\ref{figure-dotrel} and the relations from Definition~\ref{defn-category}.

In (d)+(f) the surface $\mathrm{id}(1)_{\mathrm{dot}}$ denotes an identity with an extra dot and $+1$ as an indicator.

Beware that the dot represents $\frac{1}{2}$-times a handle. This forces a sign change after composition with the cobordism $\Phi^-_+$. The reader should compare this with the relations in Definition~\ref{defn-category3}.
\end{proof}
Now we take a look at the Karoubi envelope. The discussion above shows that there is an isomorphism
\[
\none\simeq\down\oplus\up.
\]
With this notation we get
\[
\smoothinga\simeq\ddp\oplus\dup\oplus\udp\oplus\uup\;\;\text{ and }\;\;\hsmoothinga\simeq\ddpp\oplus\dupp\oplus\udpp\oplus\uupp.
\]
Recall that the standard orientation for the complex $\bn{\overcrossing}$ is (see e.g. Fig.~\ref{figure0-main})
\[
\bn{\overcrossing}=\du\xrightarrow{S(1)^{++}_{++}}\ler.
\]
In order to avoid mixing the notions of the down and up-colours and the orientations we denote this complex simply by $\bn{\overcrossing}^{++}_{++}$, i.e. standard orientations for all strings. Moreover, under the convention left=first superscript, right=second superscript, bottom=first subscript and top=second subscript, a notation like $\bn{\overcrossing}^{+-}_{-+}$ makes sense, i.e. act by $\Phi^-_+$ at the corresponding positions. The following theorem is a main observation of this section. It is worth noting again that (e) of Lemma~\ref{lem-idem} is crucial for the theorem.
\begin{thm}\label{thm-karoubi}
In $\ukobk_R$, there are sixteen chain homotopies (only four are illustrated, but it should be clear (hopefully) how the rest works)
\[
\begin{xy}
  \xymatrix{
      \bn{\overcrossing}^{++}_{++}\simeq_h\dup\oplus\udp\xrightarrow{0}\dupp\oplus\udpp, & \bn{\overcrossing}^{+-}_{-+}\simeq_h\ddp\oplus\uup\xrightarrow{0}\ddpp\oplus\uupp,\\
      \bn{\overcrossing}^{++}_{+-}\simeq_h\dup\oplus\udp\xrightarrow{0}\ddpp\oplus\uupp, & \bn{\overcrossing}^{+-}_{++}\simeq_h\ddp\oplus\uup\xrightarrow{0}\dupp\oplus\udpp.
      }
\end{xy}
\]
Moreover, similar formulas hold for $\bn{\undercrossing}$. 
\end{thm}
\begin{proof}
We use the observations from above, i.e. in the Karoubi envelope the differential of $\bn{\overcrossing}^{++}_{++}$ is a $4\times 4$-matrix of saddles. Hence, for $\bn{\overcrossing}^{++}_{++}$ we get (for simplicity write $S=S(1)^{++}_{++}$ and $S_{\mathrm{d}}$ and $S_{\mathrm{u}}$ for the saddle under the action of down and up)
\[
\begin{xy}
  \xymatrix{
       \text{\phantom{.}}\ddp\oplus\dup\oplus\udp\oplus\uup\ar[rrrr]|{\begin{pmatrix}
      S_{\mathrm{d}} & 0 & 0 & 0\\
      0 & 0 & 0 & 0\\
      0 & 0 & 0 & 0\\
      0 & 0 & 0 & S_{\mathrm{u}}
      \end{pmatrix}} & &  & &  \text{\phantom{.}}\ddpp\oplus\dupp\oplus\udpp\oplus\uupp.
      }
\end{xy}
\]
This is true, because all other saddles are killed by the orthogonality relations of the colours down and up, i.e. (b) of Lemma~\ref{lem-idem}.

Note that both non-zero saddles are invertible, i.e. their inverses are the saddles
\[
\frac{1}{2}(S\colon \hsmoothing \to \smoothing)_{\mathrm{d}}\;\;\text{and}\;\;-\frac{1}{2}(S\colon \hsmoothing \to \smoothing)_{\mathrm{u}}
\]
with only $+$ as boundary decorations. To see this one uses Lemma~\ref{lem-idem} and the neck cutting relation. Thus, we get
\[
\bn{\overcrossing}^{++}_{++}\simeq_h\dup\oplus\udp\xrightarrow{0}\dupp\oplus\udpp.
\]
To prove the rest of the statements one has to use the relation (e) of Lemma~\ref{lem-idem}, i.e. the only surviving objects change according to the action of $\Phi^-_+$. We note again that this is a very important observation, i.e. with a different action of $\Phi^-_+$ this would not be true anymore.

For $\bn{\undercrossing}^{++}_{++}$ one can simply copy the arguments from before.
\end{proof}
The following corollary is an application of the semi-local properties of our construction, i.e. we use Theorem~\ref{thm-pcomplex} and Corollary~\ref{cor-pcomplex} to avoid the usage of signs and Lemma~\ref{lem-pcomplex} to see that the involved choices do not matter up to chain isomorphisms.
\begin{cor}\label{cor-importantsaddles}
Let $T^k_D$ be a v-tangle diagram with $m$ crossings. Then $\bn{T^k_D}$ is chain homotopic to a chain complex $(C_*,c_*)$ with only $0$-differentials and objects coloured by the orientations of the resolutions of $T^k_D$, i.e. if a resolution of $\bn{T^k_D}$ is locally of the form
\[
\du\;\;\text{ or }\;\;\dd\;\;\text{ or }\;\;\uu\;\;\text{ or }\;\;\ud
\]
then $(C_*,c_*)$ is locally of the form
\[
\dup\oplus\udp\;\;\text{ or }\;\;\ddp\oplus\uup\;\;\text{ or }\;\;\ddp\oplus\uup\;\;\text{ or }\;\;\dup\oplus\udp.
\]
\end{cor}
\begin{proof}
We note that we work in the Karoubi envelope, but with Proposition~\ref{prop-karoubi} we see that we are free to do so. Moreover, as stated above, we do not care about signs or choices at this point.

Then the statement follows from Theorem~\ref{thm-karoubi} together with the Theorem~\ref{thm-pcomplex} in Sec.~\ref{sec-vkhca}. To be more precise, we copy the arguments from Theorem~\ref{thm-karoubi} for the saddles of the complex $\bn{T^k_D}$ with a $+1,-1$-indicator. Note that these saddles have an extra action of $\Phi^-_+$ at some of its boundary components. That is why the parts of $(C_*,c_*)$ are locally as illustrated above.

Moreover, the saddles with a $0$-indicator are $0$-morphisms for $\frac{1}{2}\in R$ and their local decomposition is the one given above, since they will, by construction, always be between non-alternating parts of the resolutions and due to the orthogonality relations for up and down, namely part (b) of Lemma~\ref{lem-idem}, the $\dup$ and $\udp$ parts will be therefore killed (the only possibility how they close is as the rightmost case of Fig.~\ref{figure0-order}).
\end{proof}
As an application of the Theorems~\ref{thm-nonalternating} and~\ref{thm-karoubi} above, we get the desired statement for v-link diagrams. That is, we have the following.
\begin{thm}\label{thm-leedeg}(\textbf{Degeneration}) Let $L_D$ denote a $n$-component v-link diagram. Then $\bn{L_D}_{\mathrm{Lee}}$ is homotopy equivalent (in $\ukob_R$) to a chain complex with only zero differentials and $2^n$ generators given by the $2^n$ non-alternating resolutions.

If $n=1$, i.e. $L_D$ is a v-knot diagram, then the two generators are in homology degree $0$.
\end{thm}
\begin{proof}
We will suppress the notion of the x-markers and the formal signs of the morphisms to maintain readability. Moreover, we will choose a specific orientation for the resolutions. We can do both freely because of Lemma~\ref{lem-commutativeindependence}.

The main part of the proof will be to choose the orientations in a ``good'' way and use Corollary~\ref{cor-importantsaddles}. Moreover, with Theorem~\ref{thm-karoubi}, we see that the complex will be homotopy equivalent to a complex with only $0$-differentials. Hence, the only remaining thing is to show that the number of generators will work out as claimed.

Note that, if a resolution contains a lower part of a multiplication or a upper part of a comultiplication, then by Corollary~\ref{cor-importantsaddles}, this resolution is killed, because these will always be alternating, e.g. $\du$, but will connect as the $\pm 1$ cases of Fig.~\ref{figure0-order} (the strings are closed with an even number of v-crossings). Moreover, we can ignore top and bottom parts of $\theta$, since they will always be non-alternating.

Now we define the \textit{dual graph of a resolution}, denoted $\mathcal D$, as follows. Recall that a resolution is a four valent graph without any c-crossings. Any edge of this graph is a vertex of $\mathcal D$. Two vertices are connected with a labelled edge iff they are connected by a v-crossing $\virtual$ or a $\smoothing$ (or rotations) that is a top part of a multiplication or a bottom part of a comultiplication. First edges should be labelled $v$, the second type of edges should get a labelling that corresponds to the given orientation of the resolution, that is an ``a'' for alternating orientations and a ``n'' otherwise. We will work with the simple graph of that type, i.e. remove circles or parallel edges of the same type. See Fig.~\ref{figure-dual}, i.e. the figure shows two resolutions from Fig.~\ref{figure-rasexample} and their dual graphs. Note that the leftmost $\smoothing$ of the 011 resolution is part of a $\theta$. 
\begin{figure}[ht]
     \centerline{\includegraphics[width=0.7\linewidth]{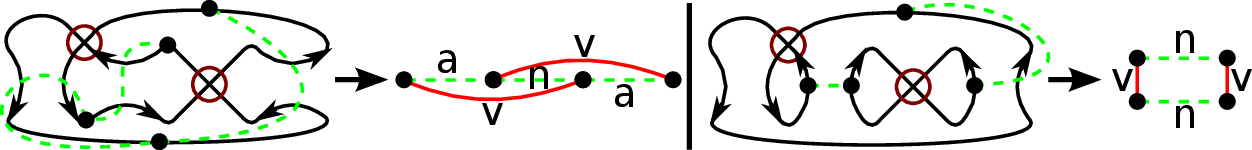}}
  \caption{The resolutions 000 and 011 and their dual graphs.}
  \label{figure-dual}
\end{figure}

The advantage of this notation is that the question of surviving resolutions simplifies to the question of a colouring of the dual graph, i.e. a colouring of the dual graph is a colouring with two colours, say red and green, such that any $v$-labelled or $a$-labelled edge has two equally coloured adjacent vertices, but any $n$-labelled edge has two equal colours at adjacent vertices. The reader should compare this to Theorem~\ref{thm-karoubi} and Corollary~\ref{cor-importantsaddles}.

Then, because of Corollary~\ref{cor-importantsaddles}, a resolution will have surviving generators iff it does not contain lower parts of multiplications or upper parts of comultiplications and, given an orientation of the resolution, it allows a colouring which has the properties described above. For example, the left resolution in Fig.~\ref{figure-dual} does not allow such a colouring, but the right one does.

Recall that the number of crossings is finite. Hence, we can choose an orientation of any resolution such that the number $m$ of alternating crossings is minimal. The rest is just a case-by-case check, i.e. we have the following three cases.
\begin{itemize}
\item[(i)] The dual graph of the resolution is a tree, i.e. no circles.
\item[(ii)] All circles in the dual graph have an even number of $v$-labelled edges.
\item[(iii)] There is one circle in the dual graph with an odd number of $v$-labelled edges.
\end{itemize}
If $m=0$, i.e. the resolution is non-alternating, we get exactly the claimed number of generators, since there are, by construction, no lower parts of multiplication or upper parts of comultiplications and the dual graph is of type (i) or (ii) and in both cases the graph can be coloured.

So let $m>0$ and let $c$ be an alternating crossing in a resolution $\mathrm{Re}$. The whole resolution is killed if the $c$ is a lower part of a multiplication or an upper part of a comultiplication. Hence, we can assume that all alternating crossings of $\mathrm{Re}$ are either top components of multiplications or bottom components of comultiplications.

So we only have to check the three cases from above. If the resolution is one of type (i), then it is possible to choose the orientations in such a way that all crossings are non-alternating, i.e. this would be a contradiction to the minimality of $m$.

If the resolution is of type (ii), then the resolution only survives, i.e. the dual graph allows a colouring, iff the number of other alternating crossings in every circle is even. But in this case one can also choose an orientation with a lower number of non-alternating crossings. Hence, we would get a contradiction to the minimality of $m$ again. An analogous argument works in the case of type (iii), i.e. contradicting the minimality of $m$ again\footnote{It is worth noting that these arguments work because of the well-known (at least to someone) fact that a graph allows a $2$-colouring iff it has no circles of odd length.}.   

Hence, only non-alternating resolutions generate non-vanishing objects and any non-alternating resolution will create exactly two of these. Thus, with Theorem~\ref{thm-nonalternating}, the statement follows.
\end{proof}
\begin{ex}\label{ex-nonalt2}
As an example how the Theorem~\ref{thm-leedeg} works consider the v-knot diagram of Example~\ref{ex-nonalt} again.

The theorem tells us that the resolution 000 should not contribute to the number of generators, i.e. it should get killed. To see this, we first note that in the Karoubi envelope there are $4^3$ different direct summands of coloured (with the idempotents down $d$ and up $u$) versions of the resolution, i.e. four for each crossing. But most of them are killed by the orthogonality of $d$ and $u$, i.e. the two components of the resolution need to have the same colour. Hence, we have the four remaining summands as shown in Fig.~\ref{figure-ras2}.
\begin{figure}[ht]
  	 \centerline{\xy
  	 (0,0)*{\includegraphics[scale=0.7]{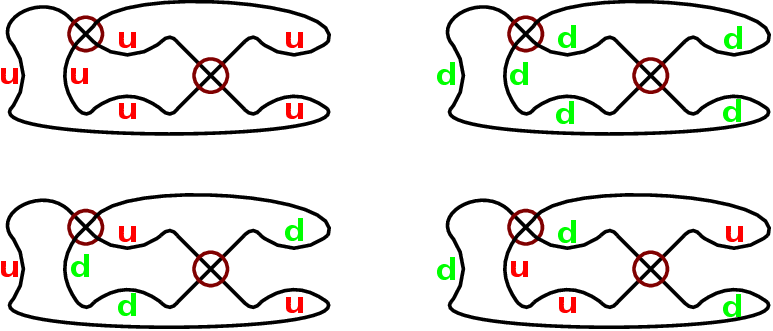}};
     (0,10.5)*{m_{*00}};
     (0,-9.5)*{m_{0*0}};
     (0,-15.5)*{m_{00*}};
     \endxy}
  \caption{The remaining four coloured versions of resolution 000.}
  \label{figure-ras2}
\end{figure}

Let us denote the three multiplications with this resolution $\gamma_{000}$ as source by
\[
m_{*00}\colon\gamma_{000}\to\gamma_{100}\;\;\text{and}\;\;m_{0*0}\colon\gamma_{000}\to\gamma_{010}\;\;\text{and}\;\;m_{00*}\colon\gamma_{000}\to\gamma_{001}.
\]
If we choose the orientation for $\gamma_{000}$ as indicated in the Fig.~\ref{figure-rasexample}, we see that
\[
m_{*00}\colon\dd\to\ril\;\;\text{and}\;\;m_{0*0}\colon\du\to\ler\;\;\text{and}\;\;m_{00*}\colon\ud\to\ril.
\]
We can now use Corollary~\ref{cor-importantsaddles} to see that the only remaining parts for the three multiplications are as follows.
\[
m_{*00}\colon\ddp\oplus\uup\to\dupp\oplus\udpp,\hspace*{0.5cm}m_{0*0},m_{*00}\colon\dup\oplus\udp\to\dupp\oplus\udpp.
\]
Hence, they pick two distinct coloured versions as illustrated in Fig.~\ref{figure-ras2}. Therefore, there are no surviving generators for the 000 resolution. It should be noted that changing for example the orientation of the leftmost v-circle in Fig.~\ref{ex-nonalt} does not affect the result, since Lemma~\ref{lem-commutativeindependence} ensures that the resulting complexes are isomorphic. And in fact such a change leads to
\[
m_{*00}\colon\ud\to\ril\;\;\text{and}\;\;m_{0*0}\colon\uu\to\ler\;\;\text{and}\;\;m_{00*}\colon\uu\to\ril.
\]
Hence, the $m_{*00}$ and the two multiplications $m_{0*0},m_{*00}$ still pick out different coloured versions of the resolution 000. Therefore, there will not be any surviving generators for this case either.
\end{ex}
We finish by using the functor $\mathcal F_{\mathrm{Lee}}$ to get the corresponding statement in the category $\RMOD$. The reader may compare this to the results in the classical case, e.g. see Proposition 2.4 in~\cite{mtv}.
\begin{prop}\label{prop-algebraicdegen}
Let $L_D$ denote a $n$-component v-link diagram. Then we have the following.
\begin{itemize}
\item[(a)] If $R=\bZ$, then there is an isomorphism
\[
H(\mathcal F_{\mathrm{Lee}}(L_D),\bZ)\cong\bigoplus_{2^n}\bZ\oplus \mathrm{Tor},
\]
where $\mathrm{Tor}$ is all torsion. Moreover, the only possible torsion is $2$-torsion.
\item[(b)] If $R=\bQ$ or $R=\bZ\left[\frac{1}{2}\right]$, then there is an isomorphism
\[
H(\mathcal F_{\mathrm{Lee}}(L_D),R)\cong\bigoplus_{2^n}R.
\]
\end{itemize}
\end{prop}
\begin{proof}
The statement (b) follows from Theorem~\ref{thm-leedeg} above. Recall that the whole construction requires that two is invertible.

For (a) recall the universal coefficients theorem. i.e. there is a short exact sequence
\[
0\rightarrow H_*(\mathcal F_{\mathrm{Lee}}(L_D),\bZ)\otimes_{\bZ} R\rightarrow H_*(\mathcal F_{\mathrm{Lee}}(L_D),R)\rightarrow \mathrm{Tor}(H_{*+1}(\mathcal F_{\mathrm{Lee}}(L_D),\bZ)),R)\rightarrow 0.
\]
Therefore, (a) follows from (b) with $R=\bQ$, since the Tor-functor will vanish in this case and from (b) with $R=\bZ\left[\frac{1}{2}\right]$. Hence, this shows the proposition.
\end{proof}
\section{Computer talk}\label{sec-vkhpro}
In this section we show some basic calculations with a computer program we have written. The program is a MATHEMATICA (see~\cite{wolf}) package called \textit{vKh.m}. There is also a notebook called \textit{vKh.nb}. Both and some calculation results are available online on the author's homepage\footnote{{\tt \tiny http://xwww.uni-math.gwdg.de/dtubben/vKh.htm}}.

The input data is a v-link diagram in a \textit{circuit notation}, i.e. the classical \textit{planar diagram notation}, but we allow v-crossings. Hence, the input data is a string of labelled $X$, i.e crossings are presented by symbols $X_{ijkl}$ where the numbers are obtained by numbering the edges of the v-link diagram, and the edges around the crossing start counting from the lower incoming edge and proceed counterclockwise. We denote such a diagram by CD[X[i,j,k,l],...,X[m,n,o,p]].

After starting MATHEMATICA and loading our package vKh.m, we type in the unknot from Fig.~\ref{figure0-big}, the classical and virtual trefoil. Our notation follows the notation of Green in his nice \textit{table of virtual knots}\footnote{J.~Green, A Table of Virtual Knots, {\tt \tiny http://www.math.toronto.edu/drorbn/Students/GreenJ/} (2004)}.

{\In Unknot:= CD[X[1,3,2,4], X[2,1,3,4]]; Knot21 := CD[X[1,3,2,4], X[4,2,1,3]];\\ Knot36 := CD[X[1,5,2,4], X[5,3,6,2], X[3,1,4,6]];}\vskip 2mm

Let us denote the elements $1,X\in A=\bZ[X]/X^2=0$ by 1=vp[i] and X=vm[i] and tensors of these elements multiplicatively. Here the module $A$ should belong to the $i$-th v-circle. Moreover, we denote by the word $a$, whose letters are from the alphabet $\{0,1,*\}$ with exactly one $*$-entry, the cobordism starting at the resolution $\gamma_{*=0}$ and going to the resolution $\gamma_{*=1}$. Let us check the different morphisms.

{\In {d2[Unknot, "0*"],d2[Unknot, "*0"], d2[Unknot, "1*"],d2[Unknot, "*1"]}}
\begingroup\Out {{vp[1] -> vm[2] vp[1] - vm[1] vp[2], vm[1] -> vm[1] vm[2]}, {vp[1] -> 0,\\vm[1] -> 0}, {vp[1] -> 0, vm[1] -> 0}, {vp[1] vp[2] -> -vp[1], vm[2] vp[1] ->\\-vm[1], vm[1] vp[2] -> -vm[1], vm[1] vm[2] -> 0}}\endgroup\vskip 2mm

We see that the two orientable morphisms are $\Delta^+_{-+}$ and $-m^{--}_-=-m^{++}_+$. With the command KhBracket[Knot,r] we generate the $r$-th module of the complex (here for simplicity without gradings). Moreover, with d[Knot][KhBracket[Knot,r]] we calculate the image of the $r$-th differential for the whole module. Let us check the output.

{\In {KhBracket[Unknot, 0], KhBracket[Unknot, 1], KhBracket[Unknot, 2]}}
\begingroup\Out {{v[0, 0] vm[1], v[0, 0] vp[1]}, {v[0, 1] vm[1] vm[2], v[0, 1] vm[2] vp[1],\\v[0, 1] vm[1] vp[2], v[0, 1] vp[1] vp[2], v[1, 0] vm[1], v[1, 0] vp[1]},\\{v[1, 1] vm[1], v[1, 1] vp[1]}}\endgroup
{\In {d[Unknot][KhBracket[Unknot, 0]], d[Unknot][KhBracket[Unknot, 1]]}}
\begingroup\Out {{v[0, 1] vm[1] vm[2], v[0, 1] vm[2] vp[1] - v[0, 1] vm[1] vp[2]}, {0,\\-v[1, 1] vm[1], -v[1, 1] vm[1], -v[1, 1] vp[1], 0, 0}}\endgroup\vskip 2mm

It is easy to check that the composition $d_1\circ d_0$ is indeed zero.

{\In d[Unknot][d[Unknot][KhBracket[Unknot, 0]]]}
\begingroup\Out {0, 0}\endgroup\vskip 2mm

Let us check this for the other two knots, too.

{\In d[Knot21][d[Knot21][KhBracket[Knot21, 0]]]}
\begingroup\Out {0, 0, 0, 0}\endgroup
{\In {d[Knot36][d[Knot36][KhBracket[Knot36, 0]]], d[Knot36][d[Knot36]\\ [KhBracket[Knot36, 1]]]}}
\begingroup\Out {{0, 0, 0, 0, 0, 0, 0, 0}, {0, 0, 0, 0, 0, 0, 0, 0, 0, 0, 0, 0}}\endgroup\vskip 2mm

Now let us check for the trefoil how the signs of the morphisms work out.

{\In {sgn[Knot36, "00*"], sgn[Knot36, "0*0"], sgn[Knot36, "*00"]}}
\begingroup\Out {1, -1, 1}\endgroup
{\In {sgn[Knot36, "01*"], sgn[Knot36, "10*"], sgn[Knot36, "0*1"], 
 sgn[Knot36, "1*0"], sgn[Knot36, "*01"], sgn[Knot36, "*10"]}}
\begingroup\Out {1, 1, 1, -1, -1, -1}\endgroup
{\In {sgn[Knot36, "11*"], sgn[Knot36, "1*1"], sgn[Knot36, "*11"]}}
\begingroup\Out {1, 1, 1}\endgroup\vskip 2mm

We observe that all of the six different faces have an odd number of signs. For example the face $F1=(\gamma_{000},\gamma_{001}\oplus\gamma_{010},\gamma_{011})$ gets a sign from the morphism $d_{0*0}$. Furthermore, the face $F2=(\gamma_{100},\gamma_{101}\oplus\gamma_{110},\gamma_{111})$ gets a sign from the morphism $d_{1*0}$.

The first face is of type 2b and the second is of type 1b. Hence, after a virtualisation the latter should have an even number of signs, but the first should have an odd number signs. Let's check this. First we define a new knot diagram which we obtain by performing a virtualisation on the second crossing of the trefoil.

{\In Knot36v := CD[X[1, 4, 2, 5], X[2, 5, 3, 6], X[3, 6, 4, 1]];}
{\In {sgn[Knot36v, "00*"], sgn[Knot36v, "0*0"], sgn[Knot36v, "*00"]}}
\begingroup\Out {1, -1, 1}\endgroup
{\In {sgn[Knot36v, "01*"], sgn[Knot36v, "10*"], sgn[Knot36v, "0*1"], 
 sgn[Knot36v, "1*0"], sgn[Knot36v, "*01"], sgn[Knot36v, "*10"]}}
\begingroup\Out {1, 1, 1, -1, -1, -1}\endgroup
{\In {sgn[Knot36v, "11*"], sgn[Knot36v, "1*1"], sgn[Knot36v, "*11"]}}
\begingroup\Out {1, -1, 1}\endgroup\vskip 2mm

Indeed only the sign of the morphism $d_{1*1}$ is different now. Hence, the face $F1$ still has an odd number, but the face $F2$ has an even number of signs. This should cancel with the extra sign of the pantsdown morphism $d_{1*1}$. And indeed:

{\In {d[Knot36v][d[Knot36v][KhBracket[Knot36v, 0]]], d[Knot36v][d[Knot36v]\\ [KhBracket[Knot36v, 1]]]}}
\begingroup\Out {{0, 0, 0, 0, 0, 0, 0, 0}, {0, 0, 0, 0, 0, 0, 0, 0, 0, 0, 0, 0}}\endgroup\vskip 2mm

Let us look at some calculation results for the four knots. The output is Betti[q,t], i.e. the dimension of the homology group in quantum degree $q$ and homology degree $t$. The unknot should have trivial homology.

{\In vKh[Unknot]}
\begingroup\Out {Betti[-1,-2] = 0, Betti[-1,0] = 0, Betti[0,-2] = 0, Betti[0,-1] = 1,\\Betti[0,0]= 0, Betti[0,1] = 1, Betti[0,1] = 0, Betti[1,0] = 0, Betti[1,2] = 0}\endgroup
\begingroup\Out 1/q + q\endgroup\vskip 2mm

For the other outputs we skip the Betti numbers. One can read them off from the polynomial: Betti[a,b] is the coefficient for $q^at^b$. The trefoil and its virtualisation have the same output (as they should).

{\In vKh[Knot21]}
\begingroup\Out 1/q^3 + 1/q + 1/(q^6 t^2) + 1/(q^2 t)\endgroup
{\In vKh[Knot36]}
\begingroup\Out 1/q^3 + 1/q + 1/(q^9 t^3) + 1/(q^5 t^2)\endgroup
{\In vKh[Knot36v]}
\begingroup\Out 1/q^3 + 1/q + 1/(q^9 t^3) + 1/(q^5 t^2)\endgroup\vskip 2mm

Let us check that the graded Euler characteristic is the Jones polynomial\footnote{To simplify the outputs we have avoided to include the orientation of the v-links in the input, i.e. every output needs a degree shift.}.

{\In Factor[(vKh[Knot21] /. t -> -1)/(q + q^-1)]}
\begingroup\Out (1 - q^2 + q^3)/q^5\endgroup
{\In Factor[(vKh[Knot36] /. t -> -1)/(q + q^-1)]}
\begingroup\Out (-1 + q^2 + q^6)/q^8\endgroup\vskip 2mm

Another observation is the following. The map $\Phi^-_+$ sends $1$ to itself, but $X$ to $-X$. Hence, there is a good change for 2-torsion. Let us check. Here Tor[q,t] denotes the $\bZ/p\bZ$-rank minus the $\bZ$-rank (both graded) of Betti[q,t]$\otimes\bZ/p\bZ$. Even the v-trefoil has 2-torsion, but no 3-torsion.

{\In vKh[Knot21,2]}
\begingroup\Out {Tor[-2,-6] = 0, Tor[-2,-4] = 0, Tor[-2,-2] = 0, Tor[-1,-4] = 1,\\Tor[-1,-2] = 0, Tor[0,-3] = 0, Tor[0,-1] = 0}\endgroup
\begingroup\Out 1/(q^4 t)\endgroup
{\In vKh[Knot21,3]}
\begingroup\Out {Tor[-2,-6] = 0, Tor[-2,-4] = 0, Tor[-2,-2] = 0, Tor[-1,-4] = 0,\\Tor[-1,-2] = 0, Tor[0,-3] = 0, Tor[0,-1] = 0}\endgroup
\begingroup\Out 0\endgroup\vskip 2mm

There seems to be a lot of 2-torsion!

{\In Knot32 := CD[X[2, 6, 3, 1], X[4, 2, 5, 1], X[5, 3, 6, 4]];}
{\In vKh[Knot32]}
\begingroup\Out 1/q^2 + 1/q + q + 1/(q^5 t^2) + 1/(q t) + q^2 t\endgroup
{\In vKh[Knot32,2]}
\begingroup\Out 1/(q^3 t) + t\endgroup\vskip 2mm

Because the virtual Khovanov complex is invariant under virtualisation, there are many examples of non-trivial v-knots with trivial Khovanov complex.

{\In Knot459 := CD[X[2, 8, 3, 1], X[4, 2, 5, 1], X[3, 6, 4, 7], X[5, 8, 6, 7]];}
{\In vKh[Knot32]}
\begingroup\Out 1/q + q\endgroup\vskip 2mm

Let us try a harder example. We mention that the faces are all anti-commutative, and hence, the composition of the differentials is zero. 

{\In Knot53 := CD[X[1, 9, 2, 10], X[2, 10, 3, 1], X[5, 4, 6, 3], X[7, 4, 8, 5],\\X[8, 7, 9, 6]];}
{\In vKh[Knot53]}
\begingroup\Out 2 + 1/q^3 + 1/q^2 + 1/q + 1/(q^7 t^3) + 1/(q^6 t^2) + 1/(q^5 t^2) \\+ 1/(q^3 t^2) + 2/(q^4 t) + 1/(q^2 t) + 1/(q t) + t/q + q^2 t + q^3 t^2\endgroup
{\In vKh[Knot53,2]}
\begingroup\Out 2/q^2 + 1/(q^5 t^2) + 1/(q^4 t) + 1/(q^3 t) + t + q t^2\endgroup
{\In \{d[Knot53][d[Knot53][KhBracket[Knot53, 0]]], d[Knot53][d[Knot53][KhBracket\\ [Knot53, 1]]], d[Knot53][d[Knot53][KhBracket[Knot53, 2]]], d[Knot53][d[Knot53]\\ [KhBracket[Knot53, 3]]]\}}
\begingroup\Out {{0, 0, 0, 0, 0, 0, 0, 0}, {0, 0, 0, 0, 0, 0, 0, 0, 0, 0, 0, 0, 0, 0, 0,\\0, 0, 0, 0, 0, 0, 0, 0, 0}, {0, 0, 0, 0, 0, 0, 0, 0, 0, 0, 0, 0, 0, 0, 0, 0, 0, 0,\\0, 0, 0, 0, 0, 0, 0, 0, 0, 0, 0, 0, 0, 0}, {0, 0, 0, 0, 0, 0, 0, 0, 0, 0, 0, 0, 0,\\0, 0, 0, 0, 0, 0, 0, 0, 0, 0, 0, 0, 0, 0, 0}}\endgroup\vskip 2mm

The virtual Khovanov complex is strictly stronger than the virtual Jones polynomial. The first example appears for v-links with seven crossings. Let's check two examples.

{\In Example1 := CD[X[1, 4, 2, 3], X[2, 10, 3, 11], X[4, 9, 5, 10], X[11, 5, 12, 6],\\X[6, 1, 7, 14], X[12, 8, 13, 7], X[13, 9, 14, 8]]; Example2 := CD[X[1, 4, 2, 3],\\X[2, 11, 3, 10], X[4, 10, 5, 9], X[14, 5, 1, 6], X[6, 12, 7, 11], X[13, 7, 14, 8],\\X[12, 8, 13, 9]]; Example3 := CD[X[1, 4, 2, 3], X[2, 11, 3, 10], X[4, 9, 5, 10],\\X[13, 5, 14, 6], X[6, 11, 7, 12], X[14, 8, 1, 7], X[12, 8, 13, 9]]; Example4 :=\\CD[X[1, 4, 2, 3], X[2, 11, 3, 10], X[4, 10, 5, 9], X[14, 5, 1, 6], X[6, 13, 7, 14],\\X[11, 7, 12, 8], X[12, 8, 13, 9]];}\vskip 2mm

So let us see what our program calculates.

{\In \{vKh[Example1], vKh[Example2], vKh[Example3], vKh[Example4]\}}
\begingroup\Out {2 + 1/q + q + 2 q^2 + 1/(q^3 t^2) + 2/(q^2 t) + q/t + 2 q t + 2 q^4 t\\+ q^3 t^2 + 2 q^5 t^2 + q^7 t^3, 2 + 1/q + q + 2 q^2 + q^3 + 1/(q^3 t^2) + 2/(q^2 t)\\+ q/t + 2 q t + q^2 t + q^3 t + 2 q^4 t + q^2 t^2 + q^3 t^2 + 2 q^5 t^2 + q^6 t^2 + q^6 t^3 + q^7 t^3, 2/q^2 + 1/q + 3 q + 1/(q^6 t^3) + 2/(q^5 t^2) + 1/(q^2 t^2) + 2/(q^3 t)\\+ 2/(q t) + t + 2 q^2 t + q^4 t^2, 1 + 2/q^2 + 2/q + 3 q + 1/(q^6 t^3) + 2/(q^5 t^2)\\+ 1/(q^4 t^2) + 1/(q^2 t^2) + 1/t + 1/(q^4 t) + 2/(q^3 t) + 2/(q t) + t + t/q + 2 q^2 t + q^3 t + q^3 t^2 + q^4 t^2}\endgroup\vskip 2mm

Good news: Example1 and Example2 have the same virtual Jones polynomial ($t=-1$), but different virtual Khovanov homology, i.e. Example2 has the six extra terms (compared to Example1) $q^2t$, $q^2t^2$, $q^3$, $q^3t$, $q^6t^2$ and $q^6t^3$. They all cancel if we substitute $t=-1$. An analogously effect happens for Example3 and Example4. Furthermore, our calculations suggest that this repeats frequently for v-knots with seven or more crossings.

The command line GausstoCD converts \textit{signed Gauss Code} to a CD representation. The signed Gauss code has to start with the first overcrossing. To get the mirror image we can use the rule from below. For example the virtual trefoil and its mirror are not equivalent.
{\In Knot21gauss := "O1-O2-U1-U2-";}
{\In GuasstoCD[Knot21gauss]}
\begingroup\Out CD[X[1, 4, 2, 3], X[2, 1, 3, 4]]\endgroup
{\In GuasstoCD[Knot21gauss] /. X[i_,j_,k_,l_] :> X[i,l,k,j]}
\begingroup\Out CD[X[1, 3, 2, 4], X[2, 4, 3, 1]]\endgroup
{\In \{vKh[GausstoCD[Knot21gauss]], 
 vKh[GausstoCD[Knot21gauss] /. X[i_, j_, k_, l_] :> X[i, l, k, j]]\}}
\begingroup\Out {q + q^3 + q^2 t + q^6 t^2, 1/q^3 + 1/q + 1/(q^6 t^2) + 1/(q^2 t)}\endgroup\vskip 2mm

\begin{figure}[ht]
\centerline{\begin{minipage}[c]{6cm}
	\includegraphics[scale=0.425]{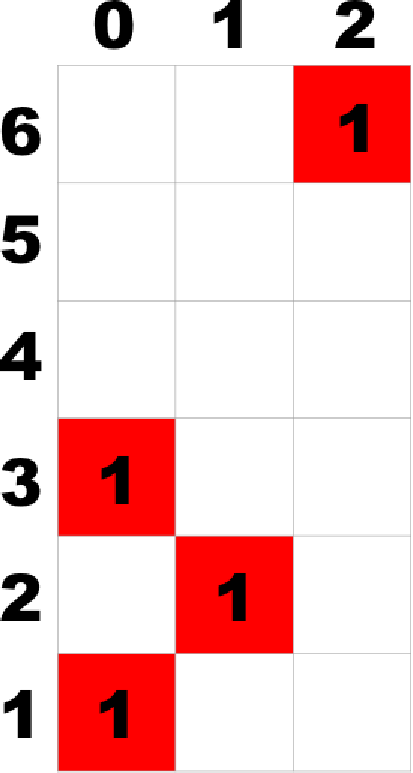}
\end{minipage}
\begin{minipage}[c]{6cm}
	\includegraphics[scale=0.425]{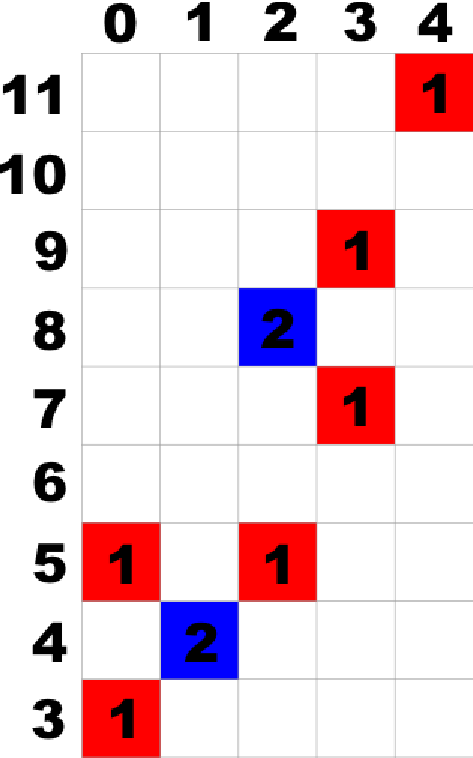}
\end{minipage}}
\caption{Left: Homology of the v-trefoil. Right: Homology of the v-knot 4.1.}
	\label{figure-calc1}
\end{figure}

We used this to calculate the virtual Khovanov homology for all different v-knots with at most five crossings. The input was the list of v-knots from Green's virtual knot table. The results are available on the author's website (as mentioned before). One could visualise the polynomial with the function \textit{Ployplot}. It creates an output as in the Fig.~\ref{figure-calc1}.
{\In Knot41 := "O1-O2-U1-U2-O3-O4-U3-U4-"; vKh[GausstoCD[Knot41]]}
\begingroup\Out q^3 + q^5 + 2 q^4 t + q^5 t^2 + 2 q^8 t^2 + q^7 t^3\\ + q^9 t^3 + 
 q^11 t^4]\endgroup\vskip 2mm
The output of this v-knot and of the mirror of the virtual trefoil is shown in the Fig.~\ref{figure-calc1}. In these pictures the quantum degree is on the y-axis and the homology degree on the x-axis.
\section{Cones, strong deformation retracts and homotopy equivalence}\label{sec-techhomalg}
In this section we have collected some well-known facts from homological algebra about (mapping) cones, strong deformation retracts and homotopy equivalences. We need these facts in this paper. In particular, we need them in the proof of Theorem~\ref{thm-geoinvarianz}. We note that anything in this section can be found in (literally any) book on homological algebra.

In this section let $\mathcal C$ denote any pre-additive category. It should be noted that this includes that the notion ``chain complex'', i.e. $d\circ d=0$, makes sense. We denote the \textit{category of chain complexes of $\mathcal C$} by $\Kom(\mathcal C)$ (in contrast to the category of \textit{bounded} chain complexes $\Kom_b(\mathcal C)$), i.e. the objects are chain complexes with chain groups in $\Ob(\mathcal C)$ and differentials in $\Mor(\mathcal C)$ and the morphisms are chain maps, i.e. sequences of elements of $\Mor(\mathcal C)$ with the standard requirements. 

We denote chain complexes by $C=(C_i,c_i), D=(D_i,d_i)\in\Kom(\mathcal C)$. With a slight abuse of terminology, we call elements of $\Mor(\mathcal C)$ simply ``maps''. All appearing indices should be elements of $\bZ$. Moreover, recall the following three definitions.
\begin{defn}\label{defn-chainhomotopy}
Let $C,D$ be two chain complexes with chain groups $C_i,D_i$ and differentials $c_i,d_i$. Let $\varphi,\varphi^{\prime}\colon C\to D$ be two chain maps. Let $h_i\colon C_i\to D_{i-1}$ be a collection of maps as illustrated below.
\begin{equation*}
\xymatrix@+3em{
{\dots} \ar[r]^{c_{i - 2}}
        & C_{i - 1}
                \ar[r]^{c_{i - 1}}
                \ar@<0.5ex>[d]^{\varphi^{\prime}_{i - 1}}
                \ar@<-0.5ex>[d]_{\varphi_{i - 1}}
                \ar[dl]|*+<1ex,1ex>{\scriptstyle h_{i - 1}}
        & C_i
                \ar[r]^{c_i}
                \ar@<0.5ex>[d]^{\varphi^{\prime}_{i}}
                \ar@<-0.5ex>[d]_{\varphi_{i}}
                \ar[dl]|*+<1ex,1ex>{\scriptstyle h_i}
        & C_{i + 1}
                \ar[r]^{c_{i + 1}}
                \ar@<0.5ex>[d]^{\varphi^{\prime}_{i + 1}}
                \ar@<-0.5ex>[d]_{\varphi_{i + 1}}
                \ar[dl]|*+<1ex,1ex>{\scriptstyle h_{i + 1}}
        & {\dots}
                \ar[dl]|*+<1ex,1ex>{\scriptstyle h_{i + 2}}\\
{\dots} \ar[r]^{d_{i - 2}}
        & D_{i - 1} \ar[r]^{d_{i - 1}}
        & D_i \ar[r]^{d_i}
        & D_{i + 1} \ar[r]^{d_{i + 1}}
        & {\dots}
}
\end{equation*}
The two chain maps $\varphi,\varphi^{\prime}\colon C\to D$ are called \textit{chain homotopic}, denoted by $\varphi\sim_h\varphi^{\prime}$, if
\[
\varphi_i-\varphi^{\prime}_i=h_{i+1}\circ c_i+d_{i-1}\circ h_i\text{ for all }i\in\bZ.
\]
Two chain complexes $C,D$ are called \textit{chain homotopic} if there exists two chain maps $\varphi\colon C\to D$ and $\psi\colon D\to C$ such that
\[
\psi_i\circ\varphi_i\sim_h\mathrm{id}_{C}\;\;\text{ and }\;\;\varphi_i\circ\psi_i\sim_h\mathrm{id}_{D}\text{ for all }i\in\bZ.
\]
Such chain maps $\varphi\colon C\to D$ and $\psi\colon C\to D$ are called \textit{homotopy equivalences}. We use the notation $C\simeq_h D$ to indicated that the chain complexes $C$ and $D$ are homotopic.
\end{defn}
\begin{defn}\label{defn-sdr}Let $\psi\colon D\to C$ be a homotopy equivalence. Assume that $\psi$ has a ``homotopy inverse'' $\varphi\colon C\to D$, that is
\[
\psi_i\circ\varphi_i=\mathrm{id}_{C}\;\;\text{ and }\;\;\varphi_i\circ\psi_i\sim_h\mathrm{id}_{D}\text{ for all }i\in\bZ,
\]
then $\psi$ is called a \textit{deformation retraction}. Moreover, if there exists a homotopy $h$ with $h\circ \varphi=0$, then $\psi$ is called a \textit{strong deformation retraction} and $\varphi$ is called an \textit{inclusion into a strong deformation retract}.
\end{defn}
A (mapping) cone of two chain complexes is defined below. Be careful: Some authors use a different sign convention.
\begin{defn}\label{defn-cone}(\textbf{The cone}) Let $C,D$ be two chain complexes with chain groups $C_i,D_i$ and differentials $c_i,d_i$. Let $\varphi\colon C\to D$ be a chain map. The \textit{cone} of $C,D$ \textit{along $\varphi$} is the chain complex $\Gamma(\varphi\colon C\to D)$ with the chain groups and differentials
\[
\Gamma_i=C_i\oplus D_{i-1}\text{ and } \gamma_i=\begin{pmatrix}-c_i & 0\\ \varphi_{i} & d_{i-1}\end{pmatrix},
\]
i.e. if the two chain complexes $C,D$ look like
\[
C,D\colon\cdots\xrightarrow{c_{i-1},d_{i-1}} C_i,D_i\xrightarrow{c_i,d_i}C_{i+1},D_{i+1}\xrightarrow{c_{i+1},d_{i+1}}C_{i+2},D_{i+2}\xrightarrow{c_{i+2},d_{i+2}}\cdots
\]
then the cone along $\varphi$ is generated by direct sums over the diagonal as shown below.
\begin{equation*}
\xymatrix@+3em{
{\dots} \ar[r]^{c_{i - 2}}
        & C_{i - 1}
                \ar[r]^{-c_{i - 1}}
                \ar[d]|{\varphi_{i - 1}}
                \ar@{.}[dl]|*+<1ex,1ex>{\scriptstyle \oplus}
        & C_i
                \ar[r]^{-c_i}
                \ar[d]|{\varphi_{i}}
                \ar@{.}[dl]|*+<1ex,1ex>{\scriptstyle \oplus}
        & C_{i + 1}
                \ar[r]^{-c_{i + 1}}
                \ar[d]|{\varphi_{i + 1}}
                \ar@{.}[dl]|*+<1ex,1ex>{\scriptstyle \oplus}
        & {\dots}
                \ar@{.}[dl]|*+<1ex,1ex>{\scriptstyle \oplus}\\
{\dots} \ar[r]^{d_{i - 2}}
        & D_{i - 1} \ar[r]^{d_{i - 1}}
        & D_i \ar[r]^{d_i}
        & D_{i + 1} \ar[r]^{d_{i + 1}}
        & {\dots}
}
\end{equation*}
\end{defn}
It is easy to check that the cone gives again a chain complex (here one has to use the signs above). The following well-known proposition concludes this section (it is so well-known that we were unable to find a reference). We need it in order to prove Theorem~\ref{thm-geoinvarianz}.
\begin{prop}\label{prop-sdr}
Let $C,D,E,F$ be four chain complexes and let
\[ 
\begin{xy}
  \xymatrix{
      C  \ar@<2pt>[d]^{f^{\phantom{\prime}}}    &   D\ar@<2pt>[d]^{g^{\prime}}   \\
      E \ar[r]_{\varphi}\ar@<2pt>[u]^{f^{\prime}}             &   F \ar@<2pt>[u]^{g} 
  }
\end{xy}
\]
be a diagram in the category $\Kom(\mathcal C)$. Assume that $f$ is an inclusion into a strong deformation retract $f^{\prime}$ and $g$ is a strong deformation retraction with inclusion $g^{\prime}$. Then
\[
\Gamma(\varphi\circ f)\simeq_h\Gamma(\varphi)\simeq_h\Gamma(g\circ \varphi).
\]
\end{prop}
\begin{proof}
In order to maintain readability, we suppress some subscripts in the following.

To show that $\Gamma(\varphi)\simeq_h\Gamma(\varphi\circ f)$ we denote their differentials by $d^{\Gamma(\varphi\circ f)}$ and $d^{\Gamma(\varphi)}$ respectively. Consider the following diagram
\begin{equation*}
\xymatrix@+3em{
{\Gamma(\varphi\circ f):\;\;\dots} \ar[r]^{}
        & C_{i + 1}\oplus F_i
                \ar[r]^{d^{\Gamma(\varphi\circ f)}}
                \ar@<2pt>[d]_{\psi^{f^{\prime}}}
        & C_{i+2}\oplus F_{i+1}
                \ar[r]^{}
                \ar@<2pt>[d]_{\psi^{f^{\prime}}}
        & {\dots}\\
{\Gamma(\varphi):\;\;\dots} \ar[r]^{}
        & E_{i + 1}\oplus F_i \ar@<2pt>[r]^{d^{\Gamma(\varphi)}}\ar@<2pt>[u]_{\phantom{.}\psi^f}
        & E_{i+2}\oplus F_{i+1} \ar[r]^{}\ar@<2pt>[u]_{\phantom{.}\psi^f}\ar@<2pt>[l]^{\psi^h}
        & {\dots},
}
\end{equation*}
where the three maps $\psi^a,\psi^b$ and $h$ are given by
\[
\psi^f=\begin{pmatrix}f & 0\\ 0& \mathrm{id}\end{pmatrix}\;\text{ and }\;\psi^{f^{\prime}}=\begin{pmatrix}f^{\prime} & 0\\ \varphi\circ h^{\prime}& \mathrm{id}\end{pmatrix}\;\text{ and }\;\psi^h=\begin{pmatrix}h & 0\\ 0& 0\end{pmatrix}.
\]
Here the map $h_i\colon E_{i}\to E_{i-1}$ should be the homotopy from the inclusion of $f$ into $f^{\prime}$. One easily checks that the diagram above is commutative and that this setting gives rise to $\Gamma(\varphi)\simeq_h\Gamma(\varphi\circ f)$. The statement $\Gamma(g\circ \varphi)\simeq_h\Gamma(\varphi)$ can be verified analogously.
\end{proof}
\section{Cubes and projective complexes}\label{sec-techcube}
In this section we define/recall some facts from homological algebra about cubes and projective complexes. We need them in the Sec.~\ref{sec-vkhca}, since we can only ensure that our assignment will be a ``projective chain complex''. That means loosely speaking that faces are ``only commutative up to a sign''.

Let $\mathrm{Cu}_n$ be a standard unit $n$-cube. We can consider this cube as a directed graph by labelling neighbouring vertices $\gamma$ by words $a_{\gamma}$ in $\{0,1\}$ of length $n$ as follows. Choose one vertex and give it the label $0\dots 0$ with $n$-entries. Any of its $n$ neighbours get a different word of length $n$ with exactly one $1$. Continue by changing exactly one entry until every vertex has a label.

For two vertices $\gamma_a,\gamma_b$ that differ by only one entry $k$ one assigns a label for the edge between them by replacing $k$ with a $*$. The edges is oriented from $\gamma_a$ to $\gamma_b$ iff $k=0$ for $\gamma_a$. We denote such an edge by $S\colon\gamma_a\to\gamma_b$. Recall that $R$ denotes a commutative, unital ring of arbitrary characteristic.
\begin{defn}\label{defn-cube}(\textbf{Cube})
An $n$-cube in a $R$-pre-additive category $\mathcal C$ is a mapping
\[
\mathrm{Cu}_n^{\mathcal C}\colon\mathrm{Cu}_n\to\mathcal C
\]
that associates each vertex $\gamma_a$ with an element $\gamma^{\mathcal C}_a\in\Ob(\mathcal C)$ and each edge $S\colon\gamma_a\to\gamma_b$ with an element $S_{a,b}^{\mathcal C}\in\Mor(\gamma^{\mathcal C}_a,\gamma^{\mathcal C}_b)$.

A morphisms of cubes $\phi\colon\mathrm{Cu}_n^{\mathcal C}\to\mathrm{Cu^{\prime}}_n^{\mathcal C}$ is a collection of morphisms for all vertices that is
\[
\{S_{a,{a^{\prime}}}^{\mathcal C}\mid \gamma_a,\gamma_{a^{\prime}}\text{ vertices of }\mathrm{Cu}_n^{\mathcal C},\mathrm{Cu^{\prime}}_n^{\mathcal C}\}.
\]
We denote the category of $n$-cubes in $\mathcal C$ by $\Cb_n(\mathcal C)$ and the category of all cubes in $\mathcal C$ by $\Cb(\mathcal C)$.

It should be noted that a morphisms between two $n$-cubes can be seen as a $n+1$-cube. Moreover, from a $n+1$ cube one can define a morphism of $n$-cubes by fixing a letter $k$ of the words associated to the vertices $\gamma_a$ and fixing the $n$-subcubes $\mathrm{Cu^0}_n^{\mathcal C}$ and $\mathrm{Cu^1}_n^{\mathcal C}$ of $\mathrm{Cu}_{n+1}^{\mathcal C}$ such that the vertices of $\mathrm{Cu^0}_n^{\mathcal C}$ have $k=0$ and the ones of $\mathrm{Cu^1}_n^{\mathcal C}$ have $k=1$. The morphism is the given by all edges of $\mathrm{Cu}_{n+1}^{\mathcal C}$ that change $k=0$ to $k=1$.
\end{defn}
\begin{defn}\label{defn-cubetype}(\textbf{Cube types})
Denote by $R^*$ all units in $R$. The category $\mathcal C_P=\mathcal C/R^*$ is called the \textit{projectivisation} of $\mathcal C$, i.e. morphisms are identified iff they differ only by a unit. A \textit{projectivisation} of a cube $\mathrm{Cu}_n^{\mathcal C_P}$ is given by the composition with the obvious (I used the word again - my bad) projection.

A face of a cube, denoted by $F$, is given by (we hope that the notation is clear)
\[
\xymatrix{
 & \gamma_{a_{01}}\ar[rd]^{S_{*1}} &\\
 \gamma_{a_{00}}\ar[rd]_{S_{*0}}\ar[ru]^{S_{0*}} &  & \gamma_{a_{11}},\\
 & \gamma_{a_{10}}\ar[ru]_{S_{1*}} &}
\]
or with an extra superscript $\mathcal C$ for a cube in $\mathcal C$. Such a face is said to be of \textit{type $a$, $c$ or $p$} if the following is satisfied.
\begin{itemize}
\item[(Type a)] We have $S^{\mathcal C}_{1*}\circ S^{\mathcal C}_{*0}=-S^{\mathcal C}_{*1}\circ S^{\mathcal C}_{0*}$ (\textit{anti-commutative}).
\item[(Type c)] We have $S^{\mathcal C}_{1*}\circ S^{\mathcal C}_{*0}=\phantom{-}S^{\mathcal C}_{*1}\circ S^{\mathcal C}_{0*}$ (\textit{commutative}).
\item[(Type p)] We have $S^{\mathcal C}_{1*}\circ S^{\mathcal C}_{*0}=uS^{\mathcal C}_{*1}\circ S^{\mathcal C}_{0*}$ for $u\in R^*$ (\textit{projective}).
\end{itemize}
Furthermore, a cube $\mathrm{Cu}_n^{\mathcal C}$ is called of \textit{type a}, \textit{type c} or \textit{type p}, if all of its faces are of the corresponding types.

A morphisms between $n$-cubes $\mathrm{Cu}_n^{\mathcal C},\mathrm{Cu^{\prime}}_n^{\mathcal C}$ is called of \textit{type a}, \textit{type c} or \textit{type p} if the corresponding $n+1$-cube is type a, type c or type p respectively.

Two cubes $\mathrm{Cu}_n^{\mathcal C}$ and $\mathrm{Cu^{\prime}}_n^{\mathcal C}$ are called \textit{p-equal} if
\[
\mathrm{Cu}_n^{\mathcal C_P}=\mathrm{Cu^{\prime}}_n^{\mathcal C_P}.
\]
We call two morphisms between $\mathrm{Cu}_n^{\mathcal C}$ and $\mathrm{Cu^{\prime}}_n^{\mathcal C}$ \textit{p-equal} if the corresponding $n+1$-cubes are p-equal.

Note that morphisms of type c and type p are closed under compositions. Hence, the category $\Cb(\mathcal C)$ has three subcategories, namely the following.
\begin{itemize}
\item[(Type a)] The subcategory $\Cb^a(\mathcal C)$ with cubes of type a and morphisms of type c.
\item[(Type c)] The subcategory $\Cb^c(\mathcal C)$ with cubes of type c and morphisms of type c.
\item[(Type p)] The subcategory $\Cb^p(\mathcal C)$ with cubes of type p and morphisms of type p.
\end{itemize}
\end{defn}
It is worth noting that we can see any cube in $\mathcal C$ as a complex $(C_*,c_*)$ (we point out that we do not say \textit{chain} complex here) by taking direct sums of vertices with the same number of $1$ in their labels and matrices of the morphisms associated to the edges between neighbouring vertices.
\begin{defn}\label{defn-chainhomotopy2}
Let $(C_*,c_*)$ and $(D_*,d_*)$ be two cubes of type p. We call two morphisms $\varphi,\varphi^{\prime}\colon C\to D$ of type p \textit{p-homotopic}, denoted by $\varphi\sim^P_h\varphi^{\prime}$, if
\[
\varphi_i-\varphi^{\prime}_i=h_{i+1}\circ c_i+u_id_{i-1}\circ h_i\text{ for all }i\in\bZ,
\]
for a backward diagonal $h$ as in Definition~\ref{defn-chainhomotopy} and units $u_i\in R^*$.

Two such cube complexes are called \textit{p-homotopic} if there exists two morphisms of type p $\varphi\colon C\to D$ and $\psi\colon D\to C$ such that
\[
\psi_i\circ\varphi_i\sim^P_h\mathrm{id}_{C}\;\;\text{ and }\;\;\varphi_i\circ\psi_i\sim^P_h\mathrm{id}_{D}\text{ for all }i\in\bZ.
\]
Such morphisms $\varphi\colon C\to D$ and $\psi\colon C\to D$ are called \textit{p-homotopy equivalences}. We denote p-homotopic complexes of type p by $C\simeq^p_h D$.
\end{defn}
\begin{defn}\label{defn-edge}
Let $\mathrm{Cu}_n$ denote an $n$-cube and let us denote the set of edges of $\mathrm{Cu}_n$ by $\mathrm{E}(\mathrm{Cu}_n)$. An \textit{edge assignment} $\epsilon$ of the cube is a map
\[
\epsilon\colon\mathrm{E}(\mathrm{Cu}_n)\to\{+1,-1\}.
\]
Let $\mathcal C$ be a $R$-pre-additive category. Then an edge assignment $\epsilon$ of the cube   $\mathrm{Cu}_n$ is called \textit{negative} (or \textit{positive}), if $\mathrm{Cu}_n^{\mathcal C}$ is a cube of type a (or of type c) after multiplying the morphism $f_e$ of the edge $e\in\mathrm{E}(\mathrm{Cu}_n)$ of the cube $\Cb(\mathcal C)$ with $\epsilon(e)$.
\end{defn}
The following lemma follows immediately from the definition.
\begin{lem}\label{lem-edge}
If $\mathrm{Cu}_n^{\mathcal C}$ is a cube of type a (or of type c), then there is a negative (or positive) edge assignment of $\mathrm{Cu}_n^{\mathcal C_P}$.
\end{lem}
\begin{proof}
Immediate from the definition.
\end{proof}
\section{Open issues}\label{sec-vkhend}
Here are some open problems that we have observed. Note that nowadays the results about classical Khovanov homology form a highly studied and rich field. So there are many more open questions related to our construction.
\begin{itemize}
\item It is quite remarkable that one has to use a ``$\wedge$-product like'' construction to define even, virtual Khovanov homology. An interpretation of this fact is missing.
\item Our complex is an extension of the classical (even) Khovanov complex. We shortly discuss a method which could lead to an extension of odd Khovanov homology, see~\cite{ors}. Even and odd Khovanov homology differ over $\bQ$ but are equal over $\bZ/2$.
\item Secondly we discuss the relationship between the virtual Khovanov complex and the categorification of the higher quantum polynomials ($n\ge 3$) due to Khovanov in~\cite{kh3} and Mackaay and Vaz in~\cite{mv1} and Mackaay, Sto\v{s}i\'{c} and Vaz in~\cite{msv} using ``cobordism like'' singular surfaces called ``foams''.
\item The results from Sec.~\ref{sec-vkhapp} could lead to an extension of the Rasmussen invariant to virtual knots.
\end{itemize}
On the second point: The reader familiar with the paper of Ozsv\'{a}th, Rasmussen and Szab\'{o} may have already identified our map
\[
F_{Kh}((\Phi^-_+\amalg\mathrm{id}^+_+)\circ\Delta^+_{++})\colon A\to A\otimes A
\]
to be the comultiplication which they use, see Sec. 1.1 in~\cite{ors}.

One main difference between the even and odd Khovanov complex is the usage of this map instead of the standard map $F_{Kh}(\Delta^+_{++})$ and the structure of an exterior algebra instead of direct sums. Furthermore, there are commutative and anti-commutative faces in the odd Khovanov complex. But because every three dimensional cube has an even number of both types of faces, there is a sign assignment which makes every face anti-commute. One major problem is the question how to handle unorientable faces, because these faces can be counted as commutative or anti-commutative. Furthermore, one should admit that faces of type 1a and 1b can be commutative or anti-commutative. Hence, there is still much work to do (life is short, but this paper is not, so...).  

On the third point: The key idea in the categorification of the $\mathfrak{sl}_n$-polynomial for $n\ge 3$ is the usage of so-called foams. This is very interesting approach due to Khovanov, Mackaay, Sto\v{s}i\'{c} and Vaz (see in their papers~\cite{kh3},~\cite{msv} and~\cite{mv1}) that was further generalised recently using categorified $q$-skew Howe duality by Lauda, Queffelec and Rose~\cite{lqr1} for $n=2,3$ and Queffelec and Rose~\cite{qr1} for $n>1$.

So in the virtual case one should use a topological construction with virtual webs and decorated, possibly non-orientable foams (immersed rather than embedded). So their concept to categorify the $\mathfrak{sl}_n$-polynomials for $n=3$ should lift to v-links. This needs further work (the sign assignment seems to be the main point), but seems to be very interesting (from some point of view). The $n>3$ case is indeed more complicated. In their paper Mackaay, Sto\v{s}i\'{c} and Vaz (see~\cite{msv}) use a special formula, the so-called Kapustin-Li formula, to find the adapted relations. But this formula only works in the orientable case and it has no straightforward extension to the non-orientable case. But hopefully the collection of relations they use is already enough to show invariance under the vRM1, vRM2, vRM3 and the mRM moves. At least in the case $n=2$ the local relations are enough to show the invariance.

In fact, Queffelec and Rose are able to avoid the Kapustin-Li formula. This is possible because all the foam relations can be obtained by only using a certain cyclotomic KL-R. It would be interesting to find a (``higher'') representation theoretical ``explanation'' for the v-link homologies.

Note that for all $n$ these constructions are now known to be related to the categorification of the $\mathfrak{sl}_n$ polynomial using matrix factorisations by Khovanov and Rozansky~\cite{kr1}. This follows from work of Lauda, Queffelec and Rose mentioned above and a paper of the author~\cite{tub4} (who - shame on him - used matrix factorisation framework instead of foams) using categorified $q$-skew Howe duality.
\vspace*{0.5cm}

Another maybe interesting point is to generalise Khovanov's arc algebra $H^n$, introduced by Khovanov in~\cite{kh4}, to virtual knots. Note that this algebra can be seen as the algebraic structure behind classical Khovanov homology.

Work by many people, for instance Brundan and Stroppel, see~\cite{bs1} as the beginning of a series of papers~\cite{bs2},~\cite{bs3},~\cite{bs4} and ~\cite{bs5}, have demonstrated that $H^n$ and its generalisations (e.g. the type $A_2$ variant was studied in~\cite{mpt},~\cite{rob2},~\cite{rob1} and~\cite{tub3} and the type $A_n$-web algebra in~\cite{mack1}. There is also a type $D$ version of the arc algebra, see~\cite{ehst1},~\cite{ehst2} and~\cite{ehst3}, and a $\mathfrak{gl}(1|1)$ variant~\cite{sar}) have, in addition to their knot theoretical origin, a beautiful combinatorial and representation theoretical structure.

A virtual generalisation could give a hint what the underlying quantum representation theory of virtual knots and links is (if there is any). At least one could hope to obtain an algebra with an interesting combinatorial structure.

\vspace{0.1in}

\noindent D.T.: { \sl \small Centre for Quantum Geometry of Moduli Spaces (QGM), University Aarhus, Denmark} 
\noindent {\small \textbf{email: dtubben@qgm.au.dk}}
\end{document}